\documentclass[12pt]{elsarticle}

\usepackage[left=1in,top=1in,right=1in,bottom=1in]{geometry}

\usepackage{calc}
\usepackage{fancyhdr}
\usepackage{amsmath}
\usepackage{amssymb}
\usepackage{amsfonts}
\usepackage{graphicx}
\usepackage[linesnumbered,boxed]{algorithm2e}
\usepackage{amsthm}
\usepackage{amscd}
\usepackage{eufrak}
\usepackage{amsmath}
\usepackage{hyperref} 
\usepackage{lscape}
\usepackage{morefloats}
\usepackage{stmaryrd}
\usepackage{color}
\usepackage{placeins}
\usepackage{subfig}
\usepackage{tikz-cd}

\usepackage{schemata}

\begin{document}

\begin{frontmatter}

\title{Divergence-Conforming Isogeometric Collocation Methods for the Incompressible Navier-Stokes Equations}

\author[1]{Ryan M. Aronson\corref{cor1}}
\ead{rmaronso@stanford.edu}

\author[2]{John A. Evans}

\cortext[cor1]{Corresponding author}

\affiliation[1]{organization={Stanford University},
postcode={94305},
city={Stanford, CA},
country={USA}}

\affiliation[2]{organization={University of Colorado Boulder},
postcode={80309},
city={Boulder, CO},
country={USA}}




\begin{abstract}

We develop two isogeometric divergence-conforming collocation schemes for incompressible flow. The first is based on the standard, velocity-pressure formulation of the Navier-Stokes equations, while the second is based on the rotational form and includes the vorticity as an unknown in addition to the velocity and pressure. We describe the process of discretizing each unknown using B-splines that conform to a discrete de Rham complex and collocating each governing equation at the Greville abcissae corresponding to each discrete space. Results on complex domains are obtained by mapping the equations back to a parametric domain using structure-preserving transformations. Numerical results show the promise of the method, including accelerated convergence rates of the three field, vorticity-velocity-pressure scheme when compared to the two field, velocity-pressure scheme.

\end{abstract}

\begin{keyword}
Isogeometric analysis \sep Collocation \sep Incompressible flow \sep Divergence-conforming discretizations \sep Velocity-pressure formulation \sep Vorticity-velocity-pressure formulation
\end{keyword}

\end{frontmatter}

\section{Introduction}

Isogeometric Analysis (IGA) is a technology \cite{hughes_IGA, cottrell_IGAbook} which replaces the standard polynomial basis functions used in traditional Finite Element Analysis (FEA) with B-splines, Non-Uniform Rational B-splines (NURBS), and other classes of splines in an aim to reduce the gap between geometry and analysis. IGA has a distinct advantage over traditional FEA due to its ability to exactly represent the geometries commonly seen in Computer-Aided Design (CAD). Moreover, the basis functions used in IGA are globally more smooth than those of FEA and it has been shown that IGA can exhibit improved accuracy and robustness over FEA. For example, higher continuity splines are shown to have optimal approximating power in the sense of Kolmogorov $n$-widths \cite{evans_nwidths} and spline discretizations have more favorable dissipation and dispersion properties than standard, high-order FEA discretizations \cite{cottrell_IGAbook}.

To improve on the complexity and implementation details of IGA, the feasibility of isogeometric collocation methods has been explored \cite{auricchio_isogeometric_coll,reali_intro}. In Galerkin IGA methods, the discrete system of equations is formed by integrating the PDE residual against the test function space. This requires numerical integration, which renders system assembly quite expensive. Collocation, on the other hand, forms the discrete system by simply requiring that the residual of the governing equations vanish at a set of discrete locations in the domain. 

Many recent studies have shown the efficacy of isogeometric collocation methods in both static and dynamic solid mechanics problems \cite{auricchio_coll_elastostatics, evans_explicit, kruse2015largedef}, and detailed comparisons between isogeometric Galerkin and collocation methods have been made \cite{schillinger_cost_comp}. In addition, recent studies have also investigated the use of mixed collocation methods for use in nearly incompressible elasticity \cite{fahrendorf_mixed, morganti_mixed}. Isogeometric collocation has even been used for acoustic problems \cite{zampieri2021wave}, computing Karhunen-Loeve expansions \cite{jahanbin2019isogeometric, mika2022comparison}, and introduced into physics-informed neural networks \cite{moller2021pinn}. Finally, a method was recently introduced for IGA collocation in immersed domains by combining with the finite cell method near the boundaries \cite{torre2023immersed}. These results all suggest that isogeometric collocation methods retain some of the improved qualities of standard IGA, while reducing the computational cost and improving sparsity of the discrete systems. 

Isogeometric collocation methods have not been as well explored in the context of fluid mechanics, though the idea of using B-spline collocation to solve incompressible fluid mechanics problems has been investigated in the past \cite{botella_collocation,johnson_Greville_coll}. In addition, spline collocation has been employed in fundamental Direct Numerical Simulation (DNS) studies of turbulent flows \cite{lee2015DNS}. However, the methods previously introduced are typically limited to simple geometries and are not divergence-conforming, meaning that the discrete velocity field does not exactly satisfy the continuity equation in incompressible flow. In addition, a mixed isogeometric collocation method has recently been proposed for use in poromechanics \cite{morganti2018poromech}, though the preliminary results were limited to one-dimensional problems. 

In the context of incompressible fluid mechanics, divergence-conforming Galerkin methods based on B-spline basis functions have been developed for both the Stokes and the Navier-Stokes equations \cite{buffa_stokes, evans_steady_NS, evans_unsteady_NS}. These methods are provably inf-sup stable and yield discrete velocity fields that are exactly pointwise divergence free, among other desirable qualities such as pressure robustness. An excellent summary of divergence-conforming methods is given by \cite{john_divconstraint}. These discretizations have prospered in areas such as turbulent flow simulation \cite{vanOpstal_divconformingVMS, evans_book, evans2020vms} and fluid-structure interaction \cite{kamensky_immersogeometric}. Moreover, efficient multigrid solvers have been developed based on these discretizations \cite{coley_multigrid}. Divergence-conforming Galerkin methods have also been developed for more advanced spline discretizations, such as hierarchical B-splines \cite{evans2020hierarchical} and LR B-splines \cite{johannessen2015divergence}.

In this paper we develop similar divergence-conforming methodologies for incompressible flow using collocation. In particular, we introduce two collocation methods, one based on the standard velocity-pressure form of the steady Navier-Stokes equations, and one based on a three field (velocity-vorticity-pressure) form of the steady Navier-Stokes equations. The latter form of the Navier-Stokes equations has recently been used to develop alternative structure-preserving finite element discretizations \cite{zhang2021mass, hanot2021arbitrary, gopalakrishnan2021minimal} and we find that collocation methods based on the resulting system of first order differential-algebraic equations returns improved convergence rates compared to the rates obtained using collocation in conjunction with the standard velocity-pressure form of the equations. In our collocation schemes, each unknown is discretized with compatible B-spline spaces that preserve the structure of the governing equations. Both collocation methods in this paper are shown to return velocity fields which are still exactly pointwise divergence free, similar to the Galerkin methods mentioned above. 

We lay out this paper as follows: In Section 2 we describe the steady form of the Navier-Stokes equations using velocity and pressure unknowns as well as vorticity, velocity, and pressure unknowns. This is followed in Section 3 by a discussion of the de Rham complex and isogeometric discrete differential forms, the tools used to develop a divergence-conforming method. Section 4 describes the collocation schemes for square domains in two dimensions. Then results are presented in the two dimensional setting in Section 5 which detail the high-order convergence rates of the methods as well as agreement with standard benchmark problems. We then discuss the necessary changes to make the methods work for cubic domains in three dimensions and illustrate that the methods performs similarly in this setting in Sections 6 and 7. Finally, we return to 2D in Section 8 and consider the Stokes equations in more complicated domains. We show that by mapping the equations and unknowns via divergence and integral preserving transformations we can also obtain results for flow in complex geometries. Section 9 summarizes these results. 

\section{Velocity-Pressure and Vorticity-Velocity-Pressure Forms of the Navier-Stokes Equations}

In this paper we consider the steady, incompressible Navier-Stokes equations on a Lipschitz open set $\Omega$ of points in either $\mathbb{R}^2$ or $\mathbb{R}^3$ when subjected to Dirichlet boundary conditions. The standard form of this problem with $d = 2,3$ is stated as follows:

\bigskip

$$
\left\{ \hspace{5pt}
\parbox{5in}{
\noindent Given $\nu \in \mathbb{R}^+$, $\textbf{f} : \Omega \rightarrow \mathbb{R}^d$, and $\textbf{g} : \partial \Omega \rightarrow \mathbb{R}^d$, find $\mathbf{u} : \Omega \rightarrow \mathbb{R}^d$ and $p : \Omega \rightarrow \mathbb{R}$ such that:

\begin{align}
    -\nu \Delta \mathbf{u} + \mathbf{u} \cdot \nabla \mathbf{u} + \nabla p &= \mathbf{f} \quad \textup{in} \quad \Omega \label{eq:vp_mom} \\
    \nabla \cdot \mathbf{u} &= 0 \quad \textup{in} \quad \Omega \label{eq:vp_cont}\\
    \mathbf{u} &= \mathbf{g} \quad \textup{on} \quad \partial \Omega. \label{eq:vp_BC}
\end{align}
}
\right.
$$

\bigskip

\noindent Equations \eqref{eq:vp_mom} and \eqref{eq:vp_cont} represent the standard forms of the momentum and mass conservation equations, while Equation \eqref{eq:vp_BC} sets the Dirichlet boundary values. Above, we denote the velocity field by $\mathbf{u}$, the kinematic pressure field by $p$, the constant kinematic viscosity by $\nu$, the applied forcing as $\mathbf{f}$, and the prescribed Dirichlet boundary values as $\mathbf{g}$.

For the purposes of this paper we will not only work with this set of equations, but also introduce vorticity $\boldsymbol{\omega}$ as a separate unknown variable and introduce $\boldsymbol{\omega} - \nabla \times \textbf{u} = 0$ as a constitutive relation. Substituting the two vector calculus identities:

\begin{equation}
    \Delta \mathbf{u} = \nabla (\nabla \cdot \mathbf{u}) - \nabla \times (\nabla \times \mathbf{u}) = -\nabla \times \boldsymbol{\omega},
\end{equation}
\begin{equation}
    \mathbf{u} \cdot \nabla \mathbf{u} = (\nabla \times \mathbf{u}) \times \mathbf{u} + \frac{1}{2} \nabla (\mathbf{u} \cdot \mathbf{u}) = \boldsymbol{\omega} \times \mathbf{u} + \frac{1}{2} \nabla (\mathbf{u} \cdot \mathbf{u}),
\end{equation}

\noindent into Equation \eqref{eq:vp_mom} when $d=3$, we arrive at the vorticity-velocity-pressure formulation of the problem in 3D:

\bigskip

$$
\left\{ \hspace{5pt}
\parbox{5in}{
\noindent Given $\nu \in \mathbb{R}^+$, $\textbf{f} : \Omega \rightarrow \mathbb{R}^3$, and $\textbf{g} : \partial \Omega \rightarrow \mathbb{R}^3$, find $\textbf{u} : \Omega \rightarrow \mathbb{R}^3$, $P : \Omega \rightarrow \mathbb{R}$, and $\boldsymbol{\omega} : \Omega \rightarrow \mathbb{R}^3$ such that:

\begin{align}
    \nu \nabla \times \boldsymbol{\omega} + \boldsymbol{\omega} \times \textbf{u} + \nabla P &= \textbf{f} \quad \textup{in} \quad \Omega \label{eq:vvp_mom_3D} \\
    \nabla \cdot \textbf{u} &= 0 \quad \textup{in} \quad \Omega \label{eq:vvp_cont_3D} \\
    \boldsymbol{\omega} - \nabla \times \textbf{u} &= 0 \quad \textup{in} \quad \Omega \label{eq:vvp_const_3D} \\
    \textbf{u} &= \textbf{g} \quad \textup{on} \quad \partial \Omega. \label{eq:vvp_BC_3D}
\end{align}
}
\right.
$$

\bigskip

\noindent Note that in the above formulation we have replaced the kinematic pressure $p$ with the total pressure $P$, which are related via $P = p + \frac{1}{2}\textbf{u}\cdot\textbf{u}$. 

For later sections in this paper it is useful to employ the component forms of the vector equations above describing the vorticity-velocity-pressure formulation. When explicitly broken into its components, the momentum conservation equation, given by Equation \eqref{eq:vvp_mom_3D}, becomes:

\begin{equation}
    \nu (\frac{\partial \omega_z}{\partial y} - \frac{\partial \omega_y}{\partial z}) + (\omega_y u_z - \omega_z u_y) + \frac{\partial P}{\partial x} = f_x,
\end{equation}
\begin{equation}
    \nu (\frac{\partial \omega_x}{\partial z} - \frac{\partial \omega_z}{\partial x}) + (\omega_z u_x - \omega_x u_z) + \frac{\partial P}{\partial y} = f_y,
\end{equation}
\begin{equation}
    \nu (\frac{\partial \omega_y}{\partial x} - \frac{\partial \omega_x}{\partial y}) + (\omega_x u_y - \omega_y u_x) + \frac{\partial P}{\partial z} = f_z,
\end{equation}

\noindent and the constitutive relation given by Equation \eqref{eq:vvp_const_3D} reads:

\begin{equation}
    \omega_x - (\frac{\partial u_z}{\partial y} - \frac{\partial u_y}{\partial z}) = 0,
\end{equation}
\begin{equation}
    \omega_y - (\frac{\partial u_x}{\partial z} - \frac{\partial u_z}{\partial x}) = 0,
\end{equation}
\begin{equation}
    \omega_z - (\frac{\partial u_y}{\partial x} - \frac{\partial u_x}{\partial y}) = 0.
\end{equation}

The above component form of the equations is also useful for considering 2D problems in the three field formulation, as we do not need to redefine operations such as cross products when the vorticity reduces to a scalar unknown. Thus we can arrive at the problem statement for 2D domains by simply removing any terms involving $\omega_x$, $\omega_y$, or any derivatives in the $z$ direction. In full, the 2D problem reads:

\bigskip

$$
\left\{ \hspace{5pt}
\parbox{5in}{
\noindent Given $\nu \in \mathbb{R}^+$, $\textbf{f} : \Omega \rightarrow \mathbb{R}^2$, and $\textbf{g} : \partial \Omega \rightarrow \mathbb{R}^2$, find $\textbf{u} : \Omega \rightarrow \mathbb{R}^2$, $P : \Omega \rightarrow \mathbb{R}$, and  $\omega : \Omega \rightarrow \mathbb{R}$ such that:

\begin{align}
    \nu \frac{\partial \omega}{\partial y} - \omega u_y + \frac{\partial P}{\partial x} = f_x\quad &\textup{in} \quad \Omega \label{eq:vvp_momx_2D} \\
    - \nu \frac{\partial \omega}{\partial x} + \omega u_x + \frac{\partial P}{\partial y} = f_y\quad &\textup{in} \quad \Omega \label{eq:vvp_momy_2D}\\
    \nabla \cdot \textbf{u} = 0 \quad &\textup{in} \quad \Omega \label{eq:vvp_cont_2D}\\
  \omega - (\frac{\partial u_y}{\partial x} - \frac{\partial u_x}{\partial y}) = 0 \quad &\textup{in} \quad \Omega \label{eq:vvp_const_2D} \\
  \textbf{u} = \textbf{g} \quad &\textup{on} \quad \partial \Omega. \label{eq:vvp_BC_2D}
\end{align}
}
\right.
$$

\bigskip

With the governing equations fully defined we can move to a more in-depth description of the collocation scheme, starting with the definition of discrete approximation spaces in the following section.

\section{The de Rham Complex and Isogeometric Discrete Differential Forms}

The first step in creating a collocation scheme is to define the sets of basis functions used to approximate the unknown variables. To construct the collocation methods presented in this paper, we leverage the de Rham complex which aids in the development of exactly divergence-conforming finite element spaces. After recalling the de Rham complex, we describe the process of constructing B-spline basis functions and conclude with the definition of B-spline spaces that conform to the de Rham complex.

\subsection{The de Rham Complex}

The de Rham complex is a cochain complex that is often used as a starting point for developing mixed finite element methods which preserve topological properties of the continuous problem and are typically more stable in practice \cite{john_divconstraint}. In 3D, it is typically written as:

\begin{equation}
\begin{tikzcd}
\mathbb{R}\arrow[r] & \Phi \arrow[r,"\nabla"] &  \boldsymbol{\Psi} \arrow[r,"\nabla \times"]&  \boldsymbol{\mathcal{V}} \arrow[r,"\nabla \cdot"] & \mathcal{Q} \arrow[r]& 0, 
\end{tikzcd}
\label{eq:deRham}
\end{equation}

\noindent where:

\begin{align}
    \Phi &:= \mathit{H}^1(\Omega), \\
    \boldsymbol{\Psi} &:= \mathbf{H}(\mathbf{curl},\Omega), \\
    \boldsymbol{\mathcal{V}} &:= \mathbf{H}(\mathbf{div},\Omega), \\
    \mathcal{Q} &:= \mathit{L}^2(\Omega).
\end{align}

In the context of fluid flow, these infinite dimensional spaces can be interpreted as the spaces of scalar potential fields ($\Phi$), vector potential fields ($\boldsymbol{\Psi}$), velocity fields ($\boldsymbol{\mathcal{V}}$), and pressure fields ($\mathcal{Q}$). This complex is exact for simply connected domains, meaning that the range of each map is the same as the null space of the following map.

For completeness, we also state the rotated 2D de Rham complex:

\begin{equation}
\begin{tikzcd}
\mathbb{R} \arrow[r] & {\Psi} \arrow[r,"\nabla^\perp"]&  \boldsymbol{\mathcal{V}} \arrow[r,"\nabla \cdot"] & \mathcal{Q} \arrow[r]& 0,
\end{tikzcd}
\label{eq:deRham2D}
\end{equation}

\noindent where:

\begin{align}
    {\Psi} &:= {H}^1(\Omega), \\
    \boldsymbol{\mathcal{V}} &:= \mathbf{H}(\mathbf{div},\Omega), \\
    \mathcal{Q} &:= \mathit{L}^2(\Omega),
\end{align}

\noindent and the rotor operator $\nabla^\perp$ acting on a scalar function $\omega$ is defined as $\nabla^\perp \omega = (\frac{\partial \omega}{\partial y}, -\frac{\partial \omega}{\partial x})$.

Not only is the topological structure of the incompressible Navier-Stokes equations embedded in the de Rham complex, but stability conditions result from the complex as well. By creating approximation spaces for the unknown variables that conform to discrete analogs of Equations \eqref{eq:deRham} and \eqref{eq:deRham2D} we can generate numerical methods which inherit these properties. More concretely, for 3D problems let the space $\Phi_h$ contain the discrete scalar potentials, $\boldsymbol{\Psi}_h$ contain the discrete vector potentials as well as the discrete vorticity, $\boldsymbol{\mathcal{V}}_h$ contain the discrete velocity, and $Q_h$ contain the discrete pressure. Then if there exist projection operators $\Pi_\Phi : \Phi \rightarrow \Phi_h $, $\Pi_{\boldsymbol{\Psi}} : \boldsymbol{\Psi} \rightarrow \boldsymbol{\Psi}_h $, $\Pi_{\boldsymbol{\mathcal{V}}} : \boldsymbol{\mathcal{V}} \rightarrow \boldsymbol{\mathcal{V}}_h $, and $\Pi_{\mathcal{Q}} : \mathcal{Q} \rightarrow \mathcal{Q}_h $ such that the following commuting diagram holds


\begin{equation}
    \begin{CD}
        \mathbb{R} @>>> \Phi @>\nabla>> \boldsymbol{\Psi} @>\nabla \times>> \boldsymbol{\mathcal{V}} @>\nabla \cdot>> \mathcal{Q} @>>> 0\\
        @. @VV\Pi_{\Phi}V @VV\Pi_{\boldsymbol{\Psi}}V @VV\Pi_{\boldsymbol{\mathcal{V}}}V @VV\Pi_{\mathcal{Q}}V @.\\
        \mathbb{R} @>>> \Phi_h @>\nabla>> \boldsymbol{\Psi}_h @>\nabla \times>> \boldsymbol{\mathcal{V}}_h @>\nabla \cdot>> \mathcal{Q}_h @>>> 0, 
    \end{CD}
    \label{eq:deRham_discrete}
\end{equation}

\noindent a Galerkin finite element method employing $\boldsymbol{\mathcal{V}}_h$ for the discrete velocity space and $\mathcal{Q}_h$ for the discrete pressure space will be inf-sup stable and will yield discrete velocity approximations that are divergence free almost everywhere \cite{evans2013darcy}. We shall prove later on that the divergence-conforming property is maintained if we utilize these spaces in our collocation scheme. 

The same holds in 2D, where we instead let $\Psi_h$ define the discrete space to which the vorticity belongs (as well as the streamfunction), $\boldsymbol{\mathcal{V}}_h$ define the discrete velocity space, and $\mathcal{Q}_h$ define the discrete pressure space. The required commuting diagram in this case is 


\begin{equation}
    \begin{CD}
        \mathbb{R} @>>> \Psi @>\nabla^\perp>> \boldsymbol{\mathcal{V}} @>\nabla \cdot>> \mathcal{Q} @>>> 0\\
        @. @VV\Pi_{{\Psi}}V @VV\Pi_{\boldsymbol{\mathcal{V}}}V @VV\Pi_{\mathcal{Q}}V @.\\
        \mathbb{R} @>>> \Psi_h @>\nabla^\perp>> \boldsymbol{\mathcal{V}}_h  @>\nabla \cdot>> \mathcal{Q}_h @>>> 0. 
    \end{CD}
    \label{eq:deRham2D_discrete}
\end{equation}

Of course, we have yet to define the specifics of how to construct discrete spaces such that these discrete complexes hold. For the purposes of this paper we will use compatible B-spline spaces, and the following section is devoted to introducing the basics of B-spline basis functions.

\subsection{Univariate and Multivariate B-Splines}

The construction of a B-spline basis in one dimension requires two objects: the degree of the basis (denoted $k$) and a series of numbers called the knot vector (denoted $\Xi = \{\xi_1, ... \xi_{n+k+1}\}$). The knots $\xi_i$ are non-decreasing and denote the locations in parametric space where the parametrization can change, similar to element boundaries in standard FEA. The number $n$ in the previous relation represents the total number of functions in the basis. The basis functions themselves are defined through the Cox-de Boor recursion: The $k = 0$ basis functions are built as

\begin{equation}
    N_{i,0}(\xi) = \begin{cases} 
                1 & \xi_i \leq \xi \leq \xi_{i+1} \\
                0 & \text{otherwise},
              \end{cases}
\end{equation}

\noindent and higher-order bases are defined through

\begin{equation}
    N_{i,k}(\xi) = \frac{\xi - \xi_i}{\xi_{i+k} - \xi_i}N_{i,k-1}(\xi) + \frac{\xi_{i+k+1} - \xi}{\xi_{i+k+1} - \xi_{i+1}}N_{i+1,k-1}(\xi).
\end{equation}

\noindent Note that in the above relations, we must utilize the convention that any occurrence of zero divided by zero is equal to zero. 

In higher dimensions (two or three for the purposes of this paper), B-spline basis functions are constructed by simply taking the tensor product of one dimensional B-spline bases in each parametric direction. Note that different polynomial degrees and knot vectors can be used in each direction. 

B-spline basis functions possess a number of useful properties for numerical method development. In particular, the smoothness of the global basis at the knot locations is controlled by the repetition of the knot value in $\Xi$. The basis at these locations is $C^{k-r}$, where $r$ is the multiplicity of the knot. Compared to a standard finite element basis, this basis has improved global continuity, which enables the use of collocation as more classical derivatives of the functions are well defined. Note that if the first and last entries in the knot vector are repeated $k+1$ times, the spline basis will become interpolatory at those locations. Such knot vectors are referred to as open knot vectors, and allow for easy specification of Dirichlet boundary conditions. For the results within this paper, all other entries in the knot vector have multiplicity one, meaning we are utilizing a B-spline basis with the highest possible order of continuity. Further, let us define the so-called regularity vector $\boldsymbol{\alpha}$ for a basis. The size of this vector is equal to the number of distinct knots, with entries equal to the polynomial degree of the basis minus the multiplicity of the corresponding knot. In terms of the regularity vector, the global basis functions are $C^{\alpha_j}$-continuous across the $\alpha_j$ unique knot. To simplify notation, the space of functions spanned by a 1D B-spline basis of degree $k$ and a provided knot vector is denoted as:

\begin{equation}
    S^k_{\boldsymbol{\alpha}} = \text{span}\{ N_{i,k} \}_{i = 1}^{n}.
\end{equation}

\noindent We extend this notation to higher dimensions by adding extra sub- and superscripts, representing the polynomial degrees and regularities in each spatial direction.

\subsection{Isogeometric Discrete Differential Forms}

Using the basics of B-spine functions above allows us to develop discrete approximation spaces for the vorticity, velocity, and pressure. These results are built upon work in the area of isogeometric discrete differential forms \cite{buffa_diff_forms,buffa_electromag}, which we will not fully develop here. The construction of these types of spaces in the context of Galerkin approximation for the Navier-Stokes equations can also be found in \cite{evans_steady_NS,evans_unsteady_NS}. For 3D problems, the B-spline spaces used to discretize our unknown fields are given by:

\begin{align}
    \Phi_h &:= \{\phi_h \in S^{k_1,k_2,k_3}_{\boldsymbol{\alpha_1},\boldsymbol{\alpha_2},\boldsymbol{\alpha_3}} \}, \label{eq:3d_potential_space} \\
    \boldsymbol{\Psi}_h &:= \{\boldsymbol{\psi}_h \in S^{k_1-1,k_2,k_3}_{\boldsymbol{\alpha_1}-1,\boldsymbol{\alpha_2},\boldsymbol{\alpha_3}} \times  S^{k_1,k_2-1,k_3}_{\boldsymbol{\alpha_1},\boldsymbol{\alpha_2}-1,\boldsymbol{\alpha_3}} \times S^{k_1,k_2,k_3-1}_{\boldsymbol{\alpha_1},\boldsymbol{\alpha_2},\boldsymbol{\alpha_3}-1}  \}, \\
    \boldsymbol{\mathcal{V}}_h &:= \{\textbf{w}_h \in S^{k_1,k_2-1,k_3-1}_{\boldsymbol{\alpha_1},\boldsymbol{\alpha_2}-1,\boldsymbol{\alpha_3}-1} \times  S^{k_1-1,k_2,k_3-1}_{\boldsymbol{\alpha_1}-1,\boldsymbol{\alpha_2},\boldsymbol{\alpha_3}-1} \times S^{k_1-1,k_2-1,k_3}_{\boldsymbol{\alpha_1}-1,\boldsymbol{\alpha_2}-1,\boldsymbol{\alpha_3}}\} \label{eq:3d_velocity_space},  \\
    \mathcal{Q}_h &:= \{q_h \in S^{k_1-1,k_2-1,k_3-1}_{\boldsymbol{\alpha_1}-1,\boldsymbol{\alpha_2}-1,\boldsymbol{\alpha_3}-1}\}. \label{eq:3d_pressure_space} 
\end{align}

\noindent It can be shown that these spaces satisfy the discrete complex in 
 Equation \eqref{eq:deRham_discrete}.

In practice we usually define $k_1 = k_2 = k_3$, and thus we can define the polynomial degree of the spline bases constructed in the above manner using a single number $k' = k_1-1 = k_2-1 = k_3-1$. This indicates that the pressure space $Q_h$ will have degree equal to $k'$ in each direction. Then, according to the above, each velocity component will have degree $k'+1$ in one direction and degree $k'$ in the other two. Similarly, the vorticity components will have degree $k'+1$ in two directions and degree $k$ in the last. 

In 2D, we define the following spline spaces:

\begin{align}
    \Psi_h &:= \{ \psi_h \in S^{k_1,k_2}_{\boldsymbol{\alpha_1},\boldsymbol{\alpha_2}} \}, \label{eq:2d_discrete_vort_space} \\
    \boldsymbol{\mathcal{V}}_h &:= \{\textbf{w}_h \in S^{k_1,k_2-1}_{\boldsymbol{\alpha_1},\boldsymbol{\alpha_2}-1} \times  S^{k_1-1,k_2}_{\boldsymbol{\alpha_1}-1,\boldsymbol{\alpha_2}} \} , \label{eq:2d_discrete_vel_space} \\
    \mathcal{Q}_h &:= \{q_h \in S^{k_1-1,k_2-1}_{\boldsymbol{\alpha_1}-1,\boldsymbol{\alpha_2}-1} \}. \label{eq:2d_discrete_pressure_space}
\end{align}

\noindent Similar to the 3D setting, these spaces are related as in Equation \eqref{eq:deRham2D_discrete}.

\section{Collocation Methods on Square Domains}

Using the discrete spaces developed above, this section focuses on the construction of collocation methods for the Navier-Stokes equations using divergence-conforming bases. Here we develop methods based on the velocity-pressure form of the Navier-Stokes equations as well as the vorticity-velocity-pressure form. As the vorticity changes between a scalar in the 2D case and a vector in the 3D case, we start by considering only square domains in 2D. This selection also lends itself to easier visualization of the methods. After briefly reviewing the form of a typical divergence-conforming isogeometric Galerkin method, we define the collocation grids for each unknown and then describe the imposition of boundary conditions. The section concludes by summarizing the form of the discrete system. 

\subsection{Review of Galerkin Methods}

We start by reviewing the form of the divergence-conforming isogeometric Galerkin methods which inspired our collocation schemes. Let us consider a problem with Dirichlet boundary conditions on the velocity for concreteness. Then we define the discrete test and trial function spaces for velocity as $\boldsymbol{\mathcal{V}}_{h,0}$ and $\boldsymbol{\mathcal{V}}_{h,\mathbf{g}}$, which are defined as the same $\boldsymbol{\mathcal{V}}_h$ from Equation \eqref{eq:2d_discrete_vel_space} with either no penetration boundary conditions strongly enforced (for the test space $\boldsymbol{\mathcal{V}}_{h,0}$) or with the normal velocity prescribed as given by the boundary data $\mathbf{g}$ at specified collocation points (for the trial space $\boldsymbol{\mathcal{V}}_{h,\mathbf{g}}$). Similarly, define the test and trial space for pressure as $\mathcal{Q}_{h,0}$, where $\mathcal{Q}_h$ is the same space as in Equation \eqref{eq:2d_discrete_pressure_space} but with the added condition that the pressure must have zero integral. Then the Galerkin formulation for the velocity-pressure form would read

\bigskip

$$
\left\{ \hspace{5pt}
\parbox{6in}{
\noindent Given $\nu \in \mathbb{R}^+$, $\textbf{f} : \Omega \rightarrow \mathbb{R}^2$, and $\textbf{g} : \partial \Omega \rightarrow \mathbb{R}^2$, find $\mathbf{u}^h \in \boldsymbol{\mathcal{V}}_{h,\mathbf{g}}$ and $p^h \in \mathcal{Q}_{h,0}$ such that, $\forall (\mathbf{w}^h, q^h) \in (\boldsymbol{\mathcal{V}}_{h,0}, \mathcal{Q}_{h,0})$:

\begin{equation}
\begin{split}
    &\int_{\Omega} (\nu \nabla \mathbf{w}^h \cdot \nabla \mathbf{u}^h + \mathbf{w}^h \cdot (\mathbf{u}^h \cdot \nabla \mathbf{u}^h) - p^h \nabla \cdot \mathbf{w}^h) d\Omega 
    \\ & - \nu \int_{\partial \Omega} ((\nabla \mathbf{u}^h \cdot \mathbf{n}) \cdot \mathbf{w}^h - \frac{C_{pen}}{h} \mathbf{u}^h \cdot \mathbf{w}^h )d\Gamma = \int_{\Omega} \mathbf{w}^h \cdot \mathbf{f} d \Omega + \nu \int_{\partial \Omega} \frac{C_{pen}}{h} \mathbf{g} \cdot \mathbf{w}^h d\Gamma \label{eq:2d_weak_nit}
\end{split}
\end{equation}
\begin{equation}
    \int_{\Omega} q^h (\nabla \cdot \mathbf{u}^h) d \Omega = 0.
\end{equation}
}
\right.
$$

\bigskip

Note that in the above we have used a Nitsche approach to enforce the tangential boundary conditions in the momentum equations. This Galerkin formulation is valid if and only if the minimum entry in the regularity vectors satisfy $\min \{ \boldsymbol{\alpha_1}-1 \} \geq 0$ and $\min \{ \boldsymbol{\alpha_2}-1 \} \geq 0$, where $\boldsymbol{\alpha_1}$ and $\boldsymbol{\alpha_2}$ are the regularity vectors from Equations \eqref{eq:2d_discrete_vort_space} - \eqref{eq:2d_discrete_pressure_space}. We can write a similar Galerkin form of the vorticity-velocity-pressure form of the Navier-Stokes equations, which would yield

\bigskip

$$
\left\{ \hspace{5pt}
\parbox{6in}{
Given $\nu \in \mathbb{R}^+$, $\textbf{f} : \Omega \rightarrow \mathbb{R}^2$, and $\textbf{g} : \partial \Omega \rightarrow \mathbb{R}^2$, find $\mathbf{u}^h \in \boldsymbol{\mathcal{V}}_{h,\mathbf{g}}$,  $P^h \in \mathcal{Q}_{h,0}$, and $\omega^h \in \Psi_h$ such that, $\forall (\mathbf{w}^h, q^h, \psi^h) \in (\boldsymbol{\mathcal{V}}_{h,0}, \mathcal{Q}_{h,0}, \Psi_h)$:

\begin{equation}
    \int_\Omega \nu \frac{\partial \omega^h}{\partial y} w_x^h d \Omega - \int_\Omega \omega^h u_y^h w_x^h d \Omega - \int_\Omega \frac{\partial w_x^h}{\partial x} P^h d \Omega + \int_\Gamma P^h w_x^h dy = \int_\omega f_x w_x^h d \Omega 
\end{equation}
\begin{equation}
    -\int_\Omega \nu \frac{\partial \omega^h}{\partial x} w_y^h d \Omega + \int_\Omega \omega^h u_x^h w_y^h d \Omega - \int_\Omega \frac{\partial w_y^h}{\partial y} P^h d \Omega + \int_\Gamma P^h w_y^h dx = \int_\omega f_y w_y^h d \Omega 
\end{equation}
\begin{equation}
    \int_\Omega (\nabla \cdot \mathbf{u}^h)q^h d \Omega = 0
\end{equation}
\begin{equation}
    \int_\Omega \omega^h \psi^h d \Omega + \int_\Omega (\frac{\partial \psi^h}{\partial x} u_y^h - \frac{\partial \psi^h}{\partial y} u_x^h) d \Omega  = \int_{\partial \Omega} \psi^h(\mathbf{g} \cdot \mathbf{s}) d\Gamma.
    \label{eq:weak_const_2d}
\end{equation}
}
\right.
$$

\bigskip

\noindent In this case the tangential velocity boundary conditions appear as natural boundary conditions in the weak form of the constitutive equation, with $\mathbf{s}$ representing the unit tangent on the boundary (oriented counter-clockwise). In contrast with the velocity-pressure Galerkin formulation, this three field Galerkin formulation is valid if and only if the minimum regularity of the discrete spaces satisfies $\min \{ \boldsymbol{\alpha_1}-1 \} \geq -1$ and $\min \{ \boldsymbol{\alpha_2}-1 \} \geq -1$ due to the reduced differential order of the strong form. 

However, we wish to pursue a collocation method inspired by the divergence-free mixed finite element form, which we describe below. Collocation imposes additional regularity requirements on the spaces of unknowns, as the unknown fields and their derivatives will be evaluated at points rather than integrated over the domain. The spaces developed above are not only divergence-conforming, meaning discrete velocity approximations will be pointwise divergence-free, but are also regular enough to use in collocation provided the polynomial degree is sufficiently large and the global basis is sufficiently regular.

\subsection{Collocation Grids}

Similarly to the Galerkin setting, each discrete unknown in the collocation schemes is assumed to lie in the corresponding space from above: The discrete velocity $\mathbf{u}^h \in \boldsymbol{\mathcal{V}}_{h,\mathbf{g}}$, the discrete pressure $p^h \in \mathcal{Q}_{h,0}$, and, when applicable, the discrete vorticity $\omega^h \in \Psi_h$. For now we ignore any boundary conditions; these will be discussed in the following section. To generate the system of equations needed to solve for the coefficients for each basis function, we define sets of collocation points of the form $\boldsymbol{\tau}_i$ for $i = 1, ..., N$ for each of the governing equations. Note that the total number of collocation points should be equal to the total number of degrees of freedom in the discretization. The full discrete system is formed by requiring the strong form of the governing equations to hold at each of the collocation points.

We choose the Greville abscissae of a B-spline space as collocation points. In one dimension, the Greville abscissae are defined by 

\begin{equation}
    \hat{\xi_i} = \frac{\xi_{i+1} + ... + \xi_{i+p}}{p},
\end{equation}

\noindent and in higher dimensions we simply take the tensor product of the Greville abscissae in each direction. By construction there will be the same number of Greville points as there are basis functions in the considered space. There are choices for collocation points other than the Greville abscissae, such as the Cauchy-Galerkin points \cite{gomez_cauchygalerkin,montardini_optimal} or the Demko abscissae \cite{demko_abs}. However, the Greville points are very easy to compute and have already been demonstrated to give satisfactory results in practical applications (see for example \cite{auricchio_coll_elastostatics,johnson_Greville_coll}). 

\begin{figure}
\centering
\subfloat[Before strong enforcement of normal velocity boundary conditions]{\label{sfig:vp_grid_a}\includegraphics[width=.4\textwidth]{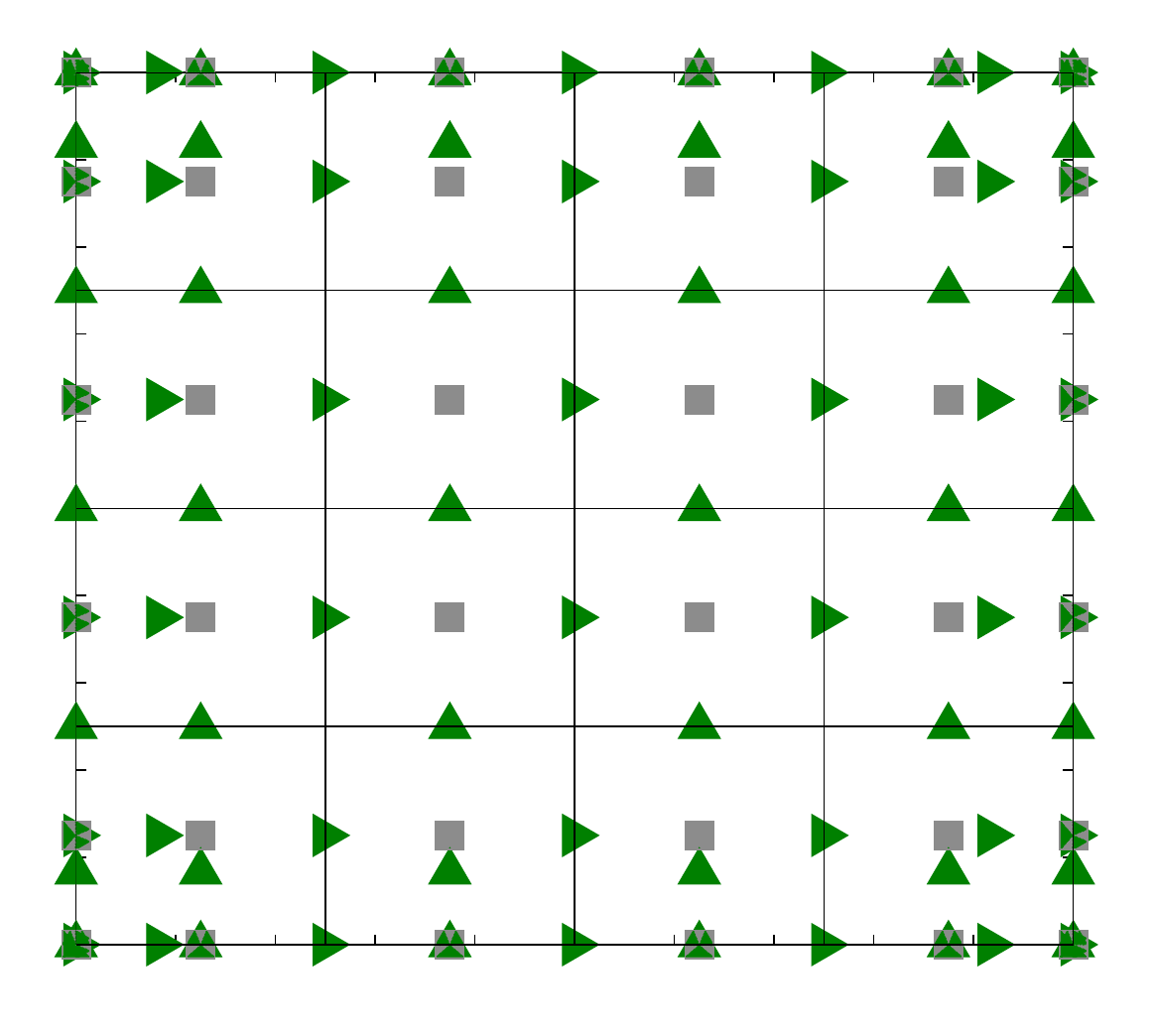}}\hfill
\subfloat[After strong enforcement of normal velocity boundary conditions]{\label{sfig:vp_grid_b}\includegraphics[width=.4\textwidth]{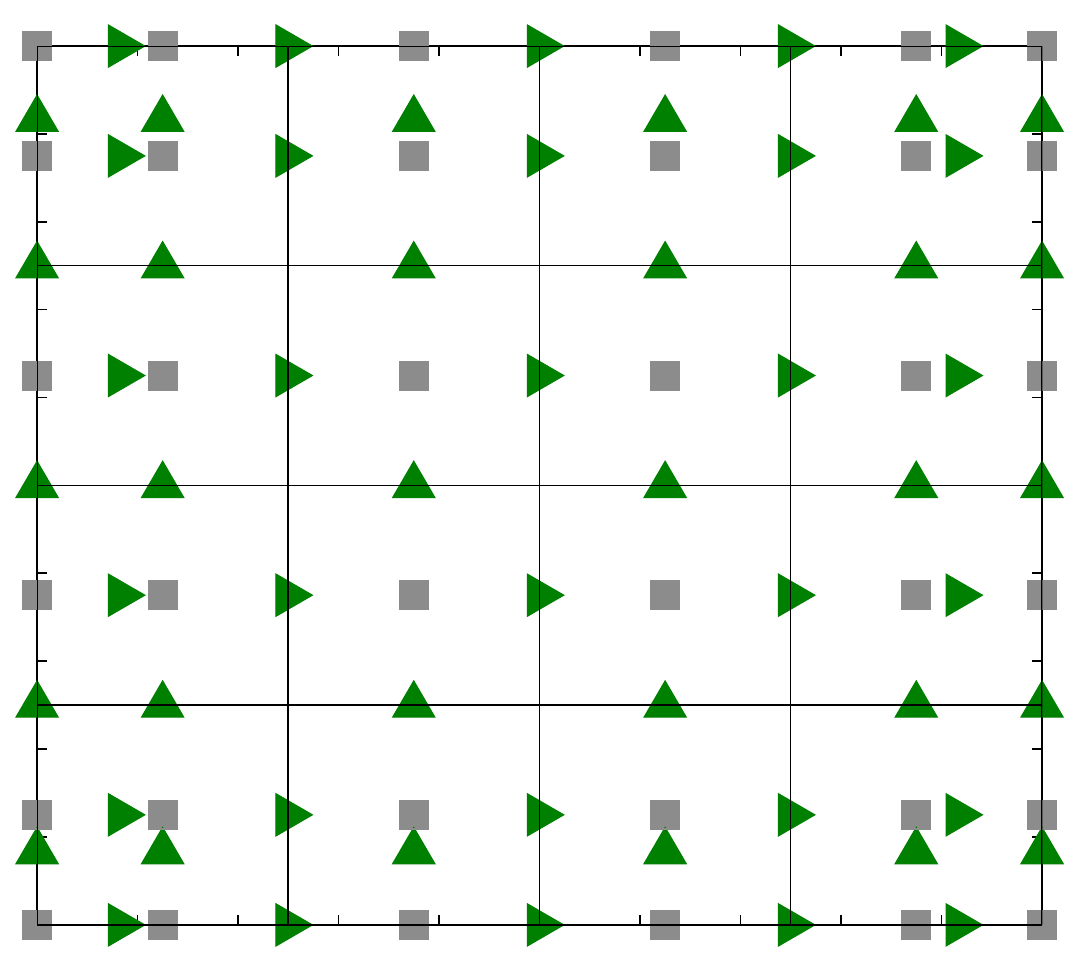}} \\
\caption{Example of collocation grid for $k' = 2$, 4 x 4 elements for the velocity-pressure scheme. Horizontal triangles represent the collocation points for the first momentum equation, vertical triangles represent the points for the second momentum equation, and squares are collocation points for the continuity equation.}
\label{fig:vp_coll_grid}
\end{figure}

\begin{figure}
\centering
\subfloat[Before strong enforcement of normal velocity boundary conditions]{\label{sfig:vvp_grid_a}\includegraphics[width=.4\textwidth]{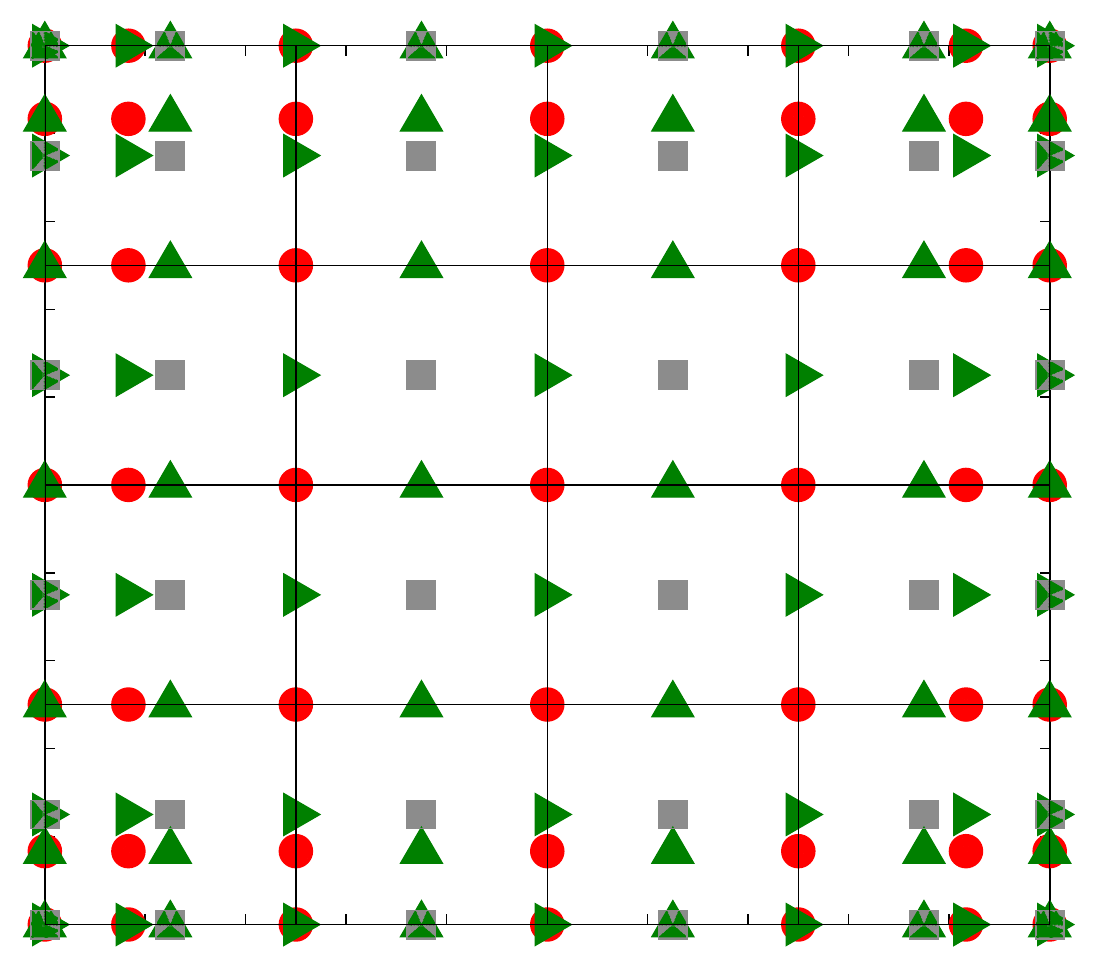}}\hfill
\subfloat[After strong enforcement of normal velocity boundary conditions]{\label{sfig:vvp_grid_b}\includegraphics[width=.4\textwidth]{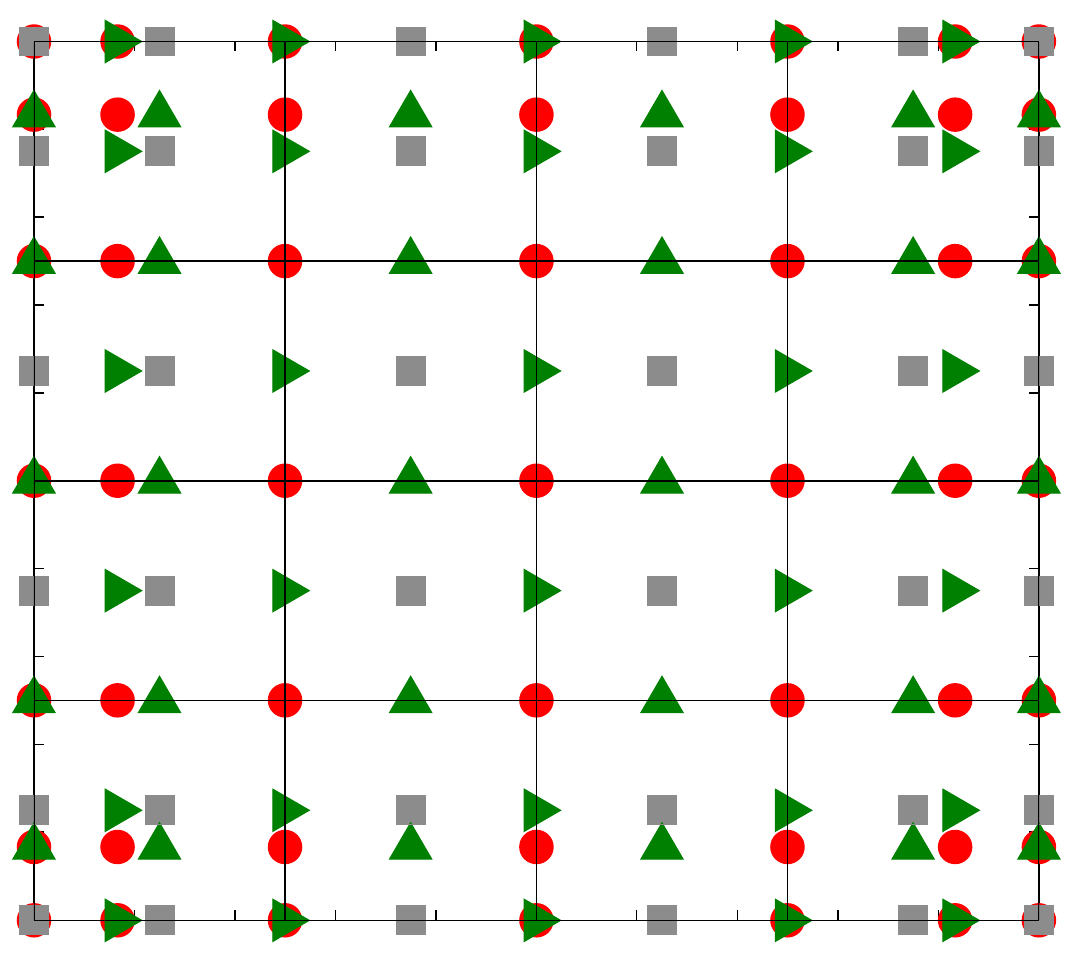}} \\
\caption{Example of collocation grid for $k' = 2$, 4 x 4 elements for the vorticity-velocity-pressure scheme. Horizontal triangles represent the collocation points for the first momentum equation, vertical triangles represent the points for the second momentum equation, squares are collocation points for the continuity equation, and circles are the points for the constitutive relation.}
\label{fig:vvp_coll_grid}
\end{figure}

As each of the discrete unknowns lie in a different B-spline space, each of the governing equations will be collocated at a different set of Greville points. In particular, for both the velocity-pressure formulation and the vorticity-velocity-pressure formulation we use the Greville abscissae associated with the basis of the $x$-velocity ($S^{k_1,k_2-1}_{\boldsymbol{\alpha_1},\boldsymbol{\alpha_2-1}}$) as collocation points for the $x$-momentum equation, the Greville abscissae associated with the basis of the $y$-velocity ($S^{k_1-1,k_2}_{\boldsymbol{\alpha_1-1},\boldsymbol{\alpha_2}}$) as collocation points for the $y$-momentum equation, and the Greville abscissae associated with the basis of the pressure ($S^{k_1-1,k_2-1}_{\boldsymbol{\alpha_1}-1,\boldsymbol{\alpha_2}-1}$) as collocation points for the continuity equation. The constitutive relation within the three field formulation is collocated at the Greville abscissae for the vorticity basis ($S^{k_1,k_2}_{\boldsymbol{\alpha_1},\boldsymbol{\alpha_2}}$). The left of Figure \ref{fig:vp_coll_grid} details an example of this construction for the velocity-pressure scheme, while the left of Figure \ref{fig:vvp_coll_grid} shows example grids for the vorticity-velocity-pressure scheme.

\subsection{Boundary Condition Enforcement}

The last unspecified aspect of the method is enforcement of the Dirichlet boundary conditions. The enforcement of the normal boundary condition is done strongly and collocation points along a boundary for the velocity component orthogonal to that boundary are removed, as the boundary condition specifies the value of the solution at these points. This is shown on the right of Figure \ref{fig:vp_coll_grid} and Figure \ref{fig:vvp_coll_grid}, which depict the same scenarios as their counterparts but with normal boundary conditions enforced.

Enforcement of the tangential boundary condition is slightly more subtle. Recall that in Equation \eqref{eq:2d_weak_nit} we utilized Nitsche's method to enforce this boundary condition. This motivates the enforcement in the velocity-pressure collocation scheme. Indeed if we take this equation and undo the integration by parts, the consistency term vanishes by construction and we are left with just the penalty terms. If we approximate the integral of the test function as done in \cite{deLorenzis_Neumann_contact} the collocated momentum equations will be of the form

\begin{equation}
     -\nu \Delta \mathbf{u}^h + \mathbf{u}^h \cdot \nabla \mathbf{u}^h + \nabla p^h + \frac{C_{pen}^2}{h^2}(\mathbf{u}^h - \mathbf{g}) = \mathbf{f},
\end{equation}

\noindent where $C_{pen}$ is a penalty constant and $h$ is the Greville mesh size perpendicular to the boundary. Note that because this construction is used to only enforce the tangential boundary conditions, this penalty term only appears in the equation for the momentum balance along the tangential direction of each boundary. 

In the vorticity-velocity-pressure scheme we do not use the same Nitsche-based approach. The method utilized here is directly inspired by the Enhanced Collocation method for enforcing Neumann boundary conditions in isogeometric collocation schemes \cite{deLorenzis_Neumann_contact}. 

We start by considering the weak form of the constitutive relation given by Equation \eqref{eq:weak_const_2d}. The final term on the left hand side is the boundary term which would be used to enforce natural boundary conditions in a Galerkin method. In a similar vein to the Enhanced Collocation approach, we can undo the integration by parts to arrive at

\begin{equation}
    \int_{\Omega} \psi^h (\omega^h - (\frac{\partial u^h_y}{\partial x} - \frac{\partial u^h_x}{\partial y})) d \Omega + \int_{\partial \Omega} \psi^h (\mathbf{u}^h \cdot \mathbf{s} - \mathbf{g} \cdot \mathbf{s}) ds = 0.
\end{equation}

By approximating the integrals of the test functions as done in \cite{deLorenzis_Neumann_contact, gomez_cauchygalerkin} we arrive at a modified strong form statement of the constitutive relation which can be collocated along the boundaries

\begin{equation}
    \omega^h - (\frac{\partial u^h_y}{\partial x} - \frac{\partial u^h_x}{\partial y}) + \frac{C_{pen}}{h}(\mathbf{u}^h \cdot \mathbf{s} - \mathbf{g} \cdot \mathbf{s}) = 0,
\end{equation}

\noindent where again $C_{pen}$ is a penalty constant and $h$ is the Greville mesh size perpendicular to the boundary.


\subsection{Final Collocated Equations}

Finally, the results of the previous sections are collected and we present the final form of the discrete equations used to solve for the discrete unknowns. The velocity-pressure scheme is considered first. Let us define $\boldsymbol{\tau}^{u_x}_\ell$ for $\ell = 1,...,M^{u_x}$ to be the set of Greville points for $S^{k_1,k_2-1}_{\boldsymbol{\alpha_1},\boldsymbol{\alpha_2-1}}$ with the points corresponding to no-penetration boundaries removed as discussed previously. Define in a similar manner $\boldsymbol{\tau}^{u_y}_\ell$ for $\ell = 1,...,M^{u_y}$, which are the Greville points of $S^{k_1-1,k_2}_{\boldsymbol{\alpha_1-1},\boldsymbol{\alpha_2}}$ with no-penetration boundary points removed. Lastly, $\boldsymbol{\tau}^{p}_\ell$ for $k = 1,...,N^{p}$ are the Greville points of $\mathcal{Q}_h$. Then the discrete 2D problem reads: 

\bigskip

$$
\left\{ \hspace{5pt}
\parbox{6in}{
\noindent Find $\textbf{u}^h \in \boldsymbol{\mathcal{V}}_{h,\mathbf{g}}$ and $P^h \in \mathcal{Q}_{h,0}$ such that:
\begin{equation}
\begin{split}
    \left(-\nu \frac{\partial^2 u_x^h}{\partial x^2} - \nu \frac{\partial^2 u_x^h}{\partial y^2} + u^h_x\frac{\partial u^h_x}{\partial x} + u^h_y\frac{\partial u^h_x}{\partial y} + \frac{\partial p^h}{\partial x} \right)(\boldsymbol{\tau}^{u_x}_\ell)  = f_x(\boldsymbol{\tau}^{u_x}_\ell) \quad \forall \boldsymbol{\tau}^{u_x}_\ell \in \Omega
    \label{eq:2D_xmom_discrete_int}
\end{split}
\end{equation}
\begin{equation}
\begin{split}
    \left(-\nu \frac{\partial^2 u_x^h}{\partial x^2} - \nu \frac{\partial^2 u_x^h}{\partial y^2} + u^h_x\frac{\partial u^h_x}{\partial x} + u^h_y\frac{\partial u^h_x}{\partial y} + \frac{\partial p^h}{\partial x} + \frac{C_{pen}^2}{h^2}(u^h_x - g_x)\right)(\boldsymbol{\tau}^{u_x}_\ell)  \\= f_x(\boldsymbol{\tau}^{u_x}_\ell) \quad \forall \boldsymbol{\tau}^{u_x}_\ell \in \partial \Omega
    \label{eq:2D_xmom_discrete_bdry}
\end{split}
\end{equation}
\begin{equation}
\begin{split}
    \left(-\nu \frac{\partial^2 u_y^h}{\partial x^2} - \nu \frac{\partial^2 u_y^h}{\partial y^2} + u^h_x\frac{\partial u^h_y}{\partial x} + u^h_y\frac{\partial u^h_y}{\partial y} + \frac{\partial p^h}{\partial y} \right)(\boldsymbol{\tau}^{u_y}_\ell) = f_y(\boldsymbol{\tau}^{u_y}_\ell) \quad \forall  \boldsymbol{\tau}^{u_y}_\ell \in \Omega
    \label{eq:2D_ymom_discrete_int}
\end{split}
\end{equation} 
\begin{equation}
\begin{split}
    \left(-\nu \frac{\partial^2 u_y^h}{\partial x^2} - \nu \frac{\partial^2 u_y^h}{\partial y^2} + u^h_x\frac{\partial u^h_y}{\partial x} + u^h_y\frac{\partial u^h_y}{\partial y} + \frac{\partial p^h}{\partial y} + \frac{C_{pen}^2}{h^2}(u^h_y - g_{y})\right)(\boldsymbol{\tau}^{u_y}_\ell) \\= f_y(\boldsymbol{\tau}^{u_y}_\ell) \quad \forall  \boldsymbol{\tau}^{u_y}_\ell \in \partial \Omega
    \label{eq:2D_ymom_discrete_bdry}
\end{split}
\end{equation} 
\begin{equation}
    \left(\frac{\partial u^h_x}{\partial x} + \frac{\partial u^h_y}{\partial y}\right)(\boldsymbol{\tau}^{p}_\ell)  = 0 \quad \forall \boldsymbol{\tau}^{p}_\ell \in \Omega \cup \partial \Omega.
\end{equation} 
}
\right.
$$

\bigskip

In the above we have split the momentum equations into expressions valid on the interior collocation points (Equations \eqref{eq:2D_xmom_discrete_int} and \eqref{eq:2D_ymom_discrete_int}) and expressions valid on the remaining boundary collocation points (Equations \eqref{eq:2D_xmom_discrete_bdry} and \eqref{eq:2D_ymom_discrete_bdry}).

Similarly, for the vorticity-velocity-pressure scheme we also define $\boldsymbol{\tau}^{\omega}_\ell$ for $\ell = 1,...,N^{\omega}$ as the Greville points of $\Psi_h$. With this scheme the discrete 2D problem reads: 

\bigskip

$$
\left\{ \hspace{5pt}
\parbox{6in}{
\noindent Find $\textbf{u}^h \in \boldsymbol{\mathcal{V}}_{h,\mathbf{g}}$, $P^h \in \mathcal{Q}_{h,0}$, and $\omega^h \in \Psi_h$ such that:
\begin{equation}
    \left( \nu \frac{\partial \omega^h}{\partial y} - \omega^hu^h_y + \frac{\partial P^h}{\partial x} \right)(\boldsymbol{\tau}^{u_x}_\ell) = f_x(\boldsymbol{\tau}^{u_x}_\ell) \quad \forall \boldsymbol{\tau}^{u_x}_\ell \in \Omega \cup \partial \Omega
\end{equation}
\begin{equation}
    \left( - \nu \frac{\partial \omega^h}{\partial x} + \omega^h u^h_x + \frac{\partial P^h}{\partial y} \right)(\boldsymbol{\tau}^{u_y}_\ell) = f_y(\boldsymbol{\tau}^{u_y}_\ell) \quad \forall \boldsymbol{\tau}^{u_y}_\ell \in \Omega \cup \partial \Omega
\end{equation}
\begin{equation}
    \left( \frac{\partial u^h_x}{\partial x} + \frac{\partial u^h_y}{\partial y} \right)(\boldsymbol{\tau}^{p}_\ell)  = 0 \quad \forall \boldsymbol{\tau}^{p}_\ell \in \Omega \cup \partial \Omega
\end{equation}
\begin{equation}
  \left(\omega^h - \frac{\partial u^h_y}{\partial x} - \frac{\partial u^h_x}{\partial y}  \right)(\boldsymbol{\tau}^{\omega}_\ell)= 0 \quad \forall \boldsymbol{\tau}^{\omega}_\ell \in \Omega
  \label{eq:2D_const_int}
\end{equation}
\begin{equation}
  \left(\omega^h - \frac{\partial u^h_y}{\partial x} - \frac{\partial u^h_x}{\partial y} +\frac{C_{pen}}{h}(\mathbf{u}^h \cdot \mathbf{s} - \mathbf{g} \cdot \mathbf{s} ) \right)(\boldsymbol{\tau}^{\omega}_\ell)= 0 \quad \forall \boldsymbol{\tau}^{\omega}_\ell \in \partial \Omega.
  \label{eq:2D_const_bdry}
\end{equation}
}
\right.
$$

\bigskip

In the three field formulation we split the constitutive law into an expression for the interior collocation points (Equation \eqref{eq:2D_const_int}) and another expression for boundary collocation points (Equation \eqref{eq:2D_const_bdry}).





Resulting from these equations are nonlinear systems of equations which we can use to solve for the unknown coefficients of velocity, pressure, and vorticity using a Newton-Raphson method. 

\subsection{Proof of Divergence Conforming Property}

From the commuting diagrams our spaces form, shown in Equations \eqref{eq:deRham_discrete} and \eqref{eq:deRham2D_discrete}, it is simple to show that both of the resulting collocation methods return an exactly pointwise divergence free velocity field. The commuting diagrams reveal that the divergence of the discrete velocity lies within the discrete pressure space, $\nabla \cdot \mathbf{u}^h \in \mathcal{Q}_h$. We can thus equivalently write the divergence of the velocity as a linear combination of the pressure basis functions:

\begin{equation}
    \nabla \cdot \mathbf{u}^h = \sum_{i = 0}^{N^p} c_i q_i,
\end{equation}

\noindent where $q_i \in \mathcal{Q}_h$ are the basis functions for the pressure and $ c_i \in \mathbb{R}$ are the coefficients. As part of the collocation scheme, we strongly enforce that the velocity field has zero divergence at a number of collocation points equal to the dimension of the discrete pressure space. This condition can be written as a system of linear equations 

\begin{equation}
    \mathbf{M} \mathbf{c} = \mathbf{0},
    \label{eq:div_free_system}
\end{equation}

\noindent where $\mathbf{M}$ is a square matrix whose entries are the pressure basis functions evaluated at each collocation point and $\mathbf{c}$ is the vector of coefficients. 

If the choice of collocation points yields a set of linearly independent equations, that is to say $\mathbf{M}$ is invertible, then we know that the solution to Equation \eqref{eq:div_free_system} is $\mathbf{c} = \mathbf{0}$, and thus the velocity field is exactly divergence free pointwise.

\section{Numerical Results on Square Domains}

The developed schemes are now tested on multiple 2D problems on the unit square. First, a manufactured vortex problem is considered to experimentally compute the convergence rates of the error and test for pressure and Reynolds number robustness. Then, the classical lid-driven cavity problem is considered at a variety of Reynolds numbers.

\subsection{Two-Dimensional Manufactured Solution}

As a first numerical experiment, we consider a vortex problem on the unit square constructed using the method of manufactured solutions. This solution was originally developed in \cite{buffa_stokes} and employs the following velocity and pressure fields:

\begin{equation}
    \tilde{\textbf{u}} = \left[
\begin{array}{c}
2e^x(-1+x)^2x^2(y^2-y)(-1+2y) \\
(-e^x(-1+x)x(-2+x(3+x))(-1+y)^2y^2)
\end{array}
\right],
\end{equation}

\begin{equation}
    \left.
\begin{array}{ccc}
\tilde{p} & = & (-424+156e+(y^2-y)(-456+e^x(456+x^2(228-5(y^2-y))+ \\ &&
2x(-228+(y^2-y))+2x^3(-36+(y^2-y))+x^4(12+(y^2-y))))).
\end{array}
\right.
\end{equation}

For the velocity-pressure scheme we define the forcing term to be

\begin{equation}
    \textbf{f} = -\nu \Delta \tilde{\textbf{u}} + \tilde{\textbf{u}} \cdot \nabla \tilde{\textbf{u}} + \nabla \tilde{p} ,
\end{equation}

\noindent while for the vorticity-velocity-pressure formulation we define the forcing term in the momentum equations to be: 

\begin{equation}
    \textbf{f} = \nu \nabla^\perp \tilde{\omega} + \tilde{\omega} \times \tilde{\textbf{u}} + \nabla \tilde{P} ,
\end{equation}

\noindent with $\tilde{\omega} = \nabla \times \tilde{\textbf{u}}$ and $\tilde{P} = \tilde{p} + \frac{1}{2}(\tilde{\mathbf{u}} \cdot \tilde{\mathbf{u}})$. Enforcing homogeneous boundary conditions and requiring that the integral of pressure is zero, it is clear to see that the velocity and kinematic pressure solutions should be equal to $\tilde{\textbf{u}}$ and $\tilde{p}$. 

\begin{figure}
\centering
\subfloat[Velocity $L^2$ error]{\label{sfig:2D_conv_lap_a}\includegraphics[width=.5\textwidth]{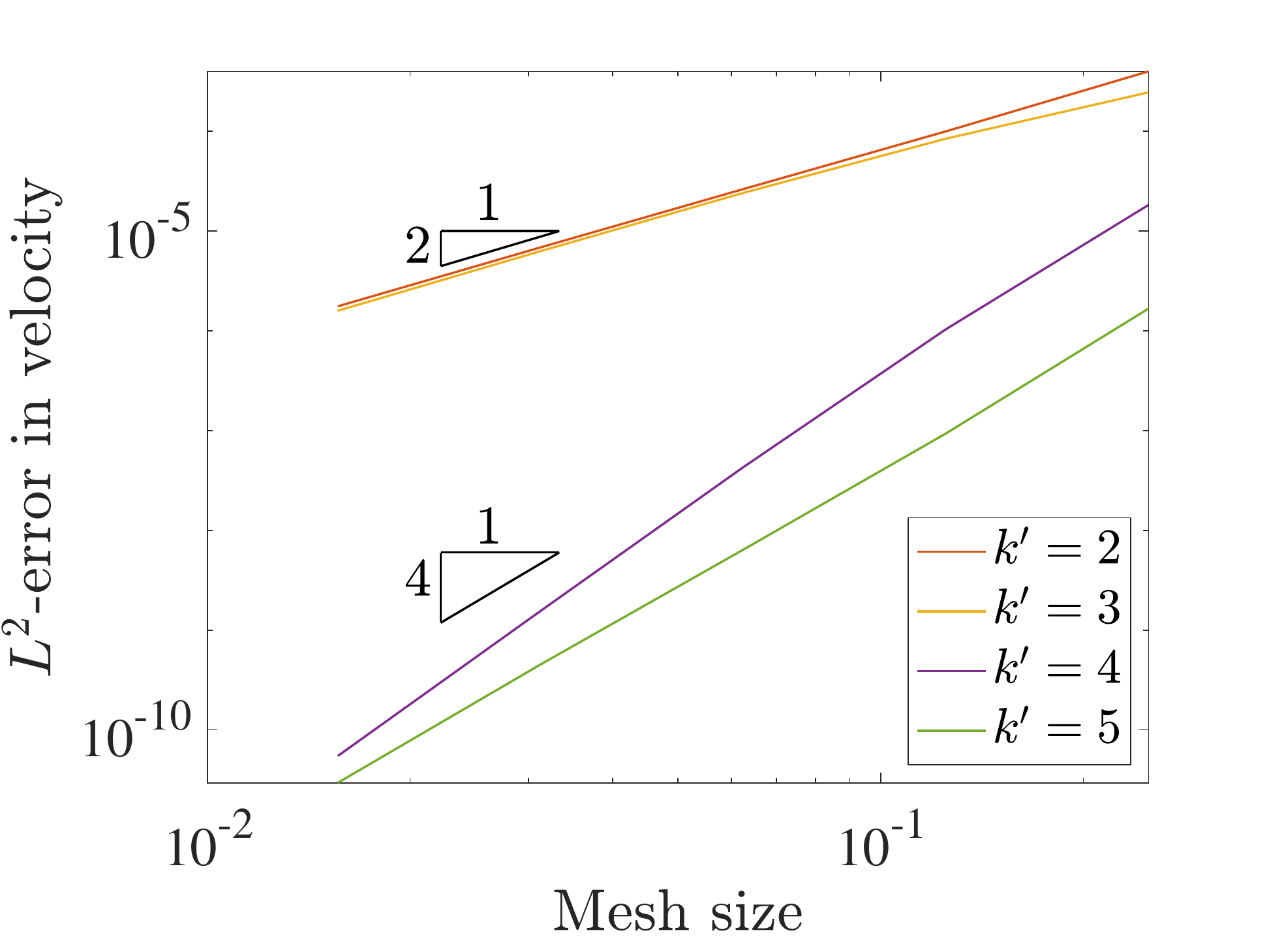}}\hfill
\subfloat[Velocity $H^1$ error]{\label{sfig:2D_conv_lap_b}\includegraphics[width=.5\textwidth]{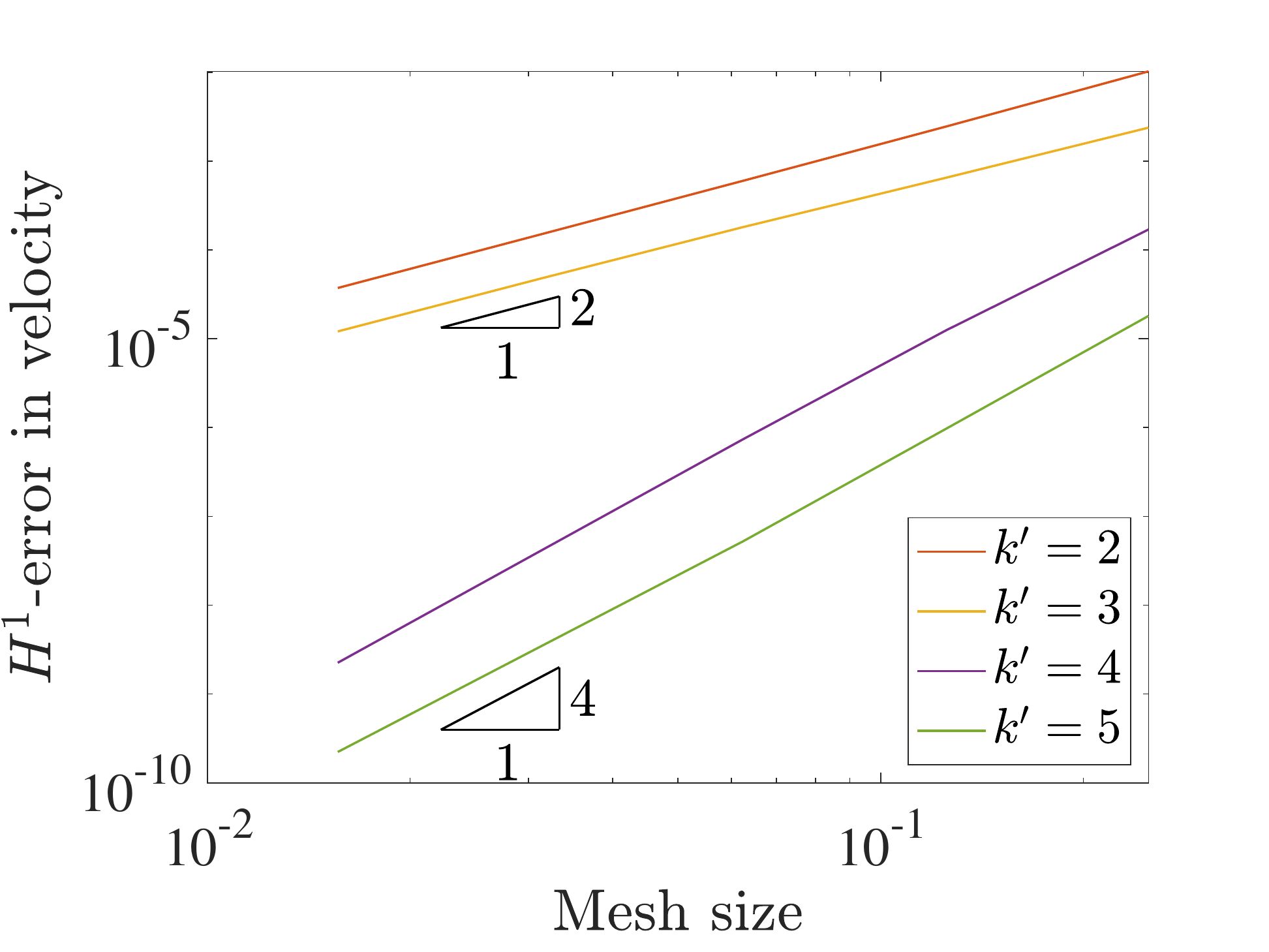}} \\
\subfloat[Pressure $L^2$ error]{\label{sfig:2D_conv_lap_c}\includegraphics[width=.5\textwidth]{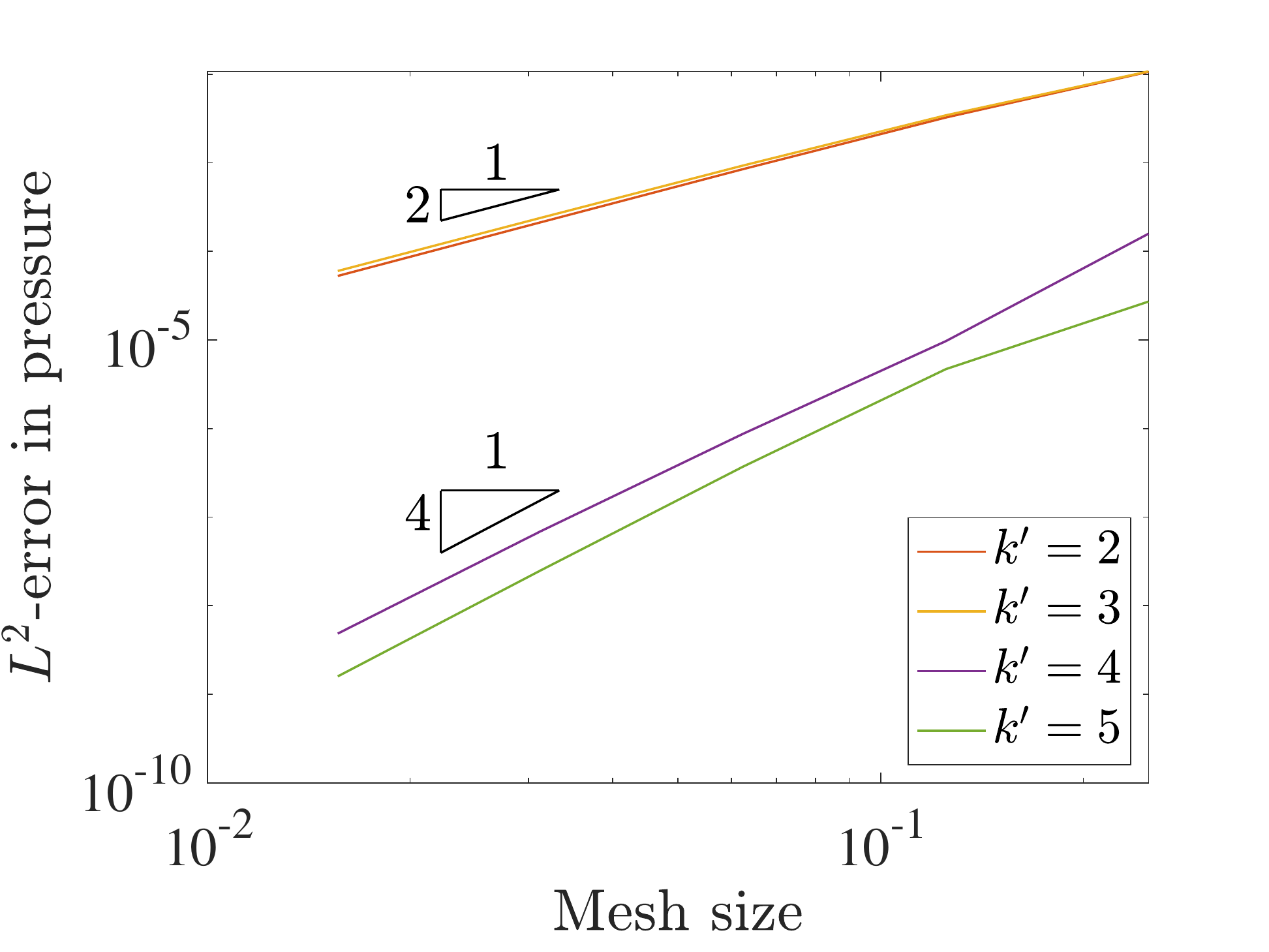}}\hfill
\subfloat[Pressure $H^1$ error]{\label{sfig:2D_conv_lap_d}\includegraphics[width=.5\textwidth]{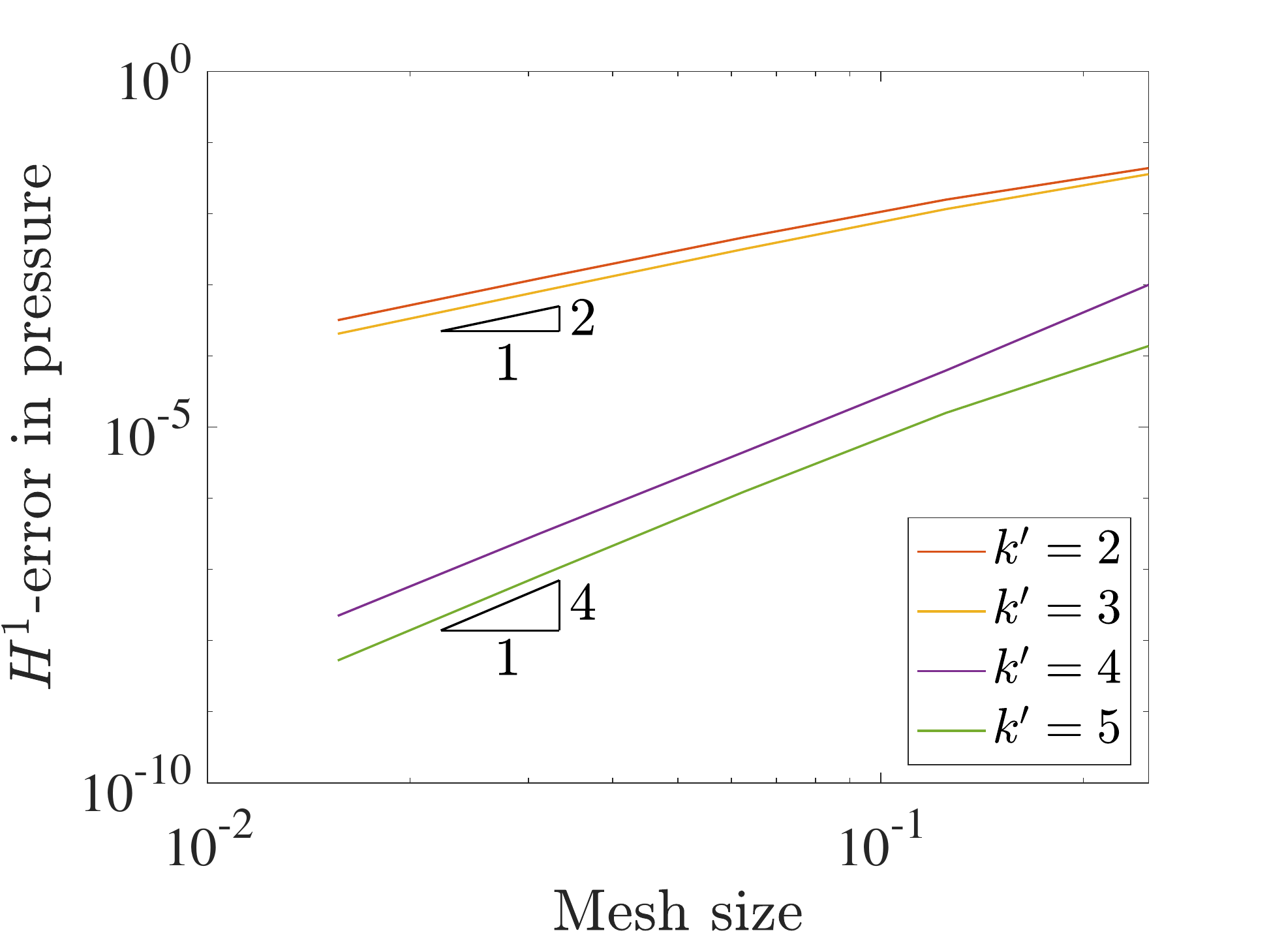}} \\
\caption{Errors in 2D manufactured vortex solution for velocity-pressure formulation}
\label{fig:2D_conv_lap}
\end{figure}

To understand the accuracy of this collocation method, we test the convergence rates on a variety of grids and with different degree B-spline bases. For this test we set the Reynolds number to $Re = \frac{1}{\nu} = 1$. We measure error using the $L^2$ norm as well as the $H^1$ semi-norm. Figure \ref{fig:2D_conv_lap} details the convergence rates of velocity and pressure when using the two field formulation. In this case the errors in both velocity and pressure converge at a rate of $k'$ when $k'$ is even and $k'-1$ for odd $k'$. These are the standard, expected rates that have been seen in other studies of isogeometric collocation, and are one and two orders suboptimal in $L^2$ for odd and even $k'$.

\begin{figure}
\centering
\subfloat[Velocity $L^2$ error]{\label{sfig:2D_conv_a}\includegraphics[width=.5\textwidth]{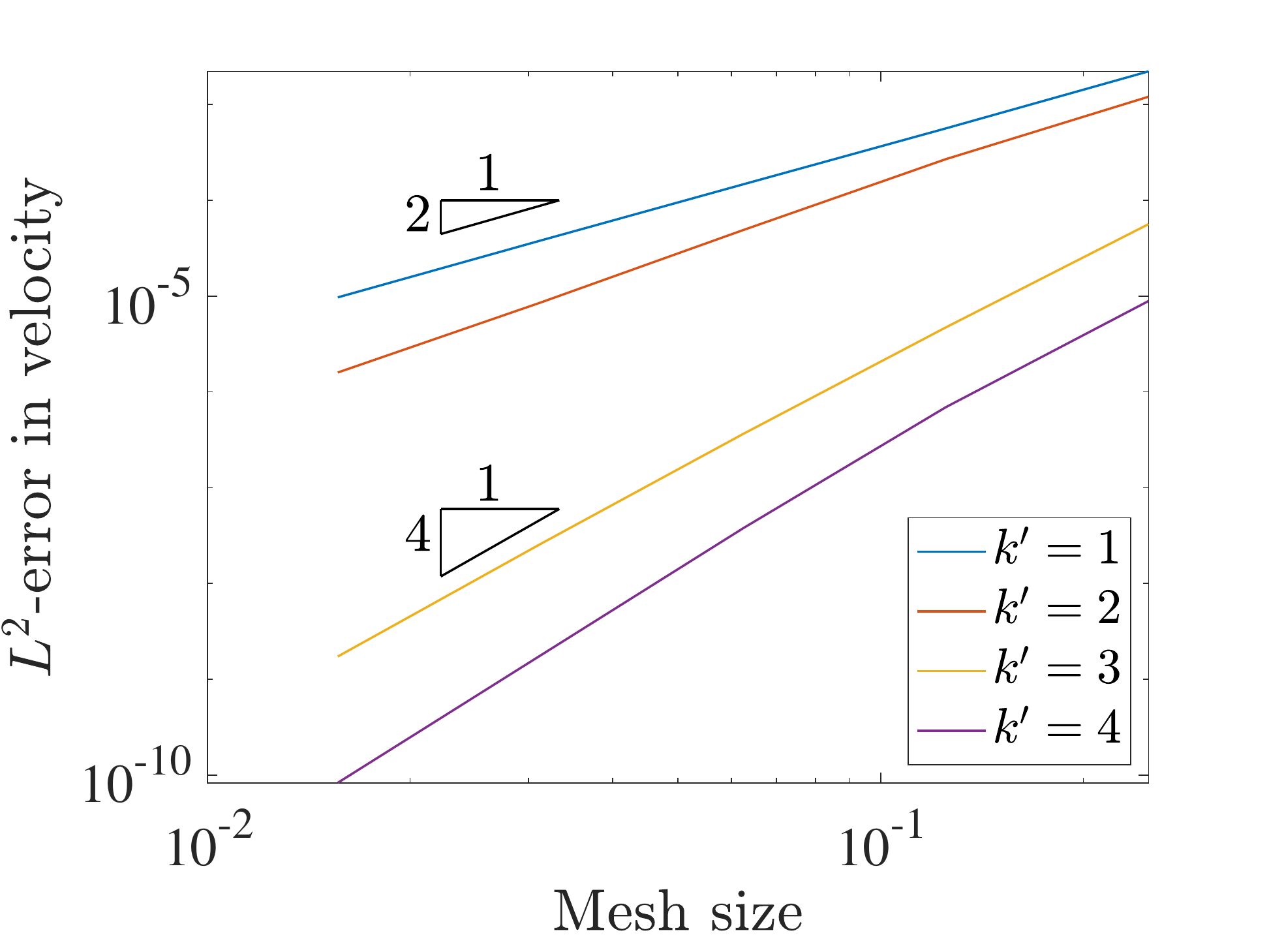}}\hfill
\subfloat[Velocity $H^1$ error]{\label{sfig:2D_conv_b}\includegraphics[width=.5\textwidth]{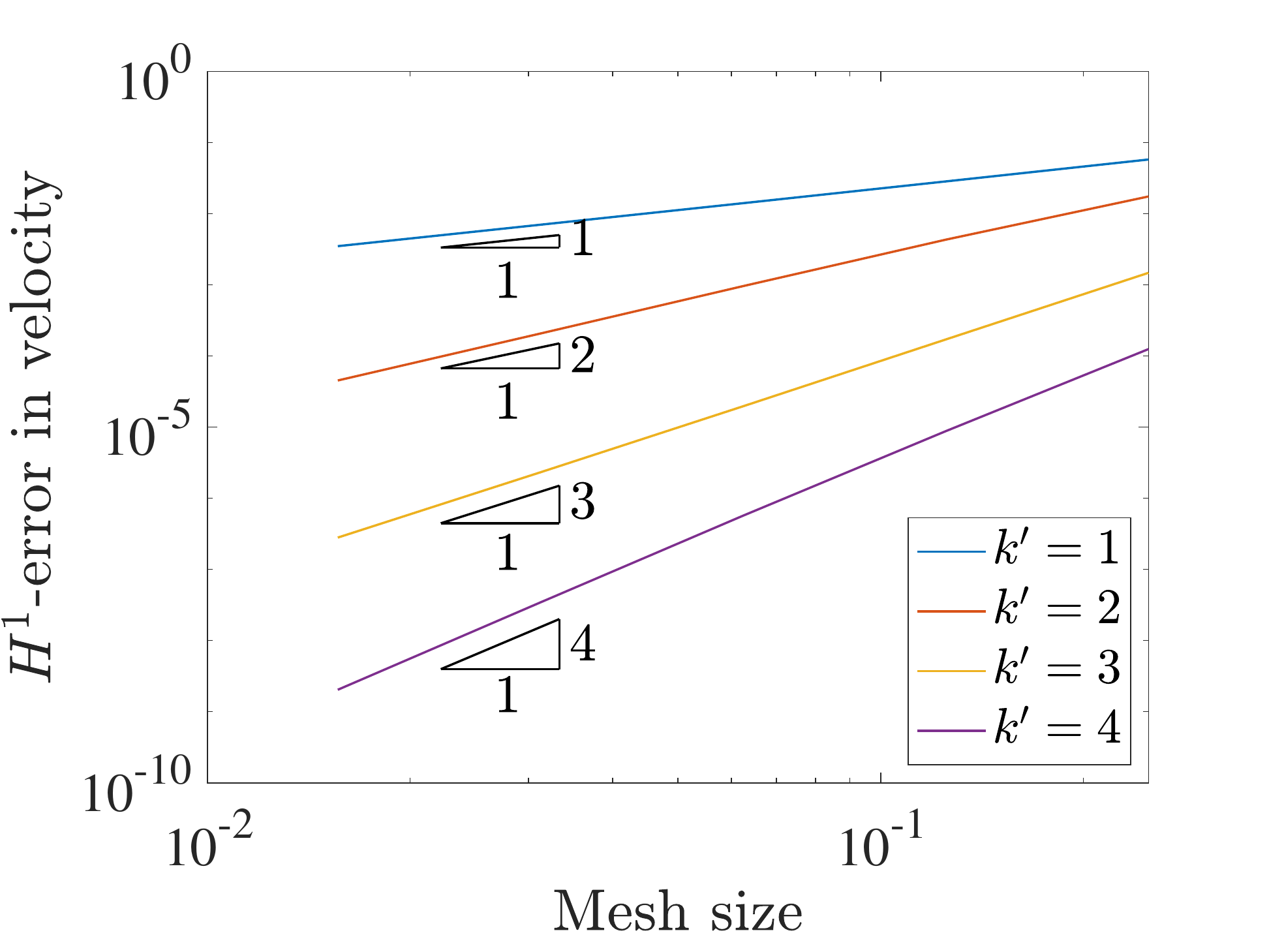}} \\
\subfloat[Pressure $L^2$ error]{\label{sfig:2D_conv_c}\includegraphics[width=.5\textwidth]{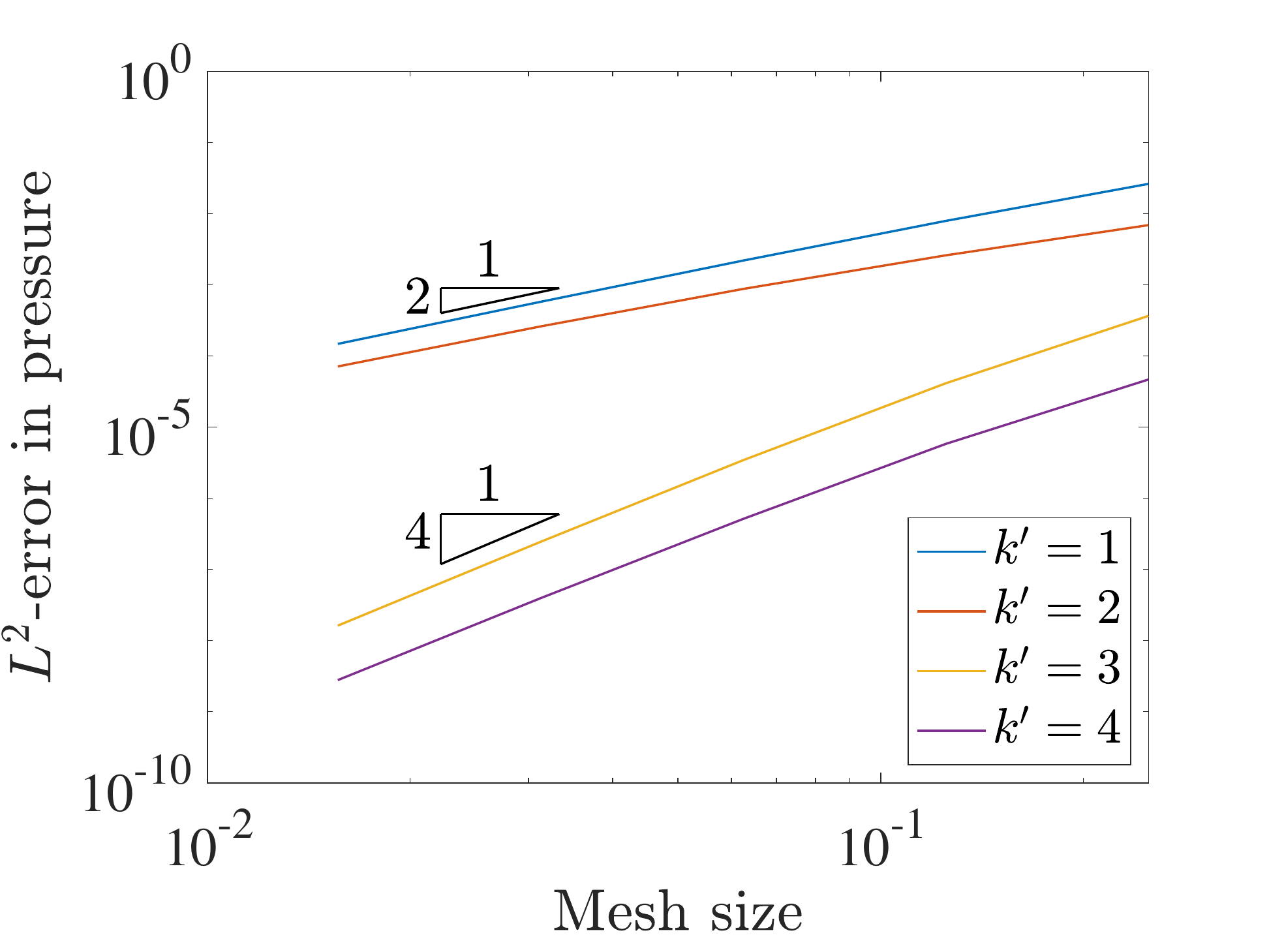}}\hfill
\subfloat[Pressure $H^1$ error]{\label{sfig:2D_conv_d}\includegraphics[width=.5\textwidth]{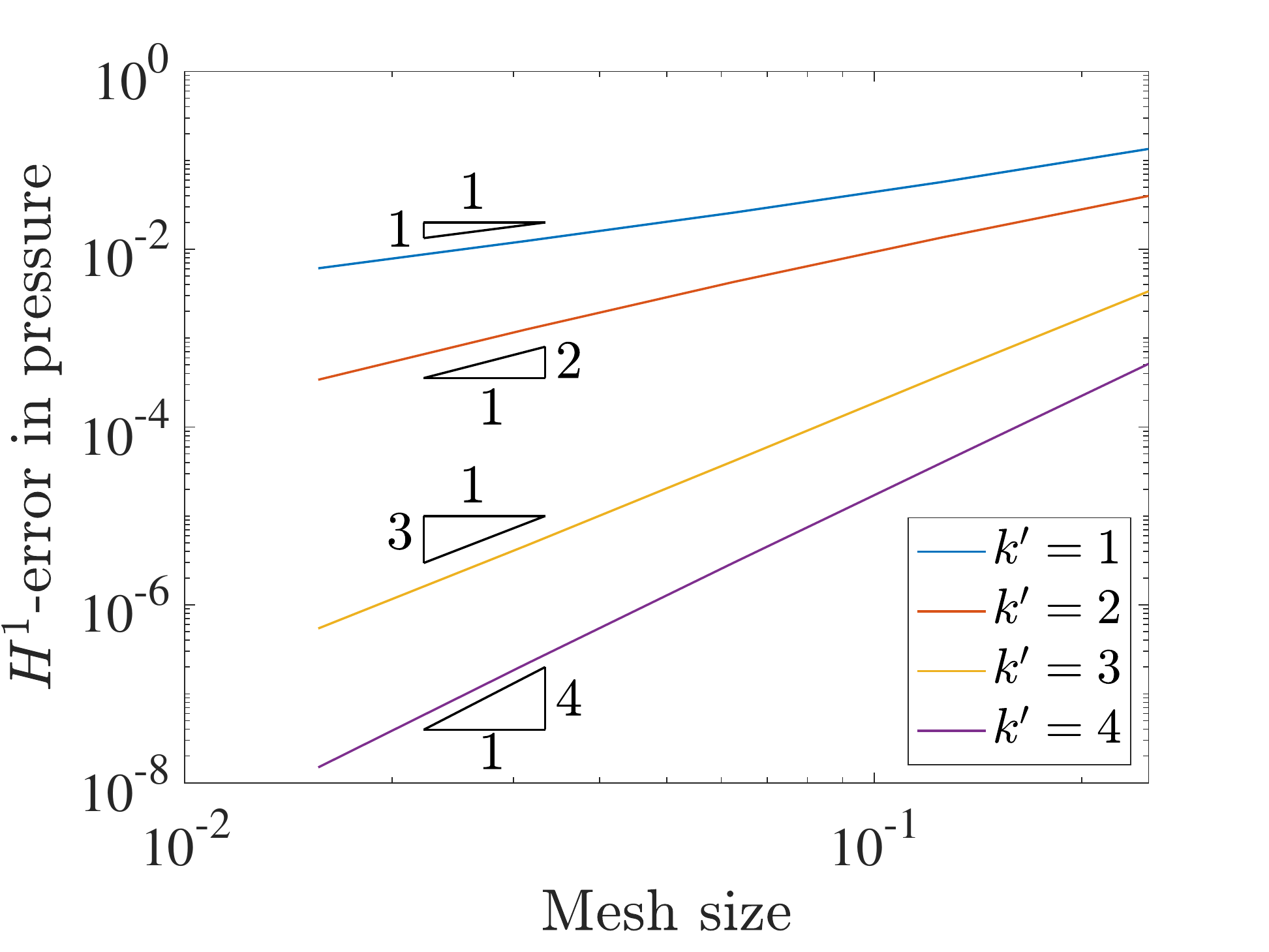}} \\
\subfloat[Vorticity $L^2$ error]{\label{sfig:2D_conv_e}\includegraphics[width=.5\textwidth]{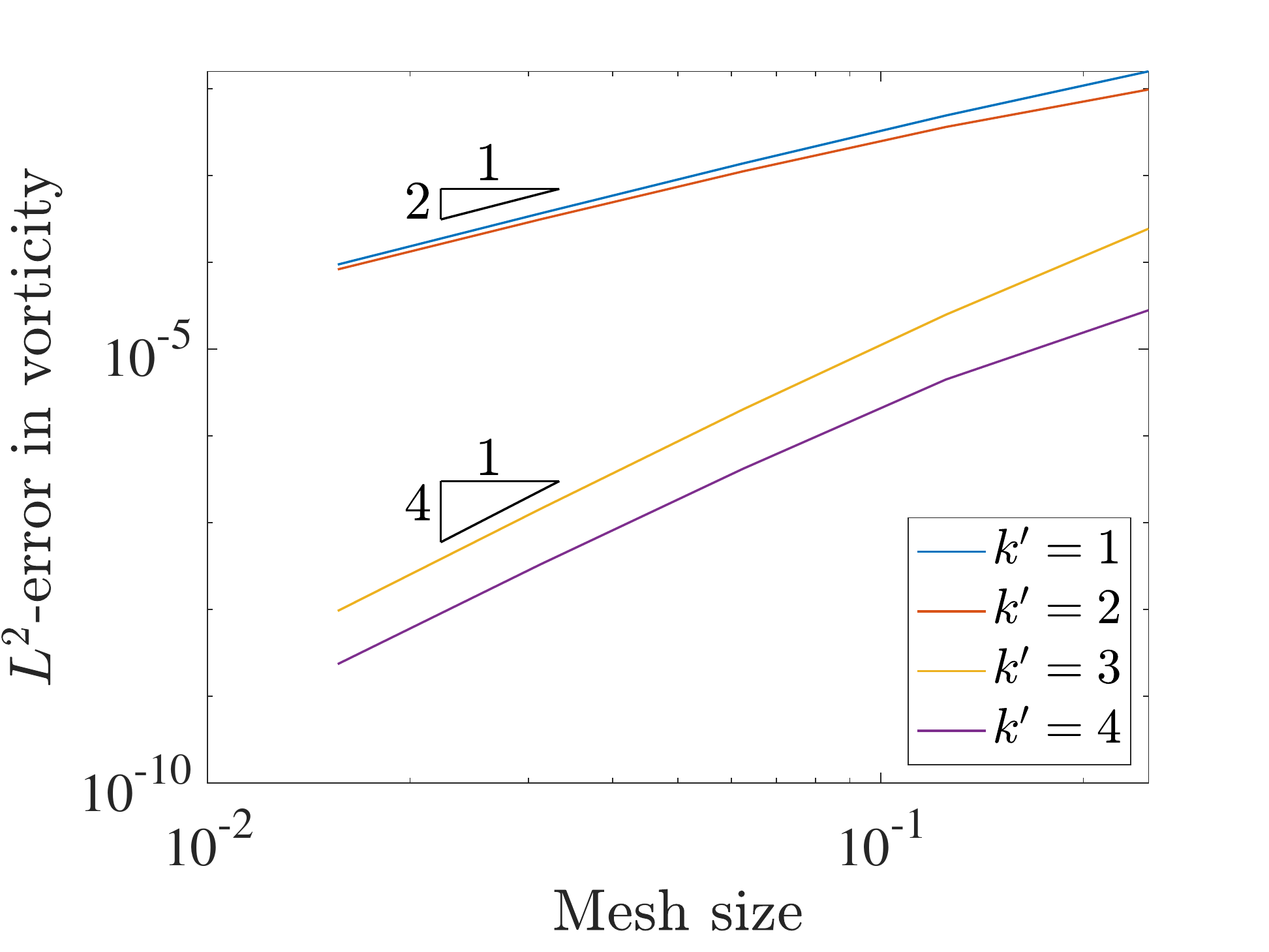}}\hfill
\subfloat[Vorticity $H^1$ error]{\label{sfig:2D_conv_f}\includegraphics[width=.5\textwidth]{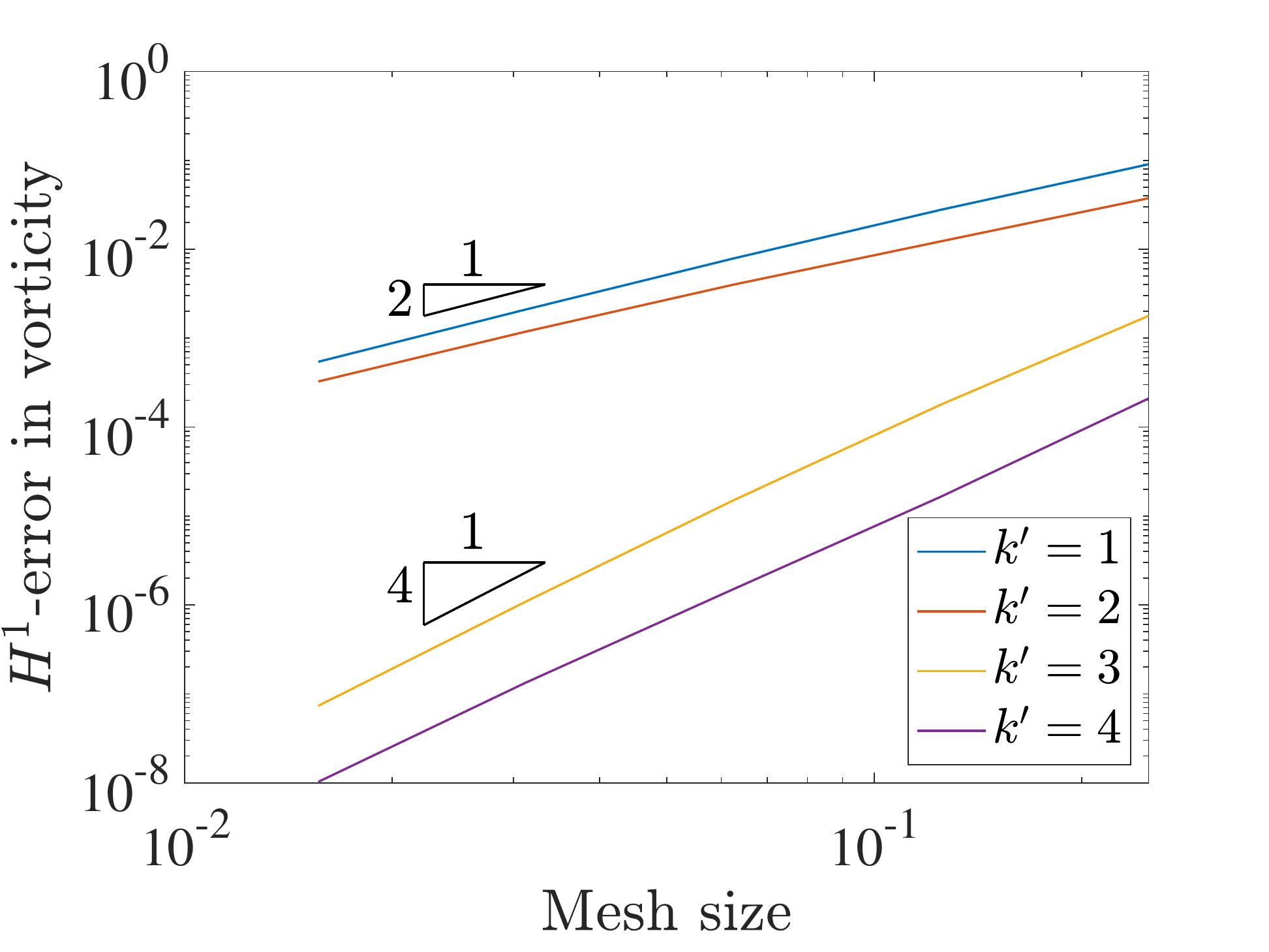}} \\
\caption{Errors in 2D manufactured vortex solution for vorticity-velocity-pressure formulation}
\label{fig:2D_conv}
\end{figure}

Figure \ref{fig:2D_conv} details the convergence of velocity, kinematic pressure, and vorticity as we refine the meshes in the three field scheme. Using this scheme, all of the unknowns converge in the $L^2$ norm at a rate of approximately $k'$ for even values of $k'$ and at a rate of $k'+1$ for odd values of $k'$. These rates match the rates achieved using even $k'$ in the two field formulation, and these rates are two orders faster for odd $k'$. In fact, this formulation recovers optimal convergence rates in the $L^2$ norm for odd $k'$. 

In the $H^1$ semi-norm, we see convergence rates of $k'$ for all polynomial degrees for the velocity and pressure. These rates are optimal in the $H^1$ semi-norm for all degrees and again are as fast or better than the corresponding velocity-pressure scheme results. Interestingly, the $H^1$ convergence of vorticity seems to be at a rate of $k'+1$ for odd $k'$ and a rate of $k'$ for even values. 

\begin{figure}
\centering
\subfloat[Velocity $L^2$ error]{\label{sfig:2D_conv_comp_a}\includegraphics[width=.5\textwidth]{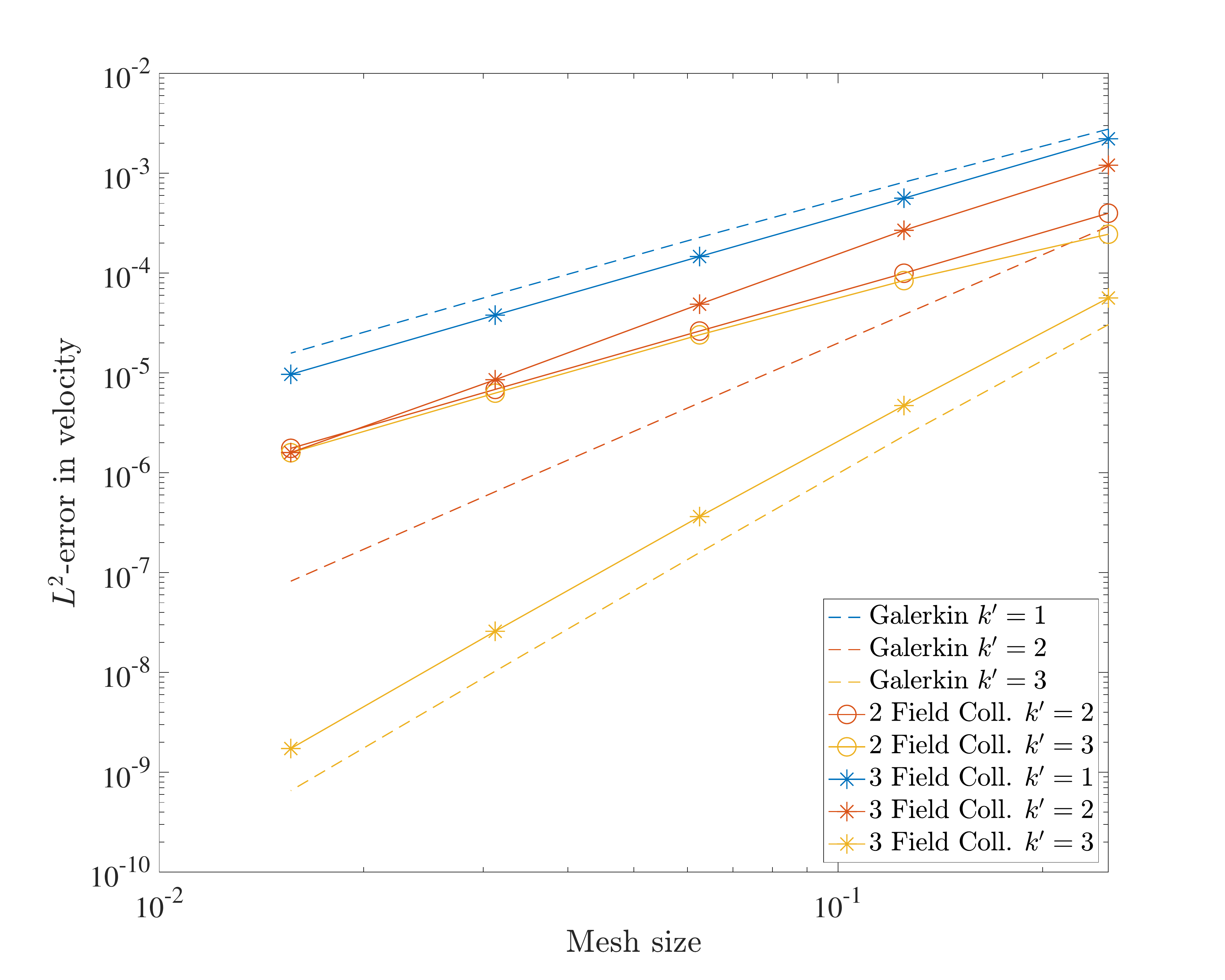}}\hfill
\subfloat[Velocity $H^1$ error]{\label{sfig:2D_conv_comp_b}\includegraphics[width=.5\textwidth]{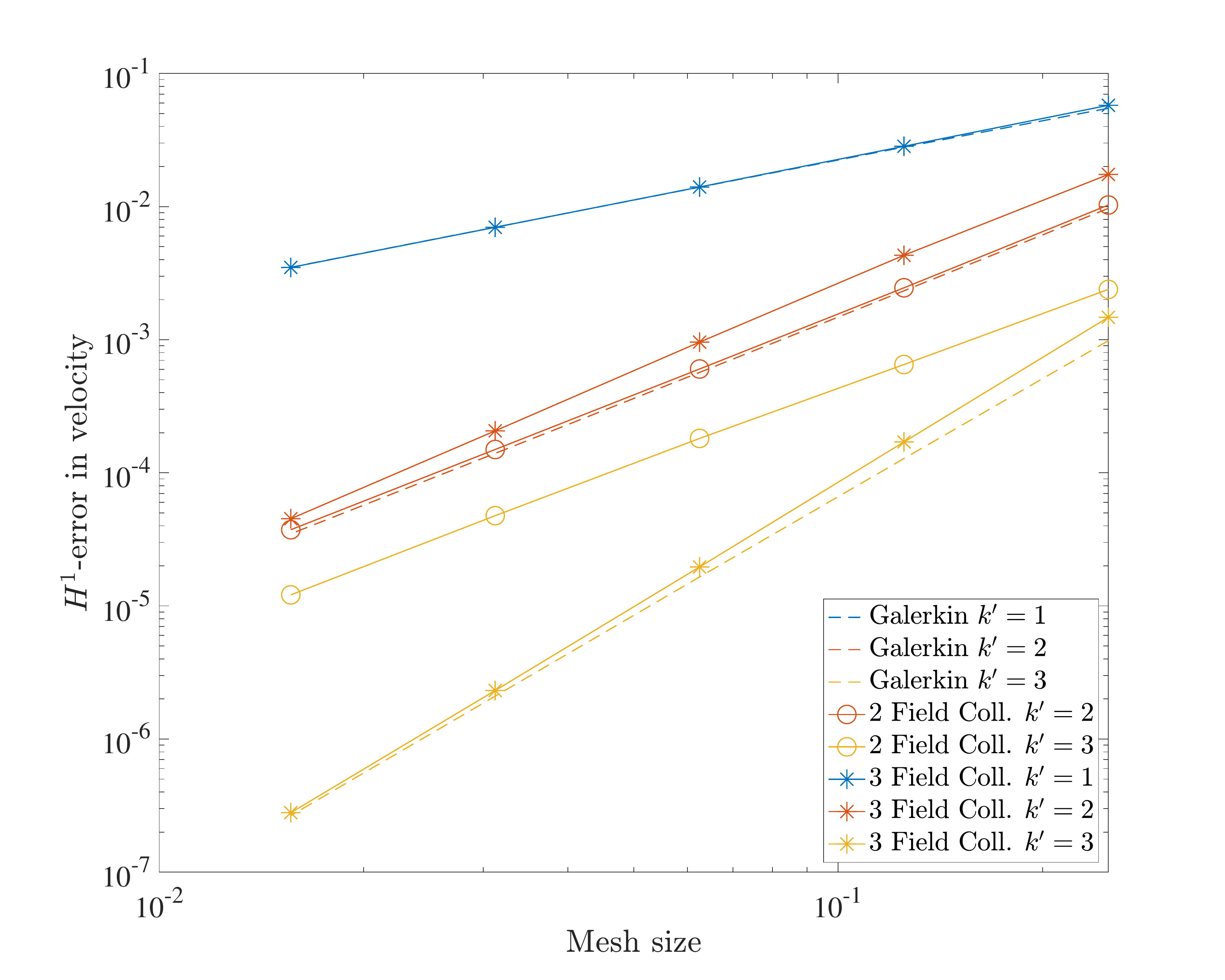}} \\
\subfloat[Pressure $L^2$ error]{\label{sfig:2D_conv_comp_c}\includegraphics[width=.5\textwidth]{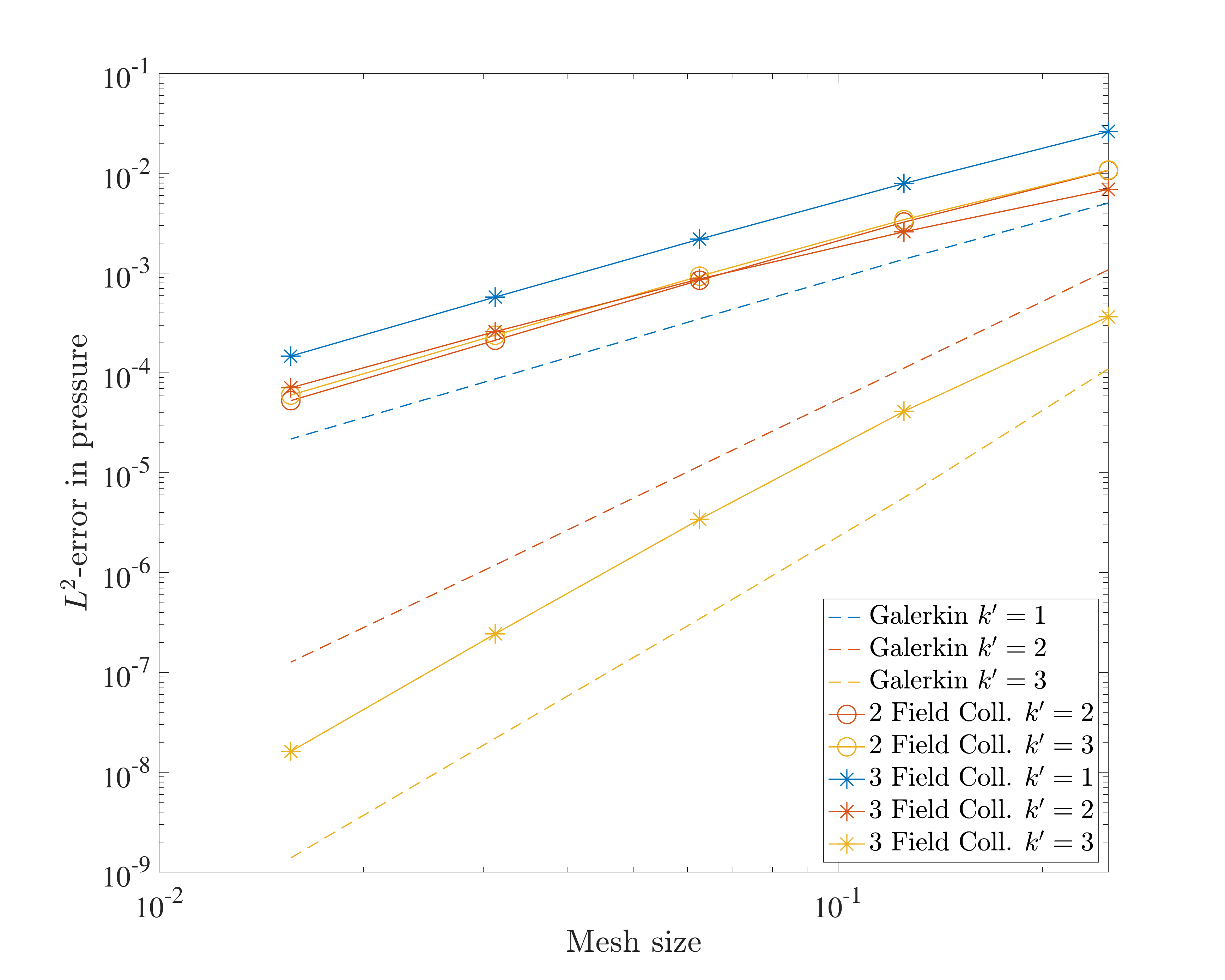}} \\
\caption{Errors in 2D manufactured vortex solution comparison}
\label{fig:2D_conv_comp}
\end{figure}

To further study our new collocation schemes, we can also directly compare the errors produced with divergence-conforming Galerkin schemes of the same orders. Figure \ref{fig:2D_conv_comp} shows the $L^2$ norm and $H^1$ semi-norm errors in velocity as well as the $L^2$ errors in pressure for both collocation schemes along with the Galerkin results for the same problem from \cite{evans_steady_NS}. This comparison highlights the severe suboptimality of the velocity-pressure results with odd $k'$. We also note that the $H^1$ errors obtained with the three field formulation nearly match the Galerkin results. 

\subsection{Pressure Robustness}

Next we perform some ancillary tests related to the manufactured solution to test some secondary robustness properties of the method. The first test relates to so-called pressure robustness \cite{john_divconstraint}. In particular, we take the kinematic pressure $\tilde{p}$ from the manufactured solution and multiply it by a scalar $\sigma$. Thus the pressure term in the forcing function $\mathbf{f}$ will also be multiplied by $\sigma$, and the exact solution to which our numerical solution should converge has the same velocity as before but with a scaled kinematic pressure field. 

\begin{figure}
\centering
\subfloat[Velocity $L^2$ error]{\label{sfig:probust_lap_a}\includegraphics[width=.5\textwidth]{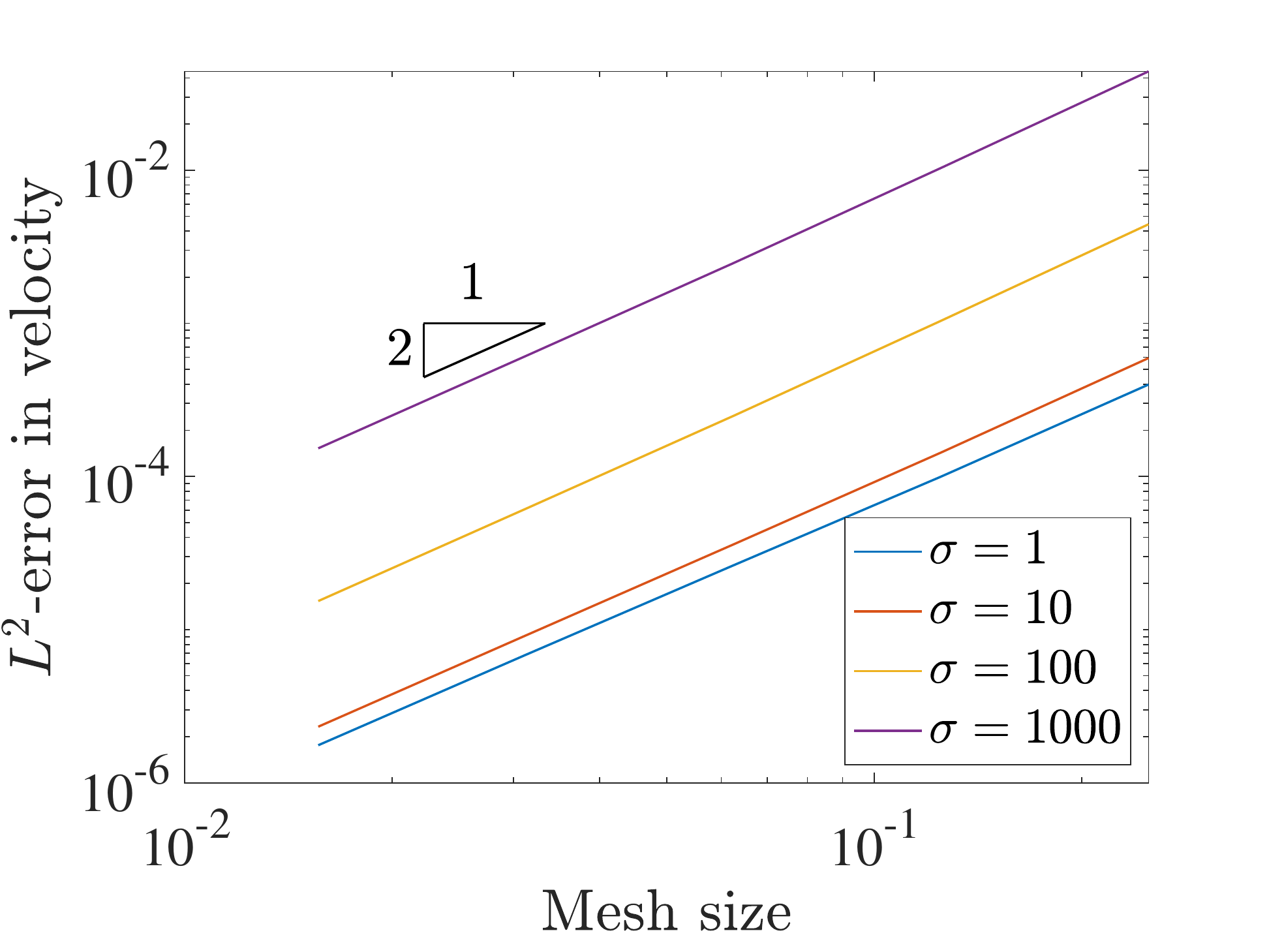}}\hfill
\subfloat[Velocity $H^1$ error]{\label{sfig:probust_lap_b}\includegraphics[width=.5\textwidth]{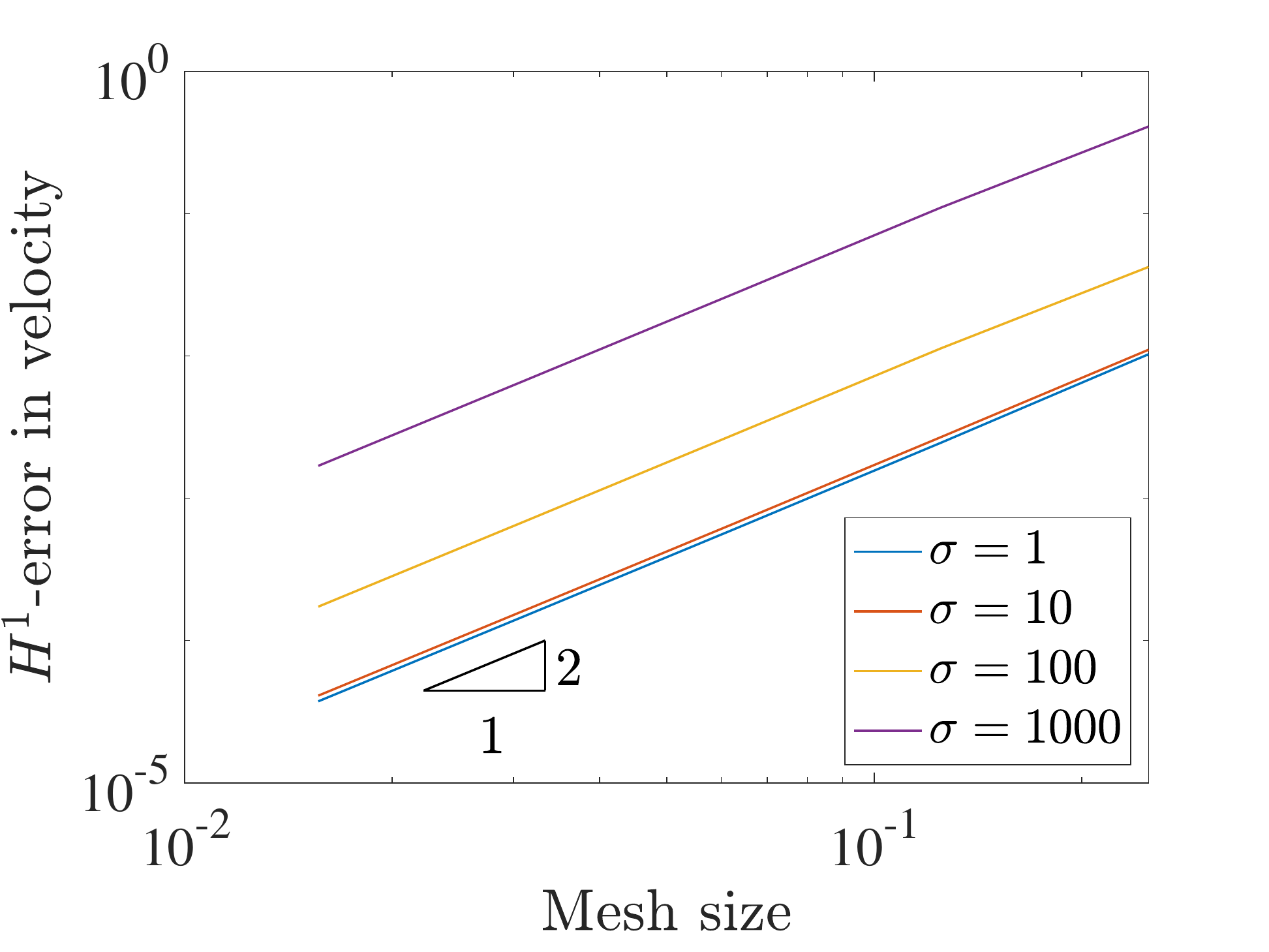}} \\
\caption{Errors in 2D manufactured vortex solution with varying pressure scaling for velocity-pressure formulation}
\label{fig:probust_lap}
\end{figure}

\begin{figure}
\centering
\subfloat[Velocity $L^2$ error]{\label{sfig:probust_a}\includegraphics[width=.5\textwidth]{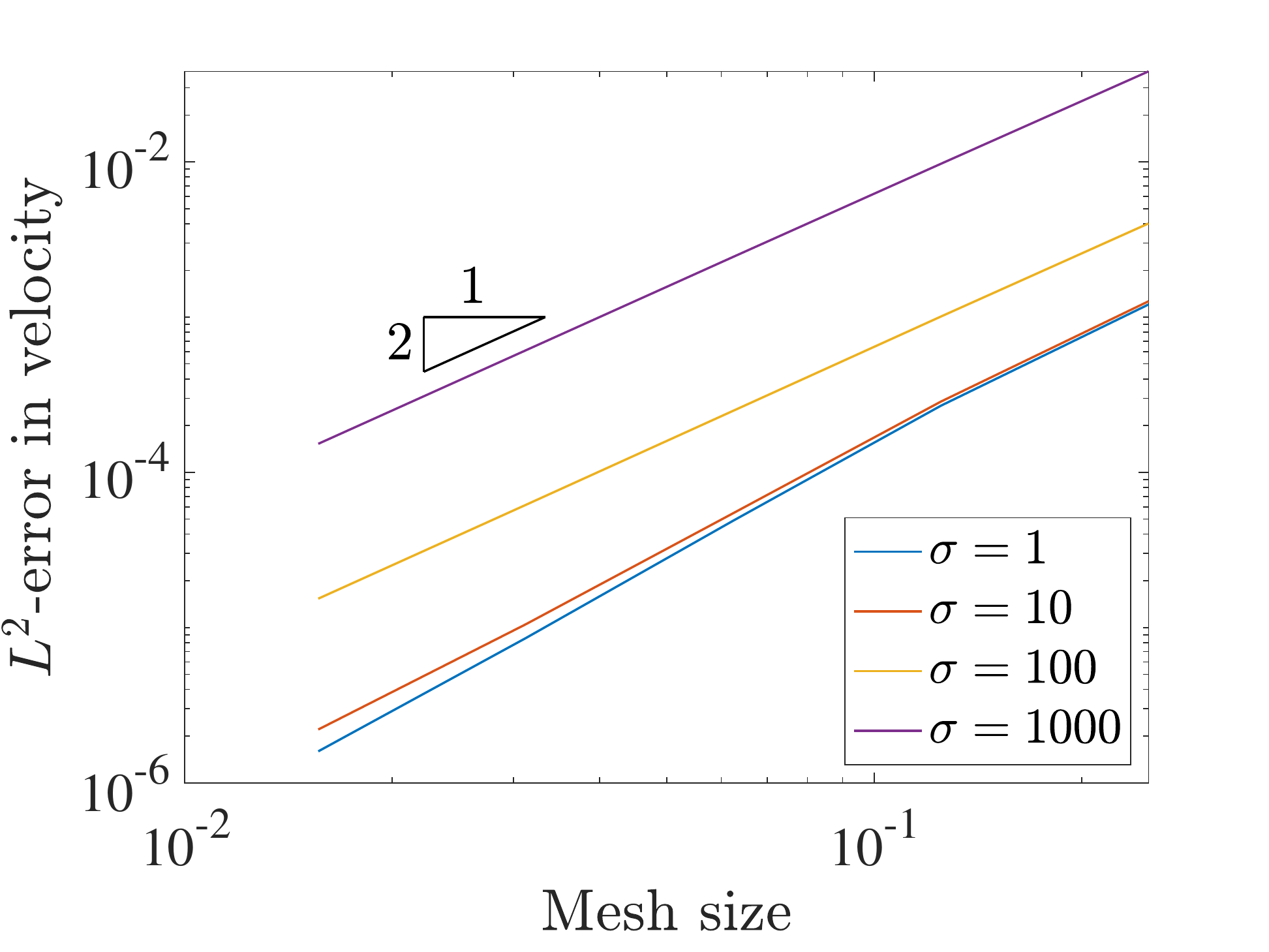}}\hfill
\subfloat[Velocity $H^1$ error]{\label{sfig:probust_b}\includegraphics[width=.5\textwidth]{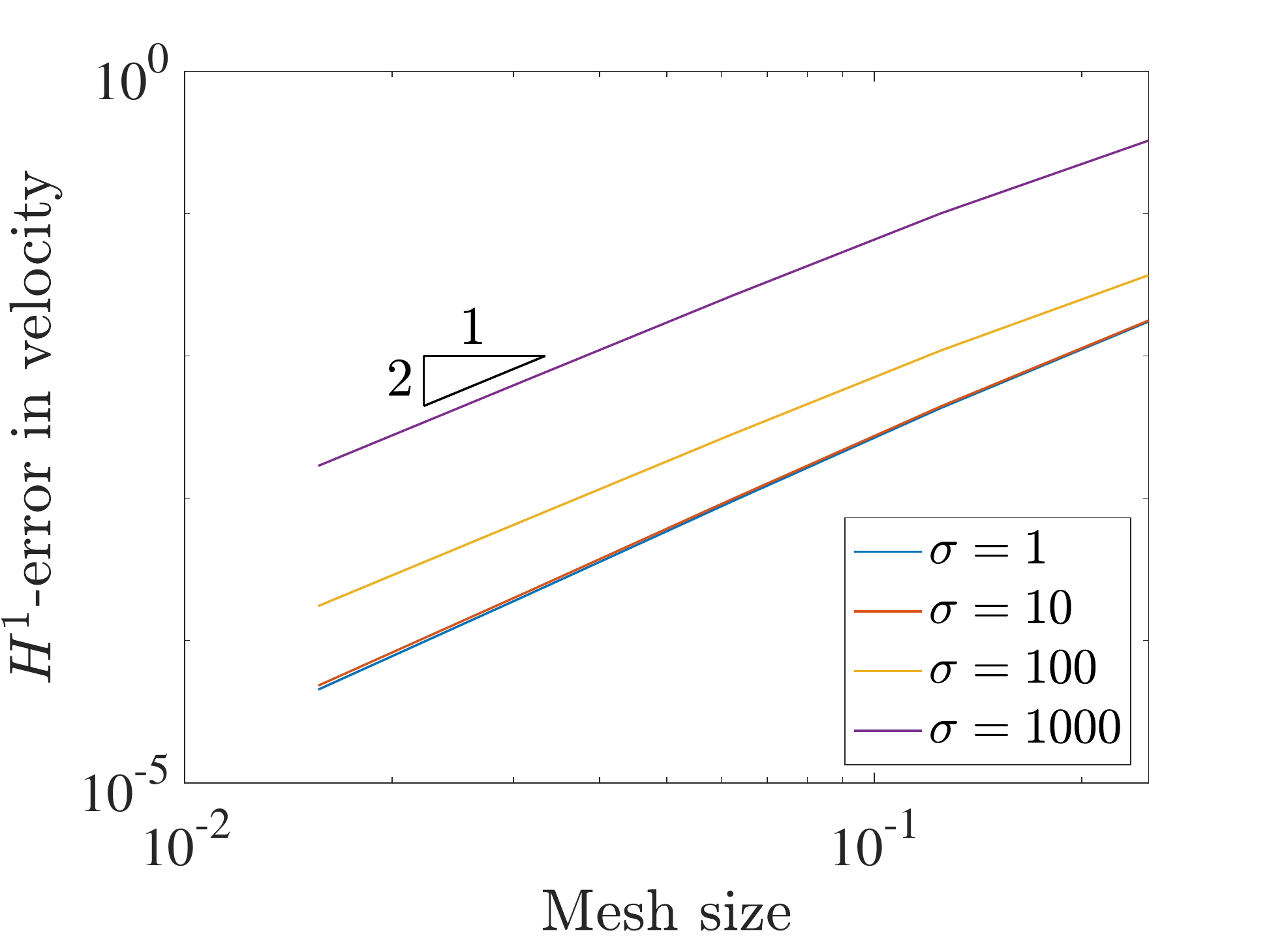}} \\
\caption{Errors in 2D manufactured vortex solution with varying pressure scaling for vorticity-velocity-pressure formulation}
\label{fig:probust}
\end{figure}

For a pressure robust method this increase in the pressure magnitude, and thus the pressure approximation errors, will not affect the velocity approximation error. Conversely, a non-pressure robust method will see its velocity errors increase as the pressure is scaled larger. Figure \ref{fig:probust_lap} shows the convergence of the velocity errors for the two field scheme with $k' = 2$ and increasing values of the scalar $\sigma$, while Figure \ref{fig:probust} shows the same for the three field formulation. Clearly the velocity error increases in both schemes as $\sigma$ increases, meaning the method is not pressure robust. This is interesting as the divergence-conforming Galerkin method upon which this work is based is pressure robust.

\subsection{Reynolds Robustness}

Similar to pressure robustness, we also want to test how the errors in the solution behave as the Reynolds number is increased. We increase the Reynolds number by decreasing the viscosity $\nu$. This affects the viscous term in the forcing vector $\mathbf{f}$, but the exact solution to the problem is identical to the original manufactured solution. 

\begin{figure}
\centering
\subfloat[Velocity $L^2$ error]{\label{sfig:rerobust_lap_a}\includegraphics[width=.5\textwidth]{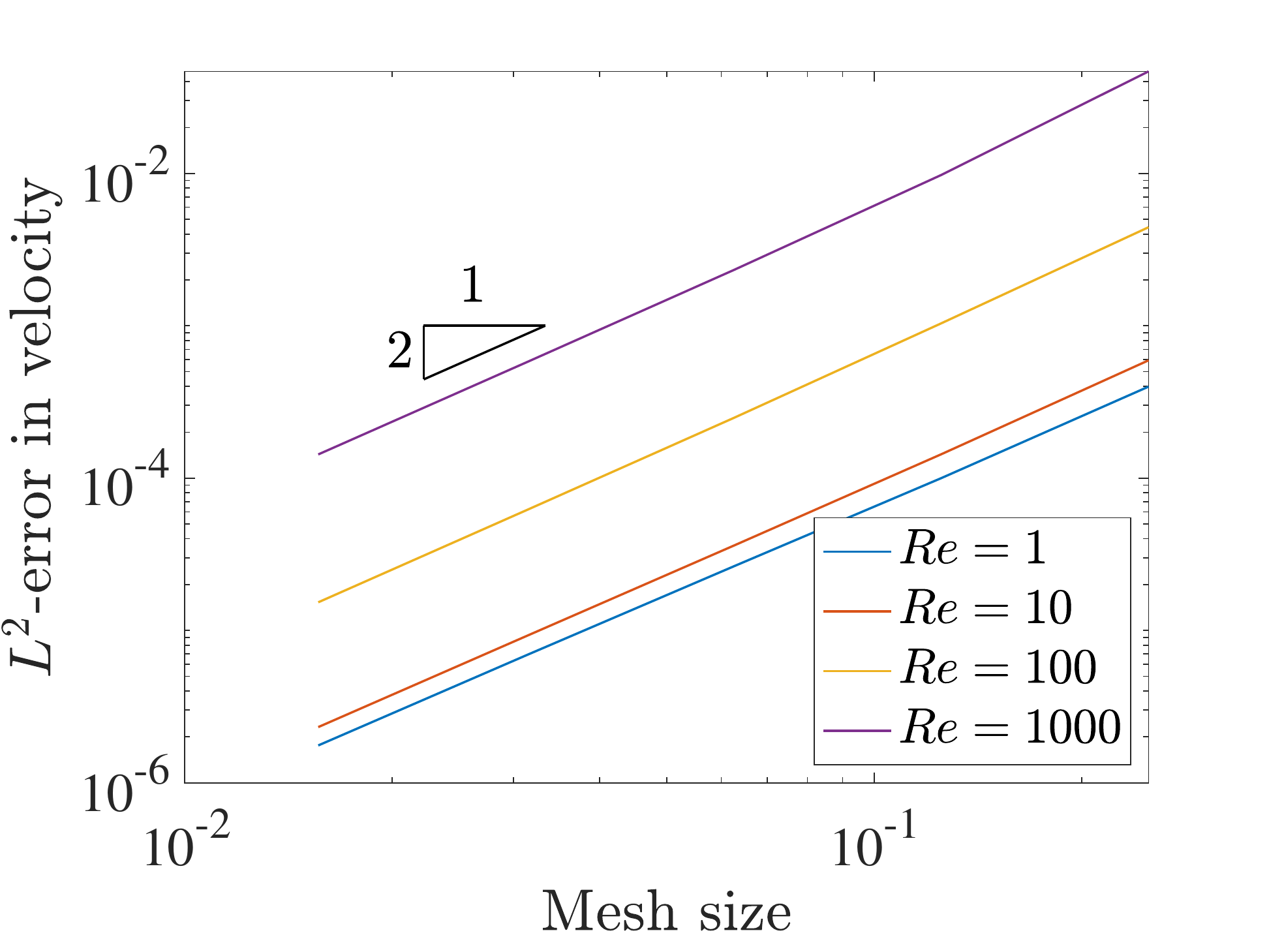}}\hfill
\subfloat[Velocity $H^1$ error]{\label{sfig:rerobust_lap_b}\includegraphics[width=.5\textwidth]{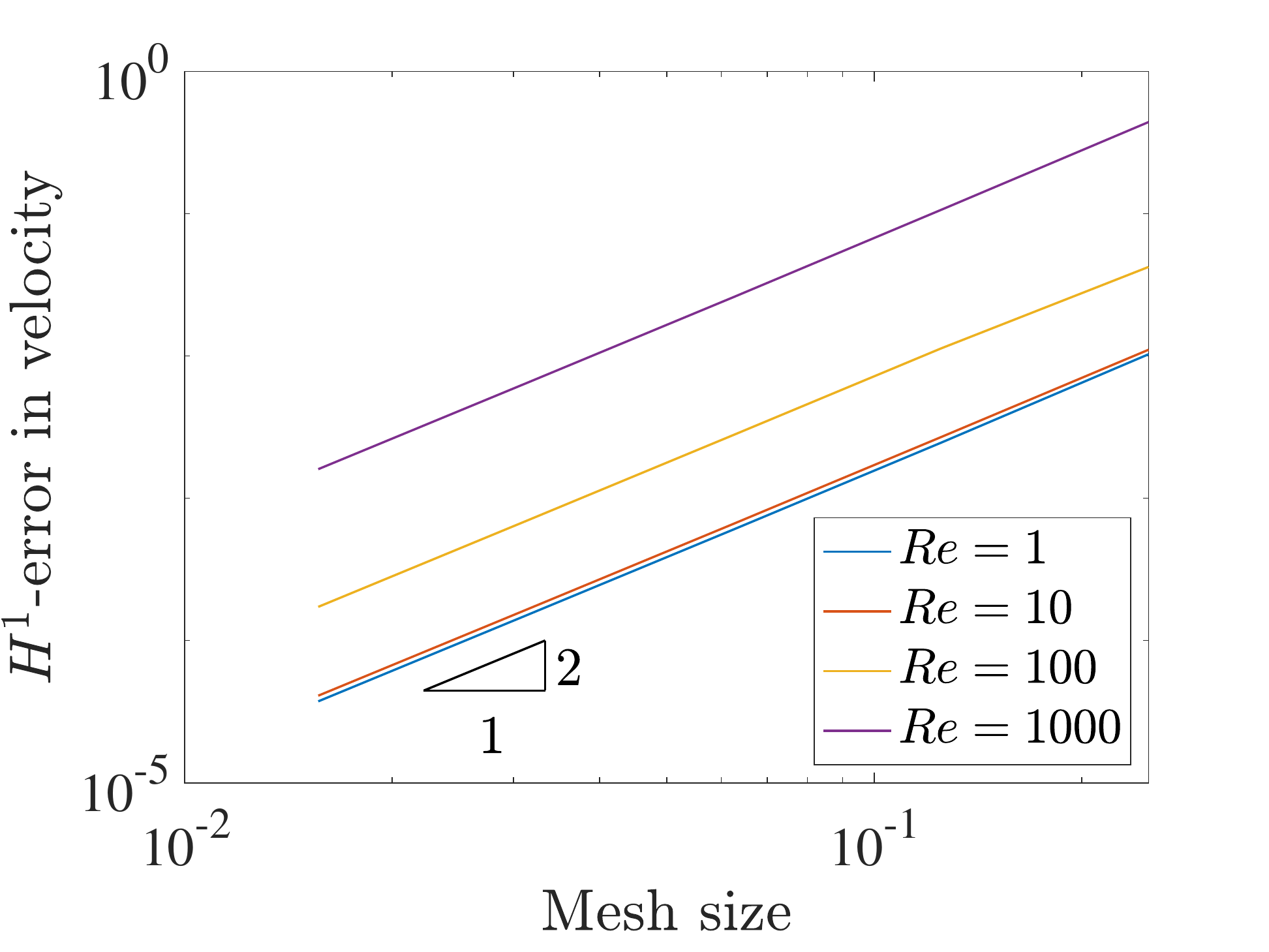}} \\
\caption{Errors in 2D manufactured vortex solution with varying Reynolds number for velocity-pressure formulation}
\label{fig:rerobust_lap}
\end{figure}

\begin{figure}
\centering
\subfloat[Velocity $L^2$ error]{\label{sfig:rerobust_a}\includegraphics[width=.5\textwidth]{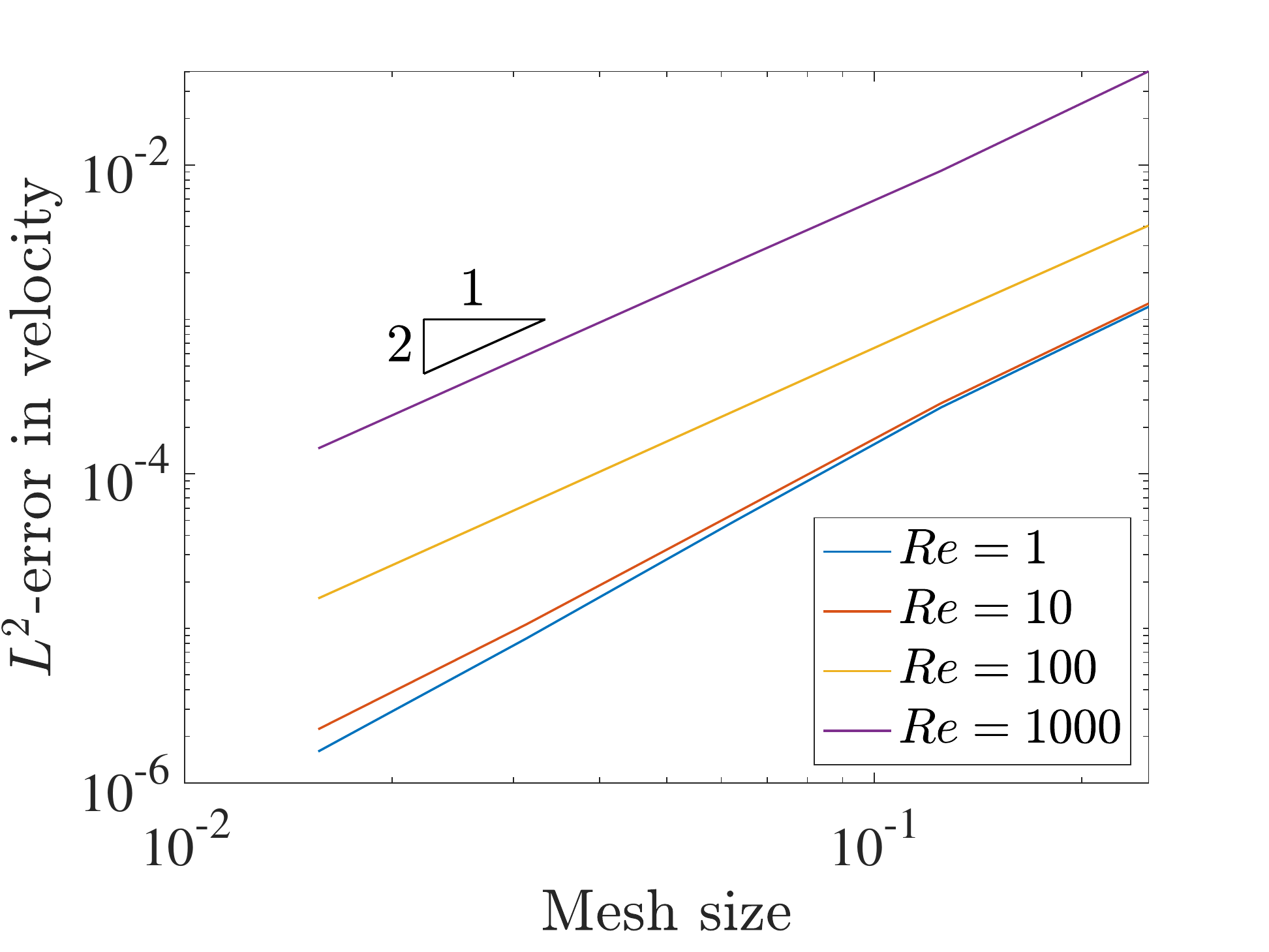}}\hfill
\subfloat[Velocity $H^1$ error]{\label{sfig:rerobust_b}\includegraphics[width=.5\textwidth]{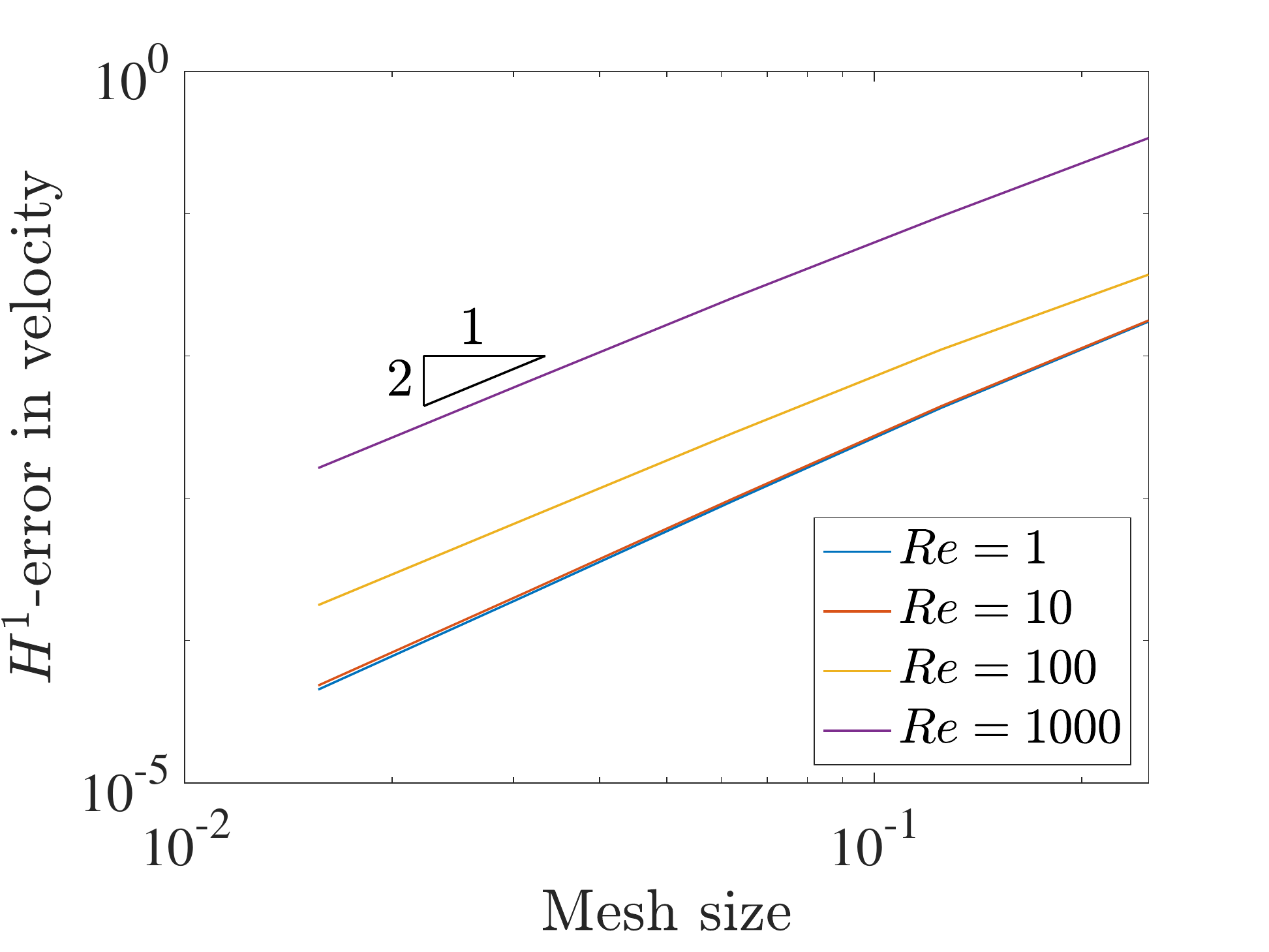}} \\
\caption{Errors in 2D manufactured vortex solution with varying Reynolds number for vorticity-velocity-pressure formulation}
\label{fig:rerobust}
\end{figure}

Figures \ref{fig:rerobust_lap} and \ref{fig:rerobust} detail the convergence of the velocity errors as the Reynolds number increases, again for $k' = 2$, in the two and three field schemes. Once again, the error in the velocity field increases as we increase the Reynolds number, in contrast to the divergence-conforming Galerkin setting, where the velocity error is agnostic to increasing Reynolds number \cite{evans_steady_NS}. 

\subsection{Two-Dimensional Lid-Driven Cavity Flow}

The next 2D numerical test problem that we consider is the square lid-driven cavity flow. The left, right, and bottom walls of the cavity remain fixed while the top wall slides in the positive $x$ direction, causing vortices to develop within the domain. Due to the inconsistency in the boundary conditions, pressure singularities exist in the corners of the domain, making this a challenging test case for a numerical scheme to properly capture.


For our simulations, we set both the speed of the top wall $U = 1$ and the wall lengths $H = 1$. The kinematic viscosity $\nu$ defines the Reynolds number through $Re = \frac{UH}{\nu} = \frac{1}{\nu}$. In particular, we consider the flows produced with $Re = 100$, $Re = 400$, and $Re = 1000$. To validate our results, we compare the centerline velocity profiles at each Reynolds number with the results from Ghia \textit{et al} \cite{ghia_cavity}.





\begin{figure}
\centering
\subfloat[$Re = 100$ velocities with $h = 1/32$]{\label{sfig:2D_cav_lap_a}\includegraphics[width=.5\textwidth]{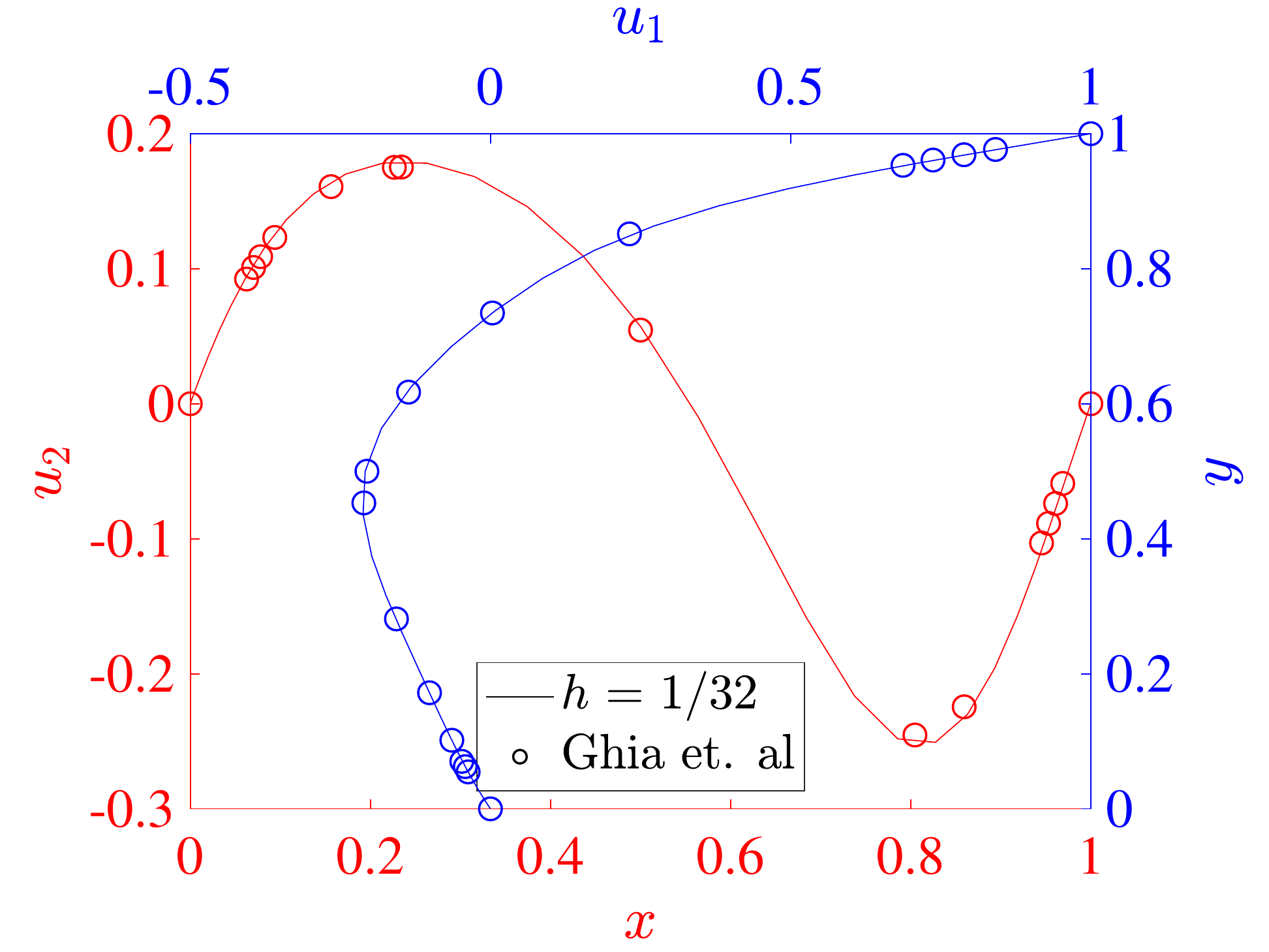}}\hfill
\subfloat[$Re = 100$ velocities with $h = 1/64$]{\label{sfig:2D_cav_lap_b}\includegraphics[width=.5\textwidth]{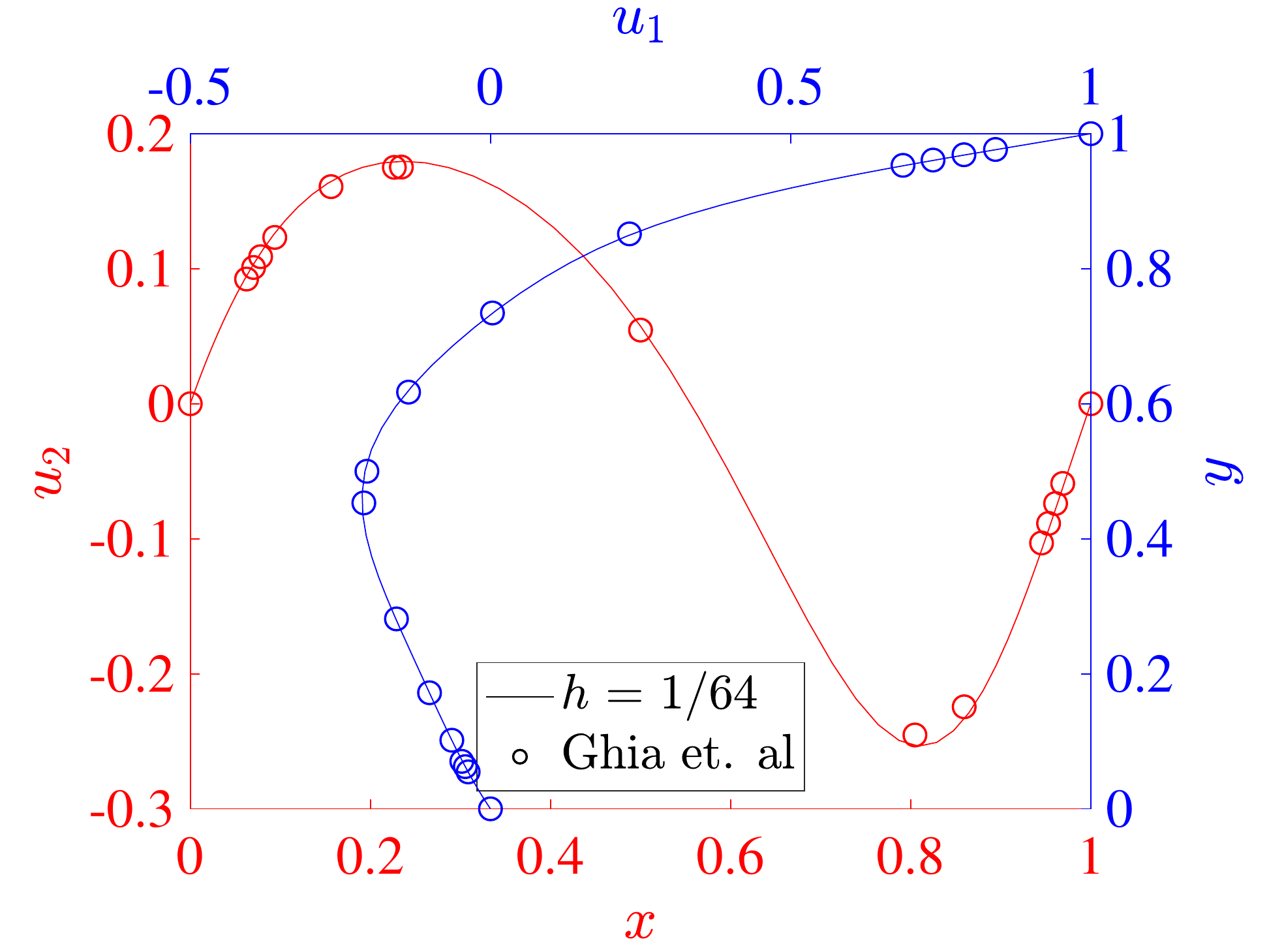}} \\
\subfloat[$Re = 400$ velocities with $h = 1/32$]{\label{sfig:2D_cav_lap_c}\includegraphics[width=.5\textwidth]{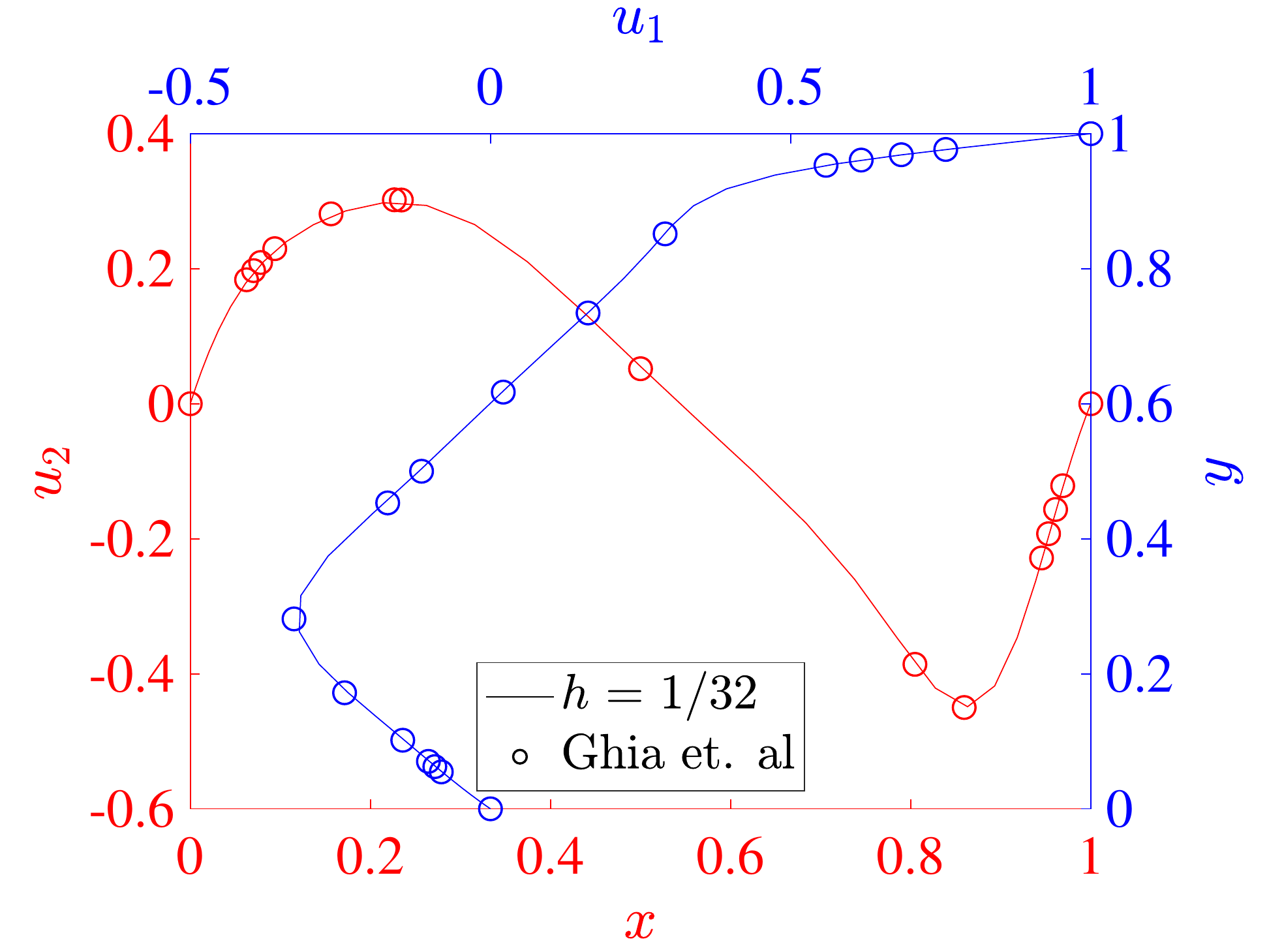}}\hfill
\subfloat[$Re = 400$ velocities with $h = 1/64$]{\label{sfig:2D_cav_lap_d}\includegraphics[width=.5\textwidth]{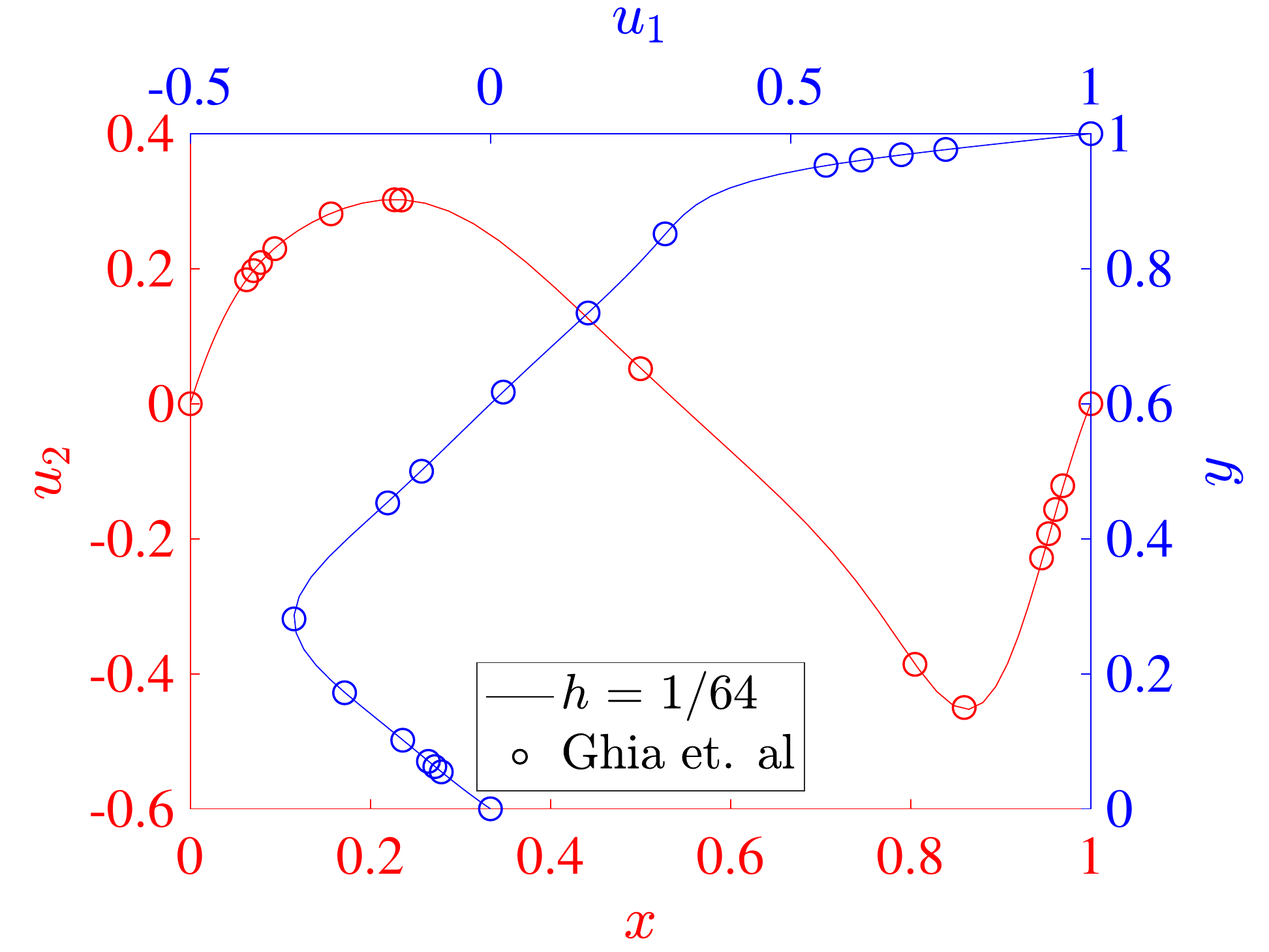}} \\
\subfloat[$Re = 1000$ velocities with $h = 1/32$]{\label{sfig:2D_cav_lap_e}\includegraphics[width=.5\textwidth]{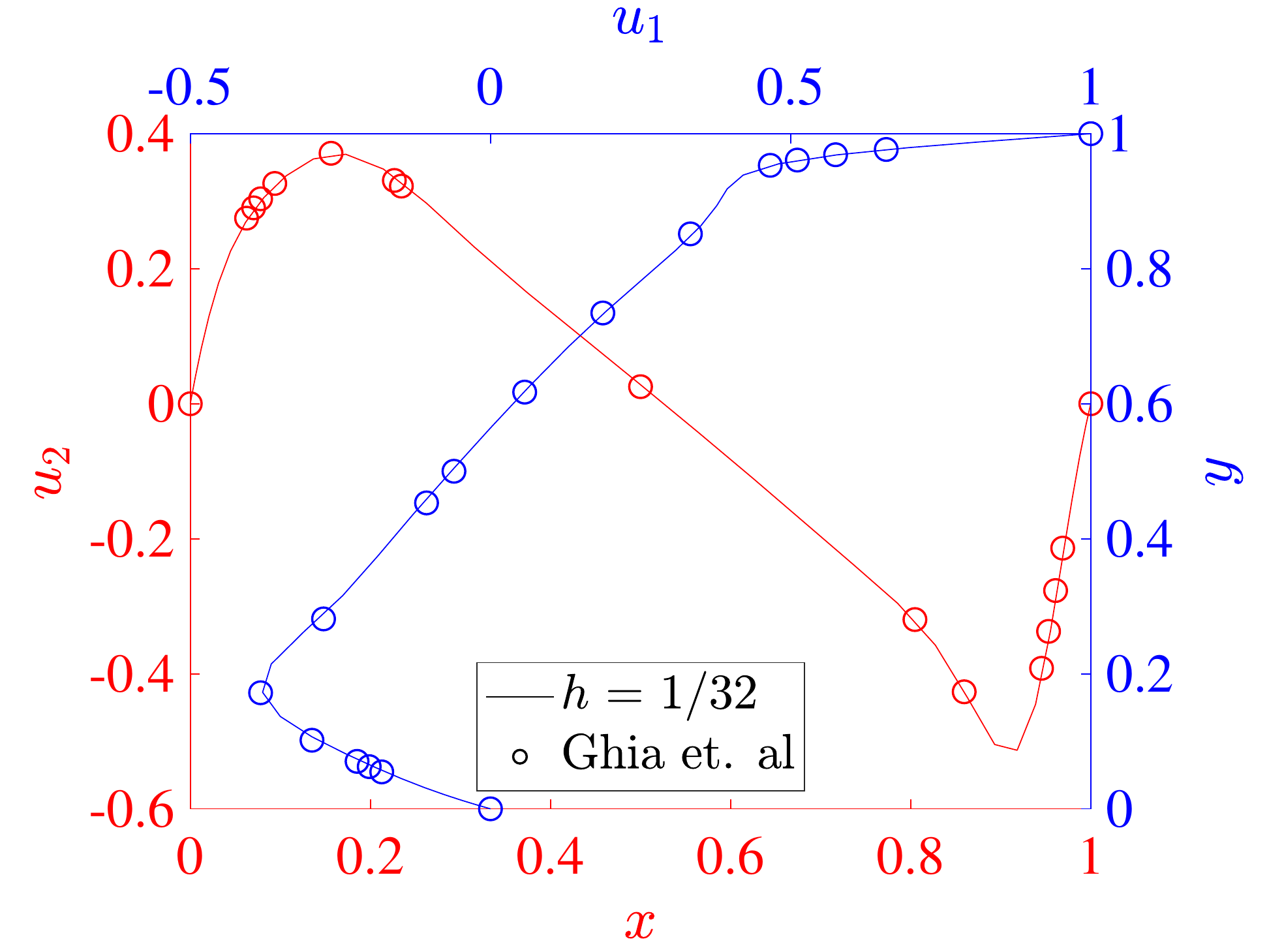}}\hfill
\subfloat[$Re = 1000$ velocities with $h = 1/64$]{\label{sfig:2D_cav_lap_f}\includegraphics[width=.5\textwidth]{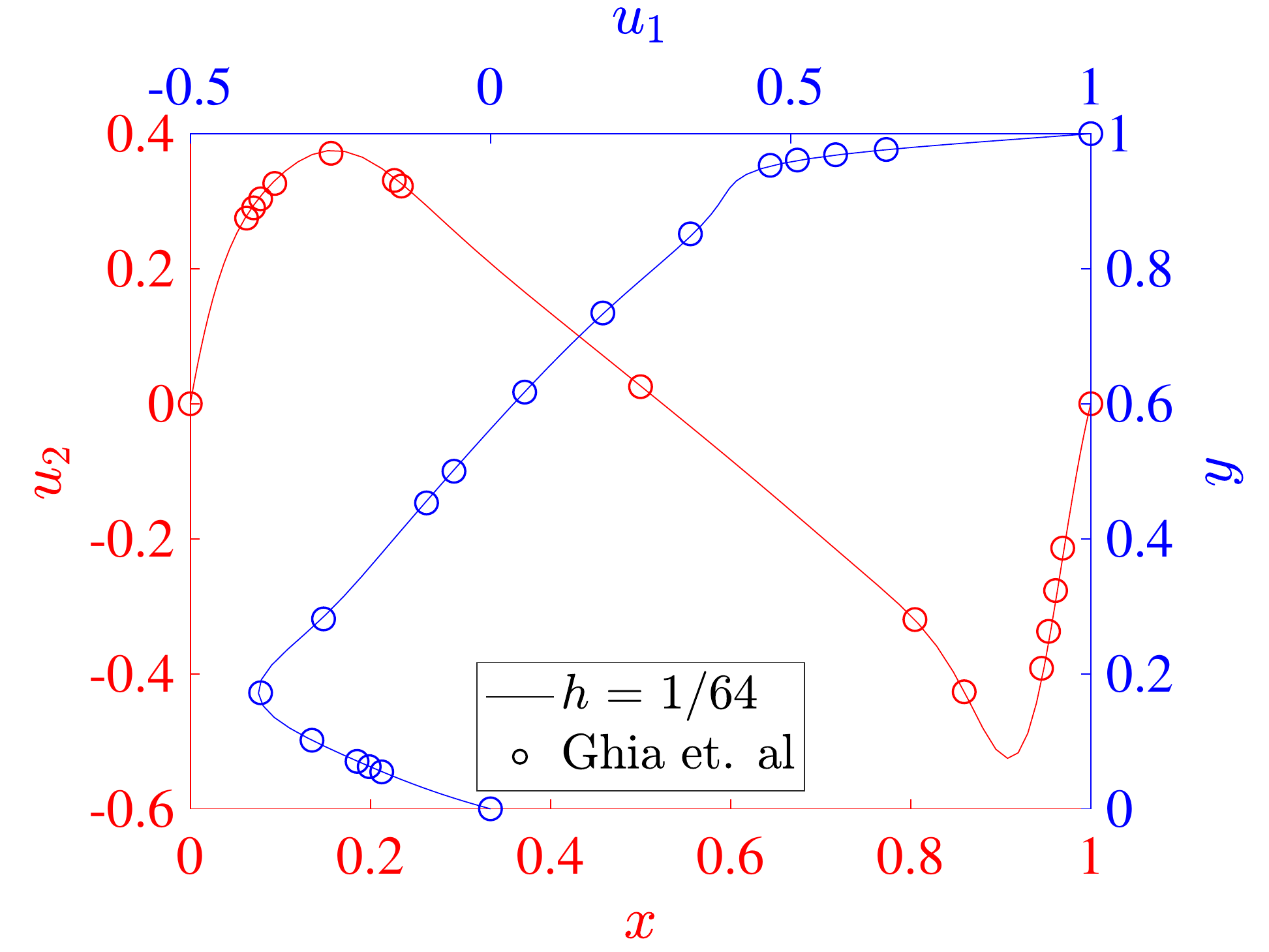}} \\
\caption{Centerline velocity profiles for 2D lid-driven cavity with velocity-pressure formulation, $k' = 2$. Red curves and axes represent the vertical velocity along the horizontal centerline, while blue curves and axes represent the horizontal velocity along the vertical centerline.}
\label{fig:2D_cav_lap}
\end{figure}

\begin{figure}
\centering
\subfloat[$Re = 100$ velocities with $h = 1/32$]{\label{sfig:2D_cav_a}\includegraphics[width=.5\textwidth]{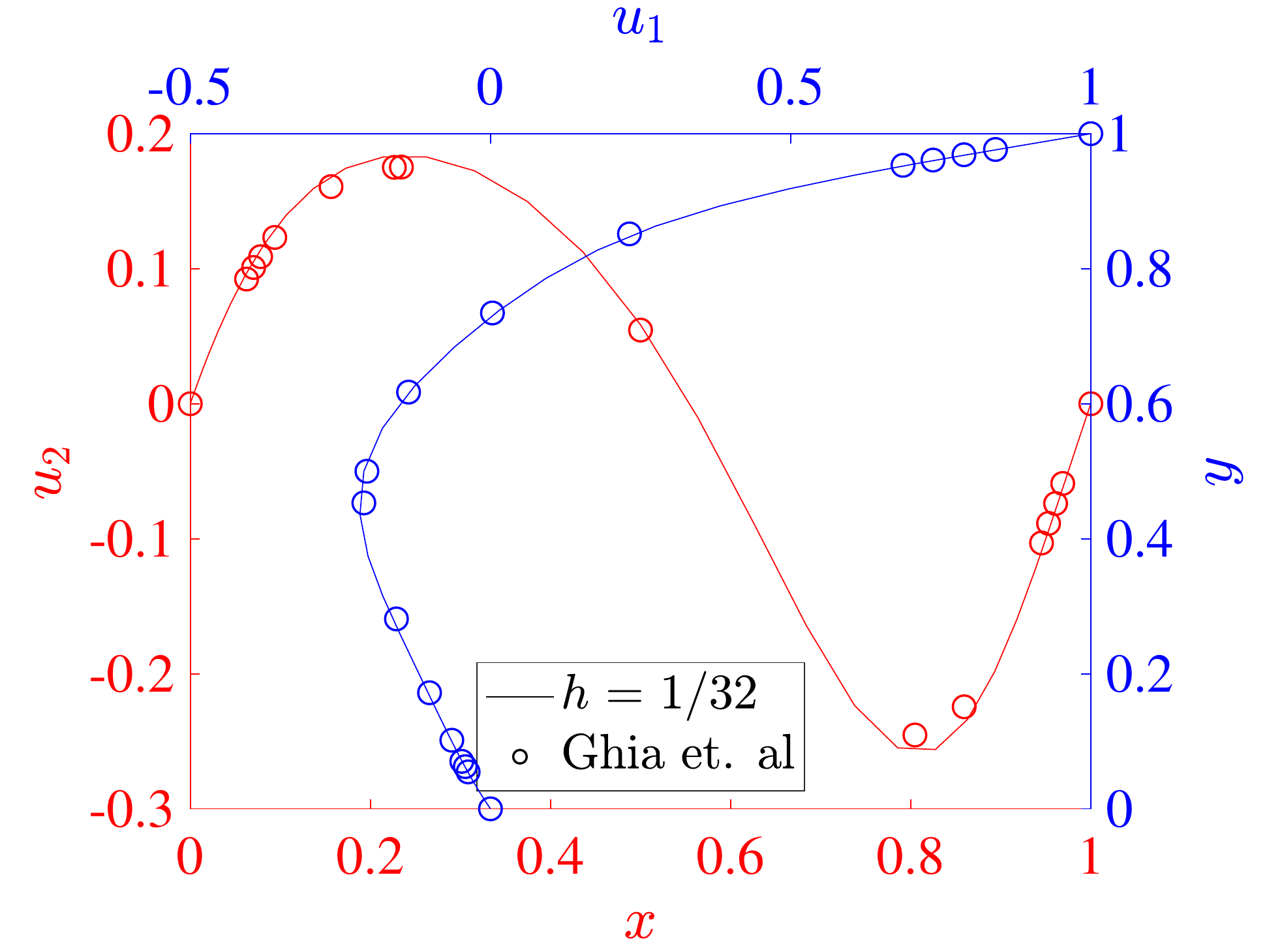}}\hfill
\subfloat[$Re = 100$ velocities with $h = 1/64$]{\label{sfig:2D_cav_b}\includegraphics[width=.5\textwidth]{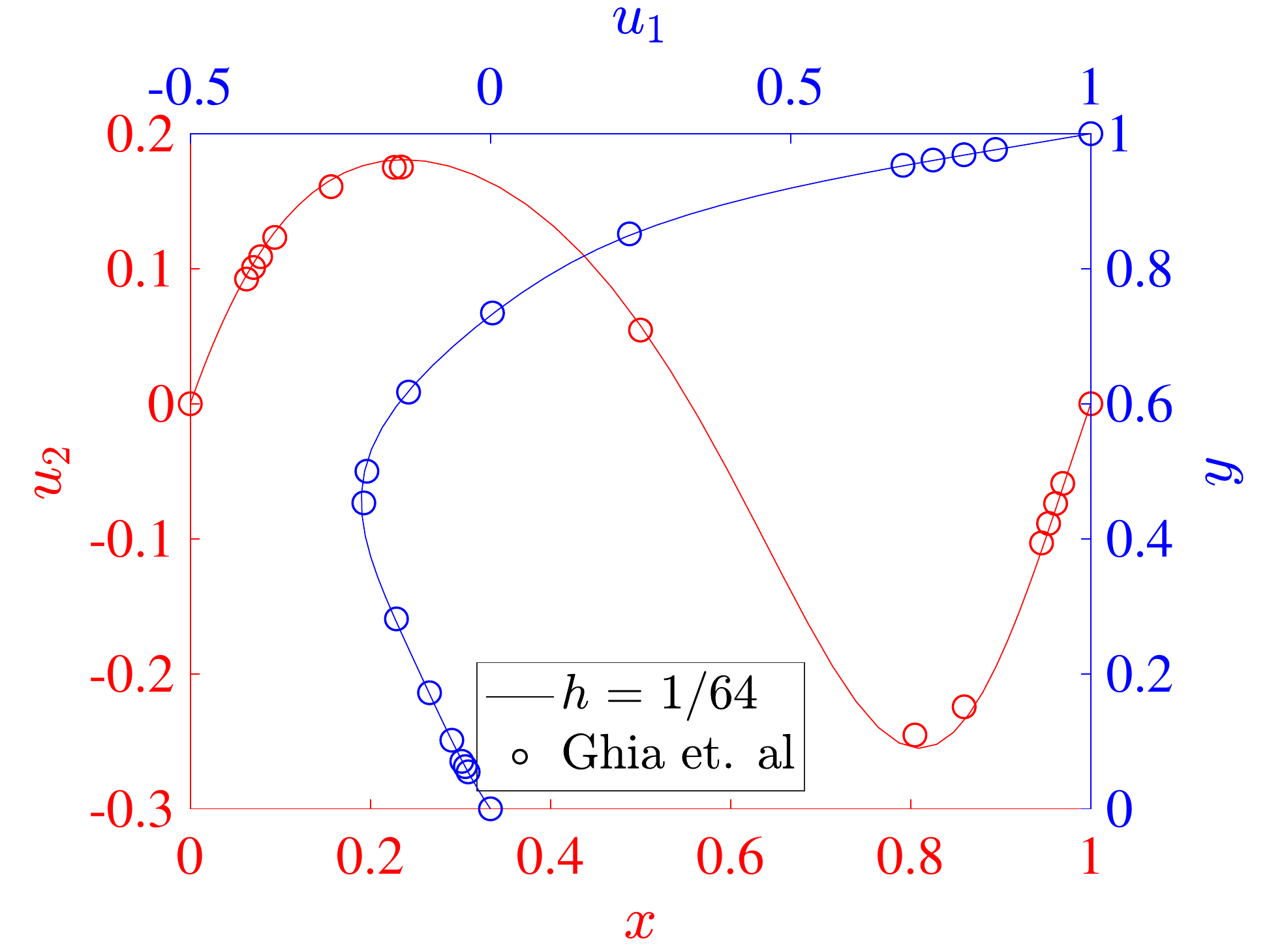}} \\
\subfloat[$Re = 400$ velocities with $h = 1/32$]{\label{sfig:2D_cav_c}\includegraphics[width=.5\textwidth]{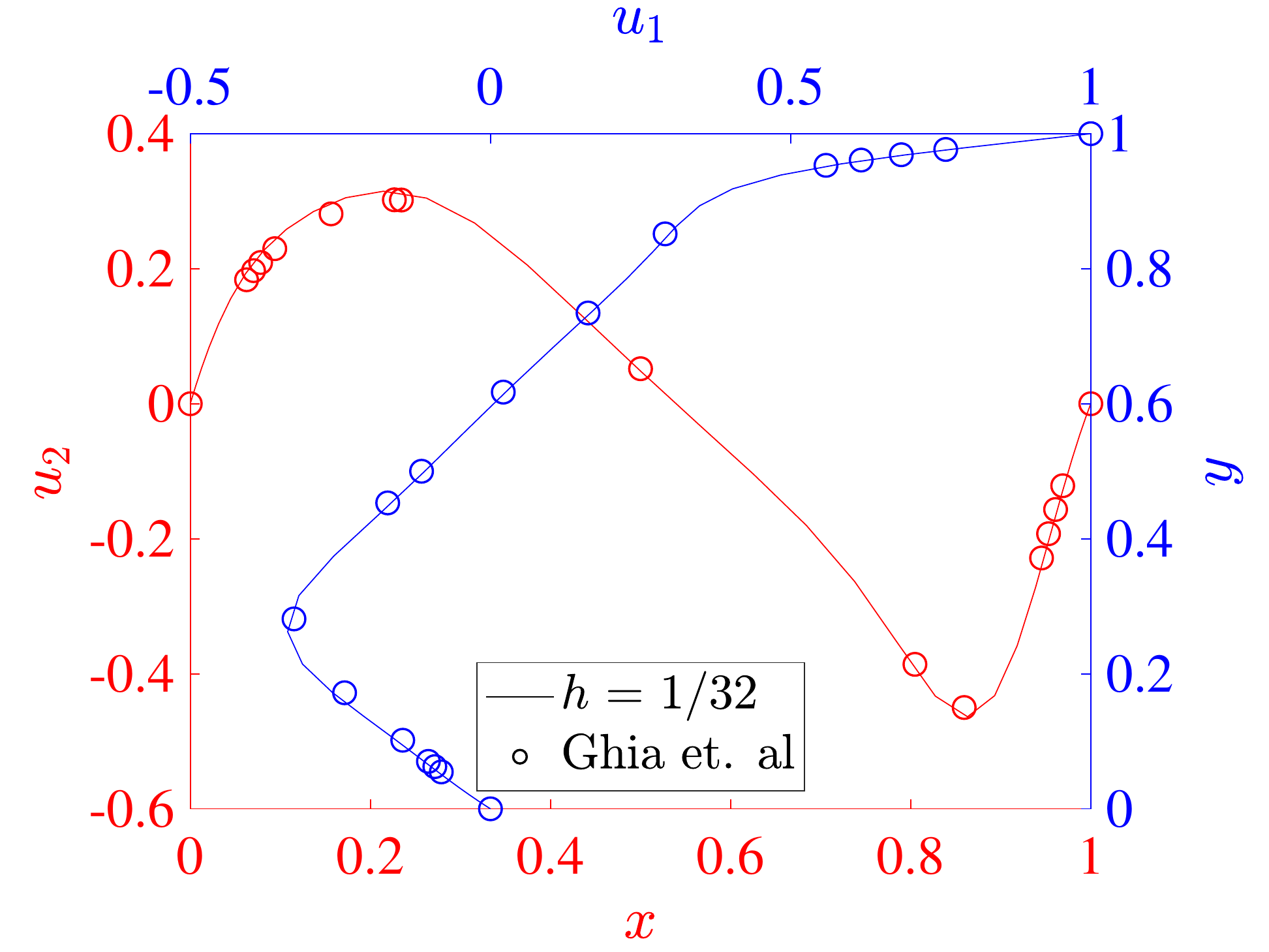}}\hfill
\subfloat[$Re = 400$ velocities with $h = 1/64$]{\label{sfig:2D_cav_d}\includegraphics[width=.5\textwidth]{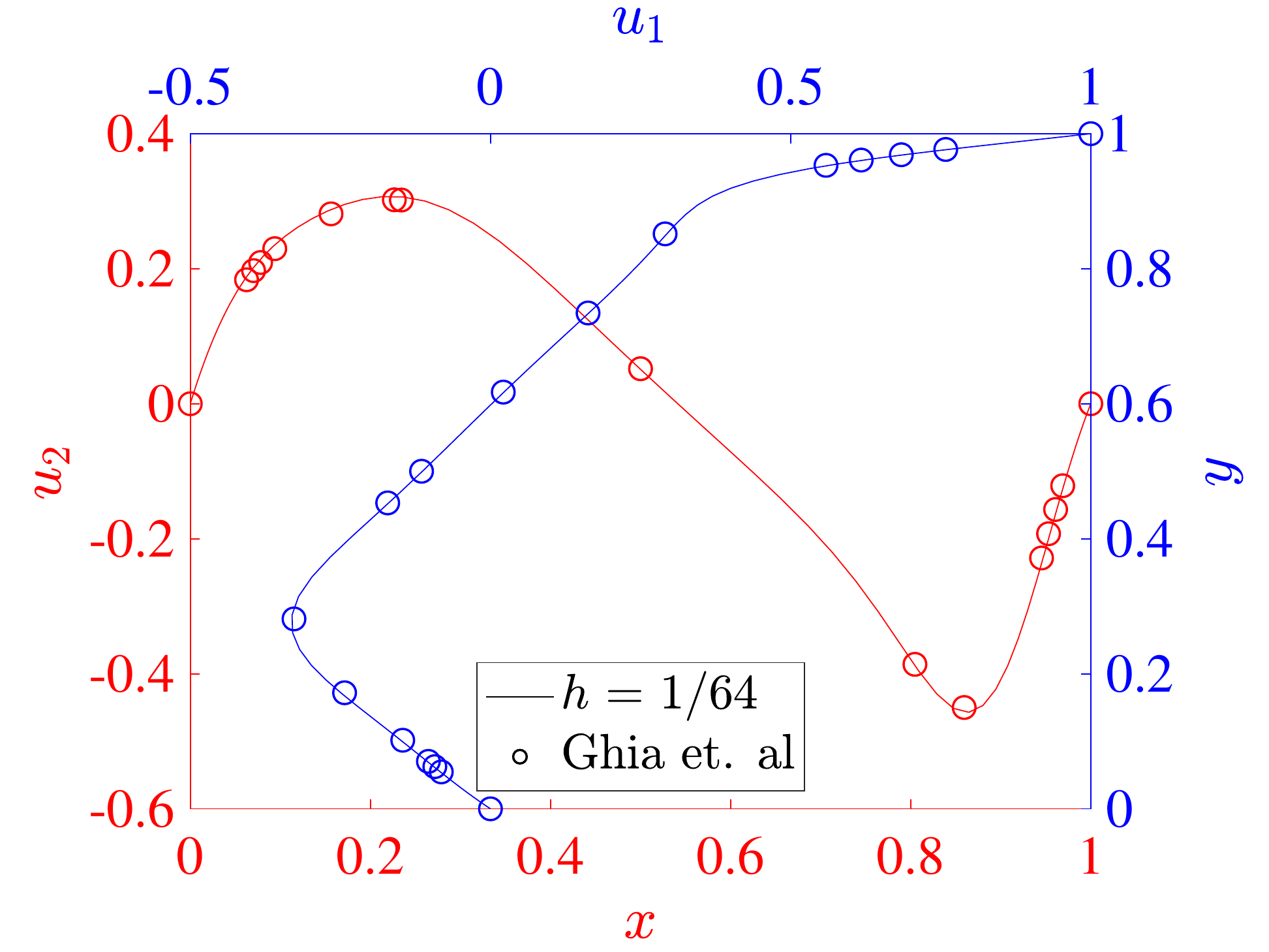}} \\
\subfloat[$Re = 1000$ velocities with $h = 1/32$]{\label{sfig:2D_cav_e}\includegraphics[width=.5\textwidth]{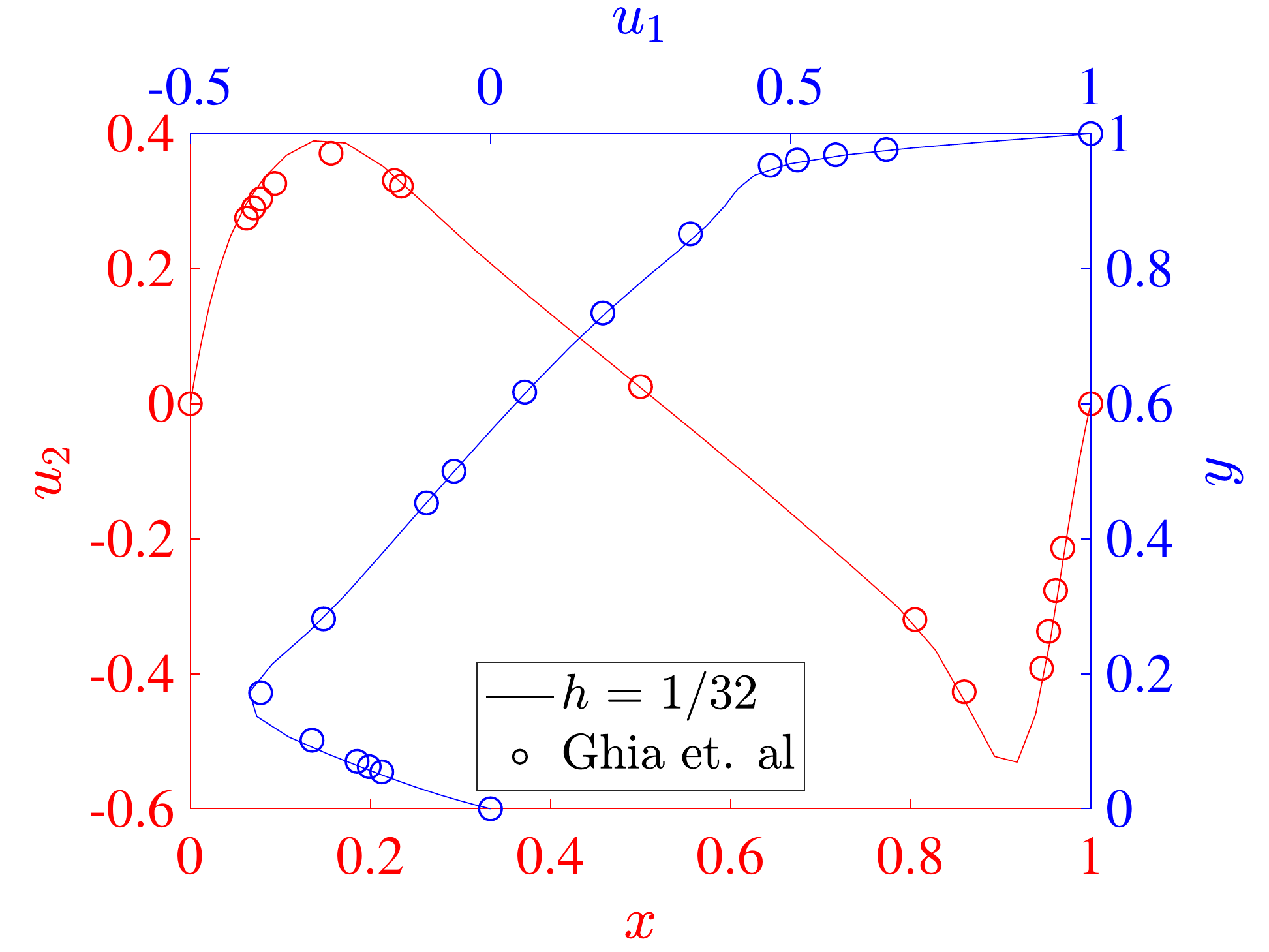}}\hfill
\subfloat[$Re = 1000$ velocities with $h = 1/64$]{\label{sfig:2D_cav_f}\includegraphics[width=.5\textwidth]{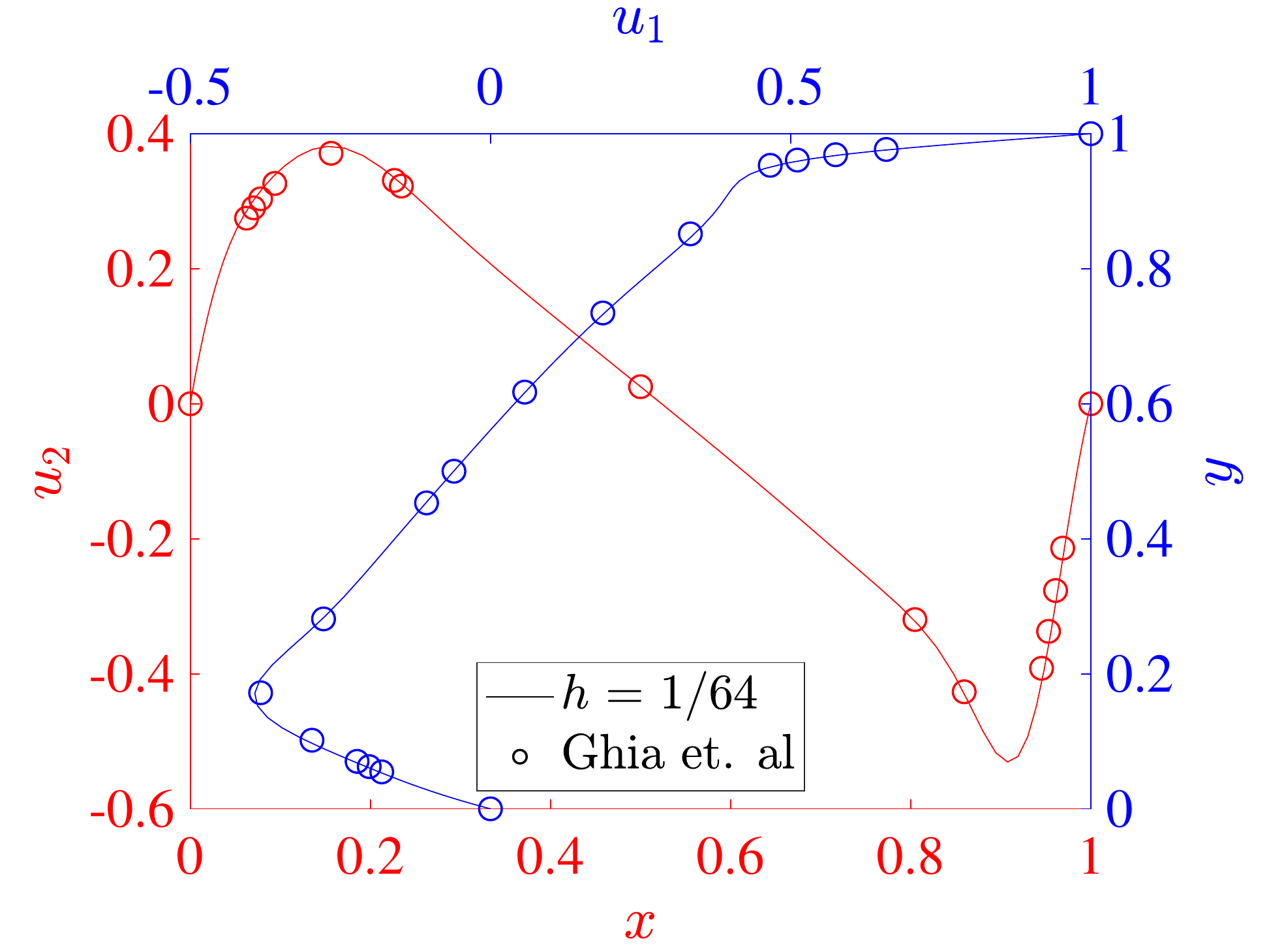}} \\
\caption{Centerline velocity profiles for 2D lid-driven cavity with vorticity-velocity-pressure formulation, $k' = 2$. Red curves and axes represent the vertical velocity along the horizontal centerline, while blue curves and axes represent the horizontal velocity along the vertical centerline.}
\label{fig:2D_cav}
\end{figure}

 Figure \ref{fig:2D_cav_lap} details the two field formulation results across the three considered Reynolds numbers and two mesh sizes: a 32 element stretched mesh and a 64 element stretched mesh. The stretched mesh is formed by setting the interior knots of the knot vectors defining the bases in each direction as 

 \begin{equation}
     \xi_i = \frac{1}{2} \left( 1 + \frac{\tanh(4ih - 2)}{\tanh(2)} \right) \quad \forall \xi_i \in \Xi, 
 \end{equation}

 \noindent where $h$ is the mesh size in each direction. Figure \ref{fig:2D_cav} shows the same results for the three field formulation. The collocation results from both schemes agree very well with the reference data in all cases, and we see that the results are converging with increasing resolution. At a Reynolds number of 100, all of our results show that the maximum and minimum values of the vertical velocity are larger in magnitude than those of Ghia \textit{et al}. This is similar to the behavior seen in the Galerkin method \cite{evans_steady_NS}, and we note that there are some inaccuracies in the Ghia data for this low Reynolds number case \cite{evans_steady_NS, botella1998benchmark}. For Reynolds number 400, the two field formulation predicts extrema in velocity that are slightly smaller than the three field predictions, which match the corresponding Galerkin results very well. This trend is also valid at a Reynolds number of 1000. Moreover, while we have used stretched meshes here, the results with a non-stretched mesh are similar.

As a more quantitative comparison, we compute the minimum horizontal velocity along the vertical centerline as well as the maximum and minimum vertical velocities along the horizontal centerline for each of simulations presented above. These results are shown for a Reynolds number of  100 in Table \ref{table:cavity_convergence}, along with the values from \cite{ghia_cavity} and \cite{botella_collocation}. These results show the inadequacy of the Ghia results at this Reynolds number, and for the most part the $k' = 2$ collocation results outperform the Ghia data when compared to the pseudospectral results. To highlight the potential possibilities of the collocation methods, we also compute results using an unstretched mesh of 8 elements in each direction and $k' = 20$ for both the two and three field formulations. While this would be essentially infeasible with a Galerkin method, as the quadrature would be prohibitively expensive, it is handled with ease by the collocation schemes. We see that these results match the pseudospectral results, even on the utilized coarse meshes. 


\begin{table}
\centering
\caption{Velocity extrema for 2D lid-driven cavity at $Re = 100$\label{table:cavity_convergence}}
\begin{tabular}{| c | c c c |}
\hline
Method & $u_{x, min}$ & $u_{y, max}$ & $u_{y, min}$\\
\hline
Collocation, 2 field formulation, $k' = 2$ and $h = 1/32$ & $-0.21348$ & $0.17941$ & $-0.25307$ \\
Collocation, 2 field formulation, $k' = 2$ and $h = 1/64$ & $-0.21389$ & $0.17953$ & $-0.25358$ \\
Collocation, 2 field formulation, $k' = 20$ and $h = 1/8$ & $-0.21404$ & $0.17957$ & $-0.25380$ \\
Collocation, 3 field formulation, $k' = 2$ and $h = 1/32$ & $-0.21800$ & $0.18392$ & $-0.25908$ \\
Collocation, 3 field formulation, $k' = 2$ and $h = 1/64$ & $-0.21511$ & $0.18075$ & $-0.25521$ \\
Collocation, 3 field formulation, $k' = 20$ and $h = 1/8$ & $-0.21404$ & $0.17957$ & $-0.25380$ \\
Pseudospectral (Ref. \cite{botella1998benchmark}) & $-0.21404$ & $0.17957$ & $-0.25380$ \\
Ghia \textit{et al.} (Ref. \cite{ghia_cavity}) & $-0.21090$ & $0.17527$ & $-0.24533$ \\
\hline
\end{tabular}
\end{table}

\section{Collocation Methods on Cubic Domains}

The previous two sections detailed the construction of the divergence-conforming collocation methods in 2D and tested their behavior numerically. In the following, we will highlight the required modifications to the methods to solve problems in 3D cubic domains.  

\subsection{Review of Galerkin Methods}










Similar to 2D, we start by reviewing the form of the divergence-conforming isogeometric Galerkin methods for 3D problems. Again assume the velocity is subject to  Dirichlet boundary conditions along the entire boundary. We then define the discrete test and trial function spaces for velocity as $\boldsymbol{\mathcal{V}}_{h,0}$ and $\boldsymbol{\mathcal{V}}_{h,\mathbf{g}}$, which are defined as the same $\boldsymbol{\mathcal{V}}_h$ from Equation \eqref{eq:3d_velocity_space} with either no penetration boundary conditions strongly enforced (for the test space $\boldsymbol{\mathcal{V}}_{h,0}$) or with the normal velocity prescribed at collocation points as given by the boundary data $\mathbf{g}$ (for the trial space $\boldsymbol{\mathcal{V}}_{h,\mathbf{g}}$). We also define the test and trial space for pressure as $\mathcal{Q}_{h,0}$, where $\mathcal{Q}_h$ is the same space as in Equation \eqref{eq:3d_pressure_space} but with the added condition that the pressure must have zero integral. Then the Galerkin formulation for the velocity-pressure form would read

\bigskip

$$
\left\{ \hspace{5pt}
\parbox{6in}{
\noindent Given $\nu \in \mathbb{R}^+$, $\textbf{f} : \Omega \rightarrow \mathbb{R}^3$, and $\textbf{g} : \partial \Omega \rightarrow \mathbb{R}^3$, find $\mathbf{u}^h \in \boldsymbol{\mathcal{V}}_{h,\mathbf{g}}$ and $p^h \in \mathcal{Q}_{h,0}$ such that, $\forall (\mathbf{w}^h, q^h) \in (\boldsymbol{\mathcal{V}}_{h,0}, \mathcal{Q}_{h,0})$:

\begin{equation}
\begin{split}
    &\int_{\Omega} (\nu \nabla \mathbf{w}^h \cdot \nabla \mathbf{u}^h + \mathbf{w}^h \cdot (\mathbf{u}^h \cdot \nabla \mathbf{u}^h) - p^h \nabla \cdot \mathbf{w}^h) d\Omega \\& - \nu \int_{\partial \Omega} (\nabla \mathbf{u}^h \cdot \mathbf{n}) \cdot \mathbf{w}^h - \frac{C_{pen}}{h} \mathbf{u}^h \cdot \mathbf{w}^h d\mathbf{A} = \int_{\Omega} \mathbf{w}^h \cdot \mathbf{f} d \Omega + \nu \int_{\partial \Omega} \frac{C_{pen}}{h} \mathbf{g} \cdot \mathbf{w}^h d\mathbf{A} \label{eq:3d_weak_nit}
\end{split}
\end{equation}
\begin{equation}
    \int_{\Omega} q^h (\nabla \cdot \mathbf{u}^h) d \Omega = 0.
\end{equation}
}
\right.
$$

\bigskip

This weak form is essentially unchanged from the 2D case, with the only major difference being that the velocity has 3 components. The vorticity-velocity-pressure Galerkin form, however, is more different. In this case the discrete problem reads

\bigskip

$$
\left\{ \hspace{5pt}
\parbox{6in}{
Given $\nu \in \mathbb{R}^+$, $\textbf{f} : \Omega \rightarrow \mathbb{R}^3$, and $\textbf{g} : \partial \Omega \rightarrow \mathbb{R}^3$, find $\mathbf{u}^h \in \boldsymbol{\mathcal{V}}_{h,\mathbf{g}}$,  $P^h \in \mathcal{Q}_{h,0}$, and $\boldsymbol{\omega}^h \in \boldsymbol{\Psi}_h$ such that, $\forall (\mathbf{w}^h, q^h, \boldsymbol{\psi}^h) \in (\boldsymbol{\mathcal{V}}_{h,0}, Q_{h,0}, \boldsymbol{\Psi}_h)$:

\begin{equation}
    \int_\Omega (\nu \nabla \times \boldsymbol{\omega}^h) \cdot \mathbf{v}^h d\Omega + \int_\Omega (\boldsymbol{\omega}^h \times \textbf{u}^h) \cdot \mathbf{v}^h d\Omega - \int_\Omega P^h(\nabla \cdot \mathbf{v}^h) d\Omega = \int_\Omega \mathbf{f} \cdot \mathbf{v}^h d \Omega
\end{equation}

\begin{equation}
    \int_\Omega (\nabla \cdot \mathbf{u}^h) q^h d\Omega = 0 
\end{equation}

\begin{equation}
    \int_\Omega (\boldsymbol{\omega}^h \cdot \boldsymbol{\psi}^h) d \Omega + \int_\Omega \mathbf{u}^h \cdot (\nabla \times \boldsymbol{\psi}^h) d\Omega - \int_{\partial \Omega} (\boldsymbol{\psi}^h \times \mathbf{g}) \cdot \mathbf{n} dA = 0 .
    \label{eq:3D_weak}
\end{equation}

}
\right.
$$

\bigskip

\noindent Again the no-slip velocity boundary conditions appear as natural boundary conditions in the weak form of the constitutive equation. 

Within the collocation schemes, the unknowns are selected to reside in the same spaces as the corresponding Galerkin scheme, as in 2D. In the following we highlight the main changes to the method for 3D problems with regards to the choice of collocation grids and boundary condition enforcement before again summarizing the final form of the discrete equations. 

\subsection{Collocation Grids}

Much like the two dimensional case, in 3D we choose to collocate at Greville abscissae and the grids are different for each of the governing equations. For both formulations the schemes for the momentum and pressure equations are essentially unchanged; each momentum equation component is collocated at the Greville abscissae of the corresponding discrete velocity component space, and the continuity equation is collocated at the Greville abscissae of the discrete pressure space. Thus the velocity-pressure formulation extends fairly trivially to 3D. 

The constitutive equation in the vorticity-velocity-pressure formulation, on the other hand, is now split into components much like how the momentum equations are treated. We choose to collocate the $x$ component of the constitutive equation at the Greville abscissae associated with the discrete $x$ vorticity space ($S^{k_1-1,k_2,k_3}_{\boldsymbol{\alpha_1}-1,\boldsymbol{\alpha_2},\boldsymbol{\alpha_3}}$), the $y$ component of the constitutive equation at the Greville abscissae associated with the discrete $y$ vorticity space ($S^{k_1,k_2-1,k_3}_{\boldsymbol{\alpha_1},\boldsymbol{\alpha_2}-1,\boldsymbol{\alpha_3}}$), and the $z$ component of the constitutive equation at the Greville abscissae associated with the discrete $z$ vorticity space ($S^{k_1,k_2,k_3-1}_{\boldsymbol{\alpha_1},\boldsymbol{\alpha_2},\boldsymbol{\alpha_3}-1}$).

\subsection{Boundary Condition Enforcement}

The no-penetration boundary condition is enforced identically to 2D case: We strongly enforce the normal velocity on face collocation points corresponding to the normal velocity component and remove these points from the set used to collocate the momentum equations. The no-slip boundary condition in the velocity-pressure scheme is also essentially enforced identically to the 2D case and again leads to equations of the form 

\begin{equation}
     -\nu \Delta \mathbf{u}^h + \mathbf{u}^h \cdot \nabla \mathbf{u}^h + \nabla p^h + \frac{C_{pen}^2}{h^2}(\mathbf{u}^h - \mathbf{g}) = \mathbf{f}.
\end{equation}

As the constitutive law relating velocity and vorticity is a vector relation in 3D, the weak enforcement of no-slip boundary conditions is slightly altered in the three field formulation. We again start by considering the weak form shown above, particularly Equation \eqref{eq:3D_weak}. The last term in this equation represents the boundary term which would be used to enforce boundary conditions by replacing terms with their prescribed values. Following the Enhanced Collocation method of \cite{deLorenzis_Neumann_contact}, the equation can be integrated by parts once again, to arrive at a strong form representation given by:

\begin{equation*}
    \int_\Omega \boldsymbol{\psi} \cdot (\boldsymbol{\omega} - \nabla \times \mathbf{u}) d \Omega + \int_{\partial \Omega} (\boldsymbol{\psi} \times \mathbf{u} - \boldsymbol{\psi} \times \mathbf{g}) \cdot \mathbf{n} dA = 0. 
\end{equation*}

\noindent Using the properties of the scalar triple product, we can re-write this as:

\begin{equation*}
    \int_\Omega \boldsymbol{\psi} \cdot (\boldsymbol{\omega} - \nabla \times \mathbf{u}) d \Omega + \int_{\partial \Omega} (\mathbf{u} \times \mathbf{n} - \mathbf{g} \times \mathbf{n}) \cdot \boldsymbol{\psi} dA = 0. 
\end{equation*}

By approximating these integrals as is done in \cite{deLorenzis_Neumann_contact, gomez_cauchygalerkin}, we arrive at a strong form statement including boundary conditions suitable for collocation:

\begin{equation}
    \boldsymbol{\omega}^h - \nabla \times \mathbf{u}^h + \frac{C_{pen}}{h} (\mathbf{u}^h \times \mathbf{n} - \mathbf{g} \times \mathbf{n}) = \mathbf{0}.
\end{equation}

\subsection{Final Collocated Equations}

Once again the entire collocation scheme based on the velocity-pressure formulation is summarized first. Let us again define $\boldsymbol{\tau}^{u_x}_\ell$ for $\ell = 1,...,M^{u_x}$ to be the set of Greville points for the basis of the $x$ velocity component ($S^{k_1,k_2-1,k_3-1}_{\boldsymbol{\alpha_1},\boldsymbol{\alpha_2-1},\boldsymbol{\alpha_3-1}}$) with the points corresponding to no-penetration boundaries removed as discussed previously. Define in a similar manner $\boldsymbol{\tau}^{u_y}_\ell$ for $\ell = 1,...,M^{u_y}$ and $\boldsymbol{\tau}^{u_z}_\ell$ for $\ell = 1,...,M^{u_z}$. The pressure Greville points are defined as $\boldsymbol{\tau}^{p}_\ell$ for $\ell = 1,...,N^{p}$. For this formulation the discrete 3D problem reads:
 
\bigskip

$$
\left\{ \hspace{5pt}
\parbox{6in}{
\noindent Find $\textbf{u}^h \in \boldsymbol{\mathcal{V}}_{h,\mathbf{g}}$ and $P^h \in \mathcal{Q}_{h, 0}$ such that:
\begin{equation}
\begin{split}
    \left(-\nu \frac{\partial^2 u_x^h}{\partial x^2} - \nu \frac{\partial^2 u_x^h}{\partial y^2} - \nu \frac{\partial^2 u_x^h}{\partial z^2} + u^h_x\frac{\partial u^h_x}{\partial x} + u^h_y\frac{\partial u^h_x}{\partial y} + u^h_z\frac{\partial u^h_x}{\partial z} + \frac{\partial p^h}{\partial x} \right) (\boldsymbol{\tau}^{u_x}_\ell) \\ = f_x(\boldsymbol{\tau}^{u_x}_\ell) \quad \forall \boldsymbol{\tau}^{u_x}_\ell \in \Omega
\end{split}
\end{equation}
\begin{equation}
\begin{split}
    \left(-\nu \frac{\partial^2 u_x^h}{\partial x^2} - \nu \frac{\partial^2 u_x^h}{\partial y^2} - \nu \frac{\partial^2 u_x^h}{\partial z^2} + u^h_x\frac{\partial u^h_x}{\partial x} + u^h_y\frac{\partial u^h_x}{\partial y} + u^h_z\frac{\partial u^h_x}{\partial z} + \frac{\partial p^h}{\partial x} + \right. \\ \left. \frac{C_{pen}^2}{h^2}(u^h_x - g_{x})\right)(\boldsymbol{\tau}^{u_x}_\ell) = f_x(\boldsymbol{\tau}^{u_x}_\ell) \quad \forall \boldsymbol{\tau}^{u_x}_\ell \in \partial \Omega
\end{split}
\end{equation}
\begin{equation}
\begin{split}
    \left(-\nu \frac{\partial^2 u_y^h}{\partial x^2} - \nu \frac{\partial^2 u_y^h}{\partial y^2} - \nu \frac{\partial^2 u_y^h}{\partial z^2} + u^h_x\frac{\partial u^h_y}{\partial x} + u^h_y\frac{\partial u^h_y}{\partial y} + u^h_z\frac{\partial u^h_y}{\partial z} + \frac{\partial p^h}{\partial y} \right)(\boldsymbol{\tau}^{u_y}_\ell) \\ = f_y(\boldsymbol{\tau}^{u_y}_\ell) \quad \forall  \boldsymbol{\tau}^{u_y}_\ell \in \Omega
\end{split}
\end{equation}
\begin{equation}
\begin{split}
    \left(-\nu \frac{\partial^2 u_y^h}{\partial x^2} - \nu \frac{\partial^2 u_y^h}{\partial y^2} - \nu \frac{\partial^2 u_y^h}{\partial z^2} + u^h_x\frac{\partial u^h_y}{\partial x} + u^h_y\frac{\partial u^h_y}{\partial y} + u^h_z\frac{\partial u^h_y}{\partial z} + \frac{\partial p^h}{\partial y} + \right. \\ \left. \frac{C_{pen}^2}{h^2}(u^h_y - g_{y})\right)(\boldsymbol{\tau}^{u_y}_\ell) = f_y(\boldsymbol{\tau}^{u_y}_\ell) \quad \forall \boldsymbol{\tau}^{u_y}_\ell \in \partial \Omega
\end{split}
\end{equation}
\begin{equation}
\begin{split}
    \left(-\nu \frac{\partial^2 u_z^h}{\partial x^2} - \nu \frac{\partial^2 u_z^h}{\partial y^2} - \nu \frac{\partial^2 u_z^h}{\partial z^2} + u^h_x\frac{\partial u^h_z}{\partial x} + u^h_y\frac{\partial u^h_z}{\partial y} + u^h_z\frac{\partial u^h_z}{\partial z} + \frac{\partial p^h}{\partial z} \right)(\boldsymbol{\tau}^{u_z}_\ell)  \\= f_z(\boldsymbol{\tau}^{u_z}_\ell) \quad \forall \boldsymbol{\tau}^{u_z}_\ell \in \Omega
\end{split}
\end{equation}
\begin{equation}
\begin{split}
    \left(-\nu \frac{\partial^2 u_z^h}{\partial x^2} - \nu \frac{\partial^2 u_z^h}{\partial y^2} - \nu \frac{\partial^2 u_z^h}{\partial z^2} + u^h_x\frac{\partial u^h_z}{\partial x} + u^h_y\frac{\partial u^h_z}{\partial y} + u^h_z\frac{\partial u^h_z}{\partial z} + \frac{\partial p^h}{\partial z} + \right. \\ \left. \frac{C_{pen}^2}{h^2}(u^h_z - g_{z})\right)(\boldsymbol{\tau}^{u_z}_\ell) = f_z(\boldsymbol{\tau}^{u_z}_\ell) \quad \forall \boldsymbol{\tau}^{u_z}_\ell \in \partial \Omega
\end{split}
\end{equation}
\begin{equation}
    \left(\frac{\partial u^h_x}{\partial x} + \frac{\partial u^h_y}{\partial y} + \frac{\partial u^h_z}{\partial z}\right)(\boldsymbol{\tau}^{p}_\ell)  = 0 \quad \forall \boldsymbol{\tau}^{p}_\ell \in \Omega \cup \partial \Omega.
\end{equation}
}
\right.
$$

\bigskip

Similarly to the velocity, in the three field formulation we also define collocation points for the vorticity component-wise. In particular, let $\boldsymbol{\tau}^{\omega_x}_\ell$ for $\ell = 1,...,N^{\omega_x}$ be the Greville points for the x component of the vorticity, and define $\boldsymbol{\tau}^{\omega_y}_\ell$ for $\ell = 1,...,N^{\omega_y}$ and $\boldsymbol{\tau}^{\omega_z}_\ell$ for $\ell = 1,...,N^{\omega_z}$ similarly. The final, discrete 3D problem for the vorticity-velocity-pressure collocation scheme reads as: 

\bigskip

$$
\left\{ \hspace{5pt}
\parbox{6in}{
\noindent Find $\textbf{u}^h \in \boldsymbol{\mathcal{V}}_{h,\mathbf{g}}$, $P^h \in \mathcal{Q}_{h, 0}$, and $\boldsymbol{\omega}^h \in \boldsymbol{\Psi}_h$ such that:
\begin{equation}
    \left(\nu (\frac{\partial \omega^h_z}{\partial y} - \frac{\partial \omega^h_y}{\partial z} ) + \omega^h_y u^h_z - \omega^h_z u^h_y + \frac{\partial P^h}{\partial x}\right)(\boldsymbol{\tau}_\ell^{u_x})  = f_x(\boldsymbol{\tau}_\ell^{u_x}) \quad \forall \boldsymbol{\tau}^{u_x}_\ell \in \Omega \cup \partial \Omega
\end{equation}
\begin{equation}
   \left(\nu (\frac{\partial \omega^h_x}{\partial z}  - \frac{\partial \omega^h_z}{\partial x} ) + \omega^h_z u^h_x - \omega^h_x u^h_z + \frac{\partial P^h}{\partial y} \right)(\boldsymbol{\tau}_\ell^{u_y})  = f_y(\boldsymbol{\tau}_\ell^{u_y}) \quad \forall \boldsymbol{\tau}^{u_y}_\ell \in \Omega \cup \partial \Omega
\end{equation}
\begin{equation}
   \left( \nu (\frac{\partial \omega^h_y}{\partial x} - \frac{\partial \omega^h_x}{\partial y} + \omega^h_x u^h_y - \omega^h_y u^h_x + \frac{\partial P^h}{\partial z} \right)(\boldsymbol{\tau}_\ell^{u_z})  = f_z(\boldsymbol{\tau}_\ell^{u_z}) \quad \forall \boldsymbol{\tau}^{u_z}_\ell \in \Omega \cup \partial \Omega
\end{equation}
\begin{equation}
    \left(\frac{\partial u^h_x}{\partial x} + \frac{\partial u^h_y}{\partial y} + \frac{\partial u^h_z}{\partial z}\right)(\boldsymbol{\tau}^{p}_\ell)  = 0 \quad \forall \boldsymbol{\tau}^{p}_\ell \in \Omega \cup \partial \Omega
\end{equation}
\begin{equation}
\begin{split}
    &\left(\omega^h_x - (\frac{\partial u^h_z}{\partial y} - \frac{\partial u^h_y}{\partial z} ) \right)(\boldsymbol{\tau}_\ell^{\omega_x}) = 0 \quad \forall \boldsymbol{\tau}^{\omega_x}_\ell \in \Omega
\end{split}
\end{equation}
\begin{equation}
\begin{split}
    &\left(\omega^h_x - (\frac{\partial u^h_z}{\partial y} - \frac{\partial u^h_y}{\partial z} ) + \right. \\ &\left. \frac{C_{pen}}{h}( (u^h_y - g_{y})n_z - (u^h_z - g_{z}) n_y )\right)(\boldsymbol{\tau}_\ell^{\omega_x}) = 0 \quad \forall \boldsymbol{\tau}^{\omega_x}_\ell \in \partial \Omega
\end{split}
\end{equation}
\begin{equation}
\begin{split}
   &\left(\omega^h_y  - (\frac{\partial u^h_x}{\partial z}  - \frac{\partial u^h_z}{\partial x} ) \right) (\boldsymbol{\tau}_\ell^{\omega_y}) = 0 \quad \forall \boldsymbol{\tau}^{\omega_y}_\ell \in \Omega
\end{split}
\end{equation}
\begin{equation}
\begin{split}
   &\left(\omega^h_y  - (\frac{\partial u^h_x}{\partial z}  - \frac{\partial u^h_z}{\partial x} ) + \right. \\ & \left. \frac{C_{pen}}{h}( (u^h_z - g_{z})n_x - (u^h_x - g_{x}) n_z ) \right) (\boldsymbol{\tau}_\ell^{\omega_y}) = 0 \quad \forall \boldsymbol{\tau}^{\omega_y}_\ell \in \partial \Omega
\end{split}
\end{equation}
\begin{equation}
\begin{split}
    &\left( \omega^h_z  - (\frac{\partial u^h_y}{\partial x} - \frac{\partial u^h_x}{\partial y} ) \right)(\boldsymbol{\tau}_\ell^{\omega_z}) = 0 \quad \forall \boldsymbol{\tau}^{\omega_z}_\ell \in \Omega
\end{split}
\end{equation}
\begin{equation}
\begin{split}
    &\left( \omega^h_z  - (\frac{\partial u^h_y}{\partial x} - \frac{\partial u^h_x}{\partial y} ) + \right. \\ & \left. \frac{C_{pen}}{h}( (u^h_x - g_{x})n_y - (u^h_y - g_{y}) n_x) \right)(\boldsymbol{\tau}_\ell^{\omega_z}) = 0 \quad \forall \boldsymbol{\tau}^{\omega_z}_\ell \in \partial \Omega.
\end{split}
\end{equation}
}
\right.
$$

\bigskip

\section{Numerical Results on Cubic Domains}

To verify that the schemes properly extend into 3D, two sample problems are considered. First, a manufactured solution gives even more insight into the convergence properties of the methods. Then the three-dimensional lid-driven cavity problem is considered and the results are compared with established literature. 

\subsection{Three-Dimensional Manufactured Solution}

In 3D, we also start our numerical studies by considering a manufactured solution. In this case, the exact solution represents the flow around a single vortex filament within the unit cube. We define a potential function as

\begin{equation}
    \tilde{\boldsymbol{\phi}} = \left[
    \begin{array}{c}
         x(x-1)y^2(y-1)^2z^2(z-1)^2  \\
         0 \\
         x^2(x-1)^2y^2(y-1)^2z(z-1)
    \end{array}
    \right],
\end{equation}

\noindent through which we can define the velocity field as 

\begin{equation}
    \tilde{\textbf{u}} = \nabla \times \tilde{\boldsymbol{\phi}},
\end{equation}

\noindent and the vorticity as 

\begin{equation}
    \tilde{\boldsymbol{\omega}} = \nabla \times \tilde{\textbf{u}}.
\end{equation}

\noindent Finally, we specify the pressure field as 

\begin{equation}
    \tilde{p} = \sin(\pi x)\sin(\pi y) - \frac{4}{\pi^2}.
\end{equation}

For the velocity-pressure scheme we define the forcing term on the right hand sign of the momentum equations as

\begin{equation}
    \textbf{f} = -\nu \Delta \tilde{\textbf{u}} + \tilde{\textbf{u}} \cdot \nabla \tilde{\textbf{u}} + \nabla \tilde{p}, 
\end{equation}

\noindent while for the vorticity-velocity-pressure scheme the forcing term is given by

\begin{equation}
    \textbf{f} = -\nu \Delta \tilde{\textbf{u}} + \tilde{\boldsymbol{\omega}} \times \tilde{\textbf{u}} + \nabla \tilde{P}.
\end{equation}

\noindent Once again we enforce homogeneous Dirichlet boundary conditions everywhere and require that the kinematic pressure field has zero average. With these conditions the discrete solution should again converge to the quantities above with mesh refinement.

\begin{figure}
\centering
\subfloat[Velocity $L^2$ error]{\label{sfig:3D_conv_lap_a}\includegraphics[width=.5\textwidth]{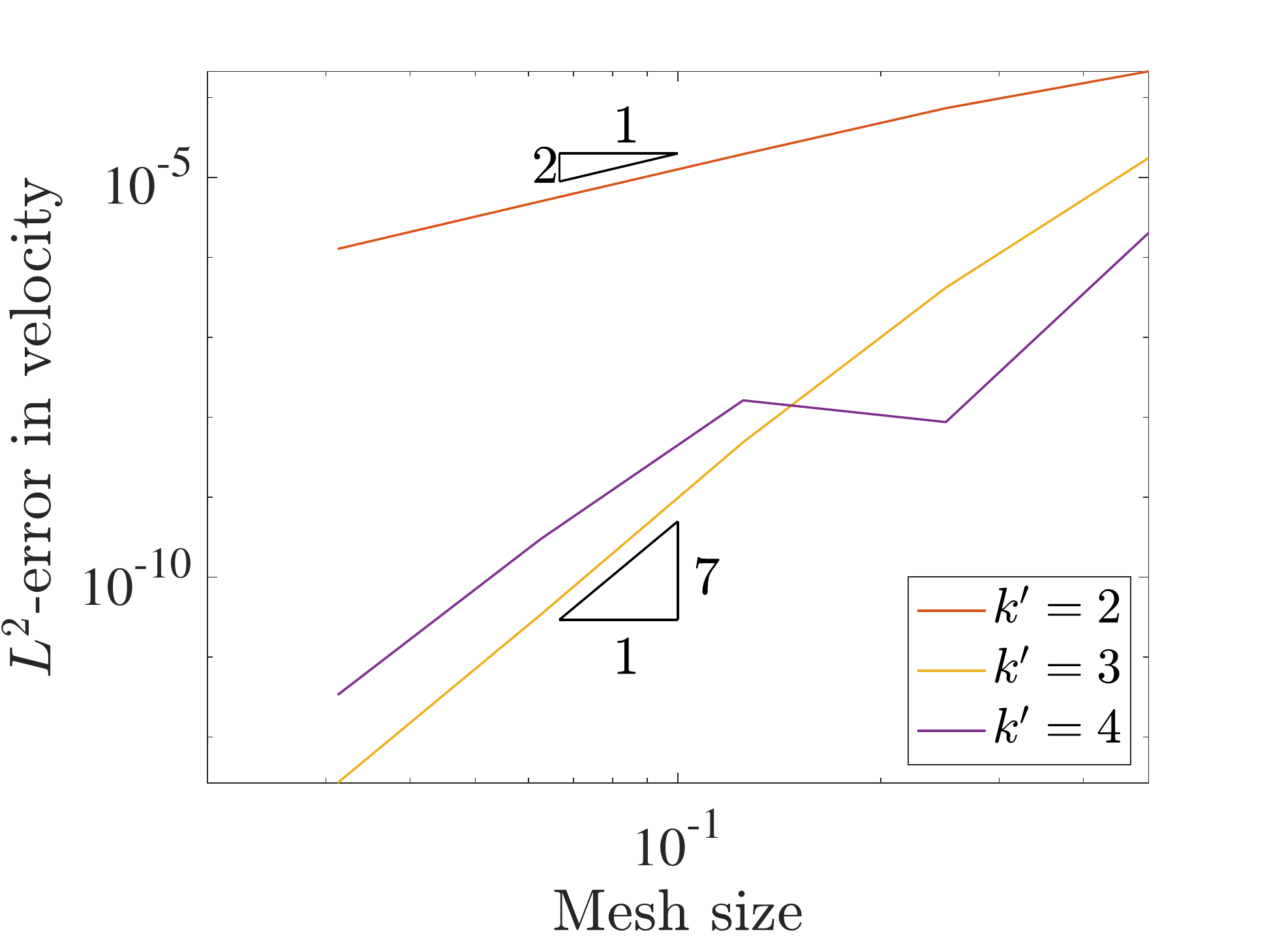}}\hfill
\subfloat[Velocity $H^1$ error]{\label{sfig:3D_conv_lap_b}\includegraphics[width=.5\textwidth]{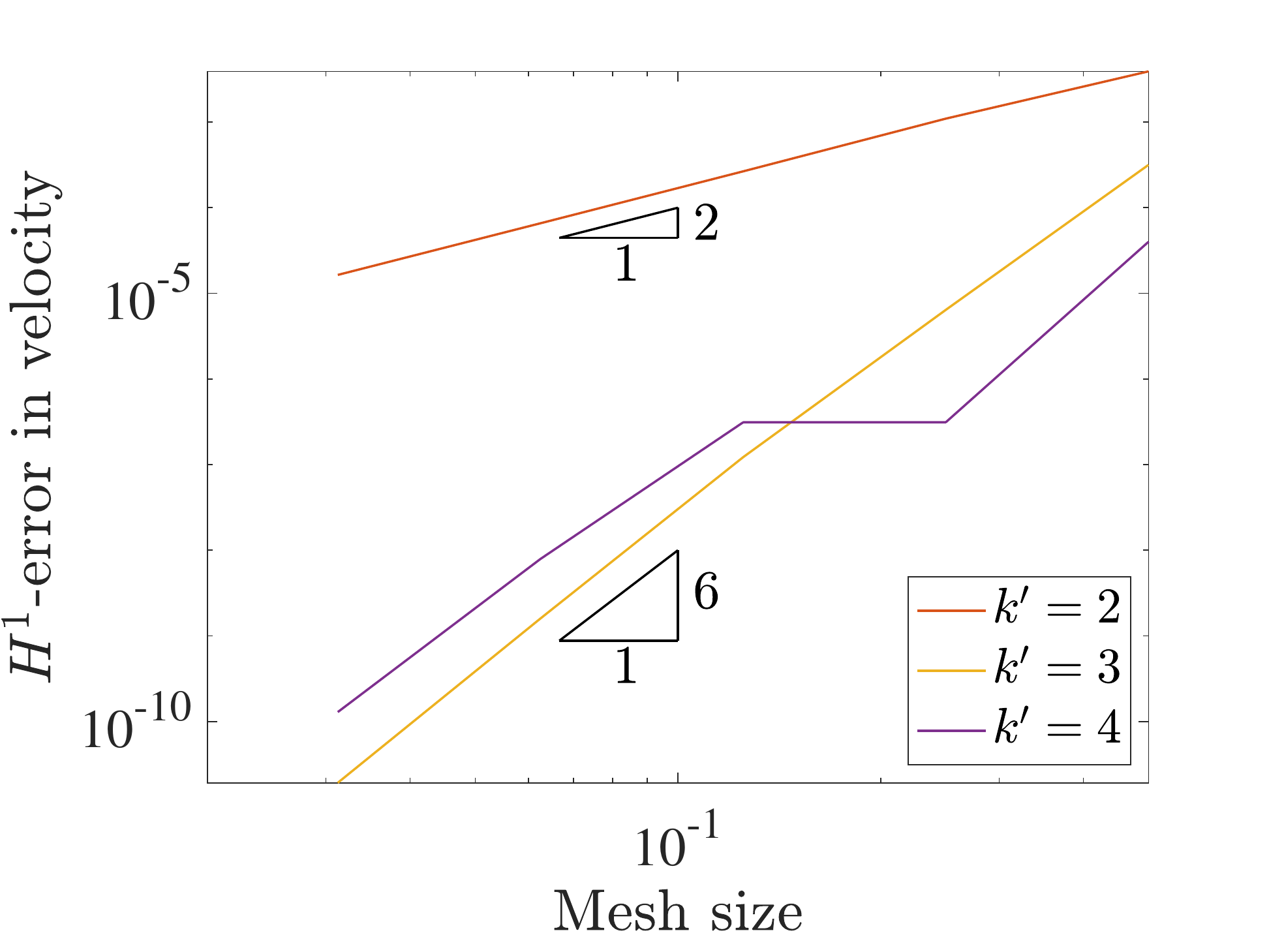}} \\
\subfloat[Pressure $L^2$ error]{\label{sfig:3D_conv_lap_c}\includegraphics[width=.5\textwidth]{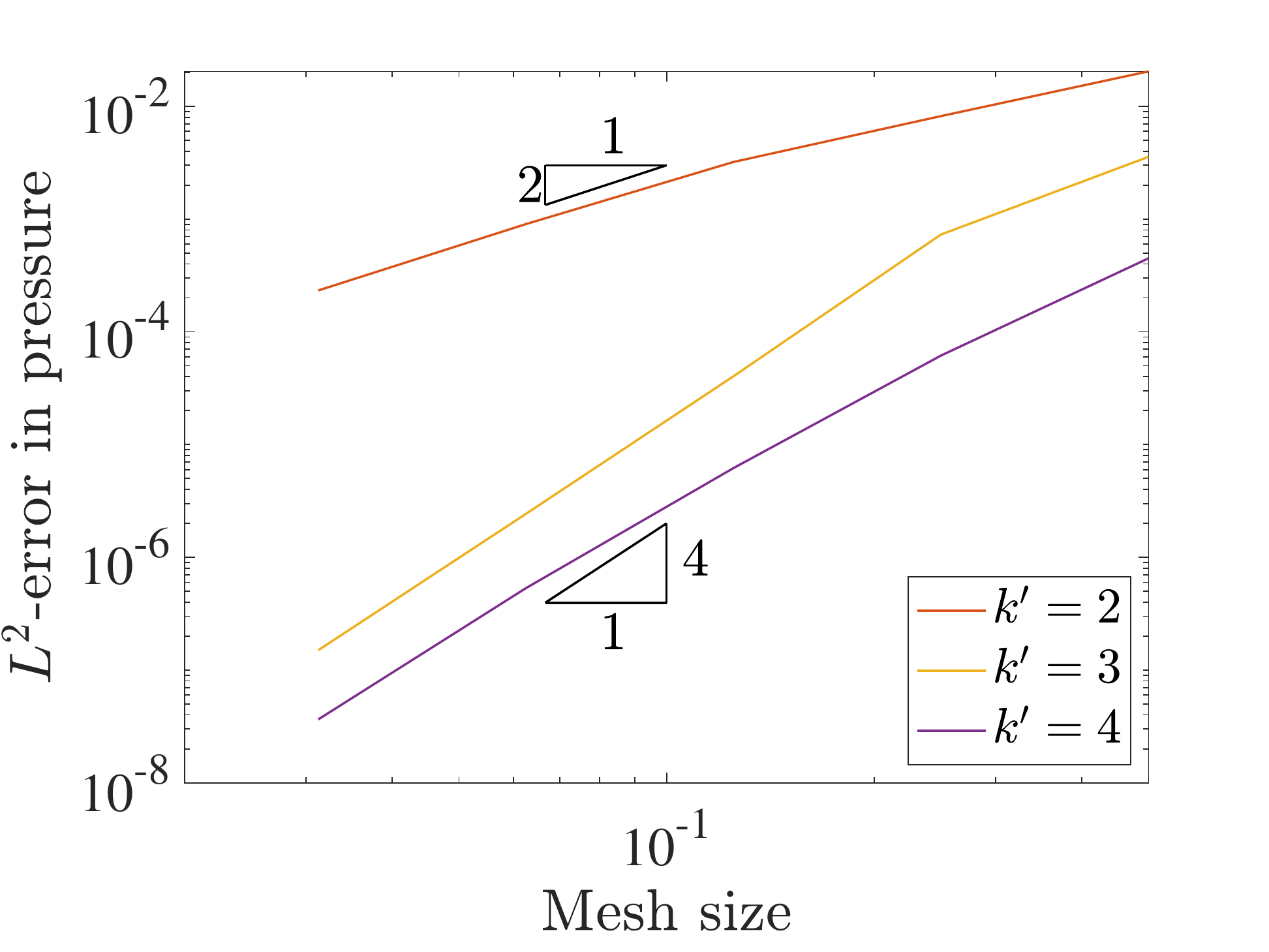}}\hfill
\subfloat[Pressure $H^1$ error]{\label{sfig:3D_conv_lap_d}\includegraphics[width=.5\textwidth]{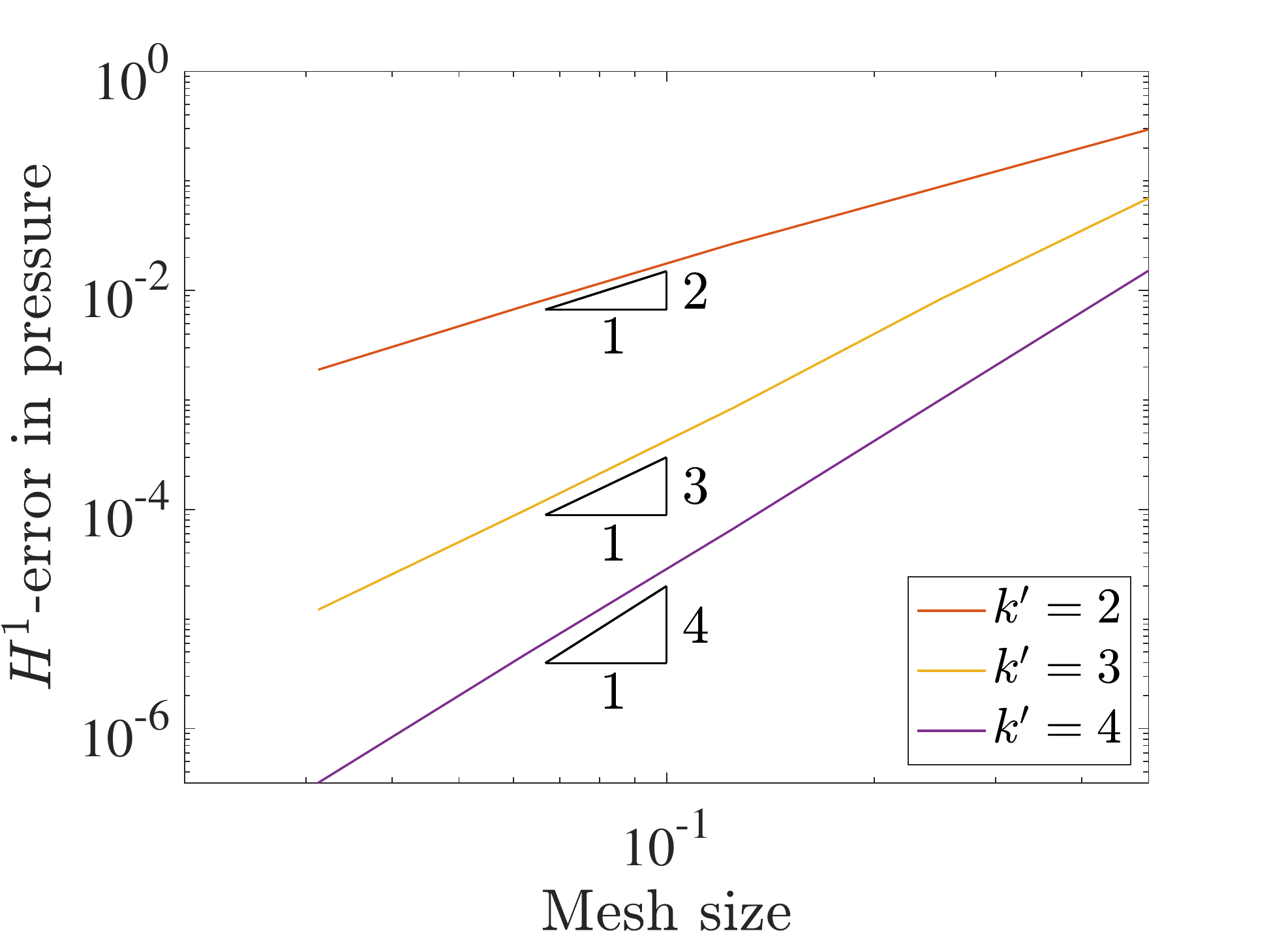}} \\
\caption{Errors in 3D manufactured vortex solution for velocity-pressure formulation}
\label{fig:3D_conv_lap}
\end{figure}
\begin{figure}
\centering
\subfloat[Velocity $L^2$ error]{\label{sfig:3D_conv_a}\includegraphics[width=.5\textwidth]{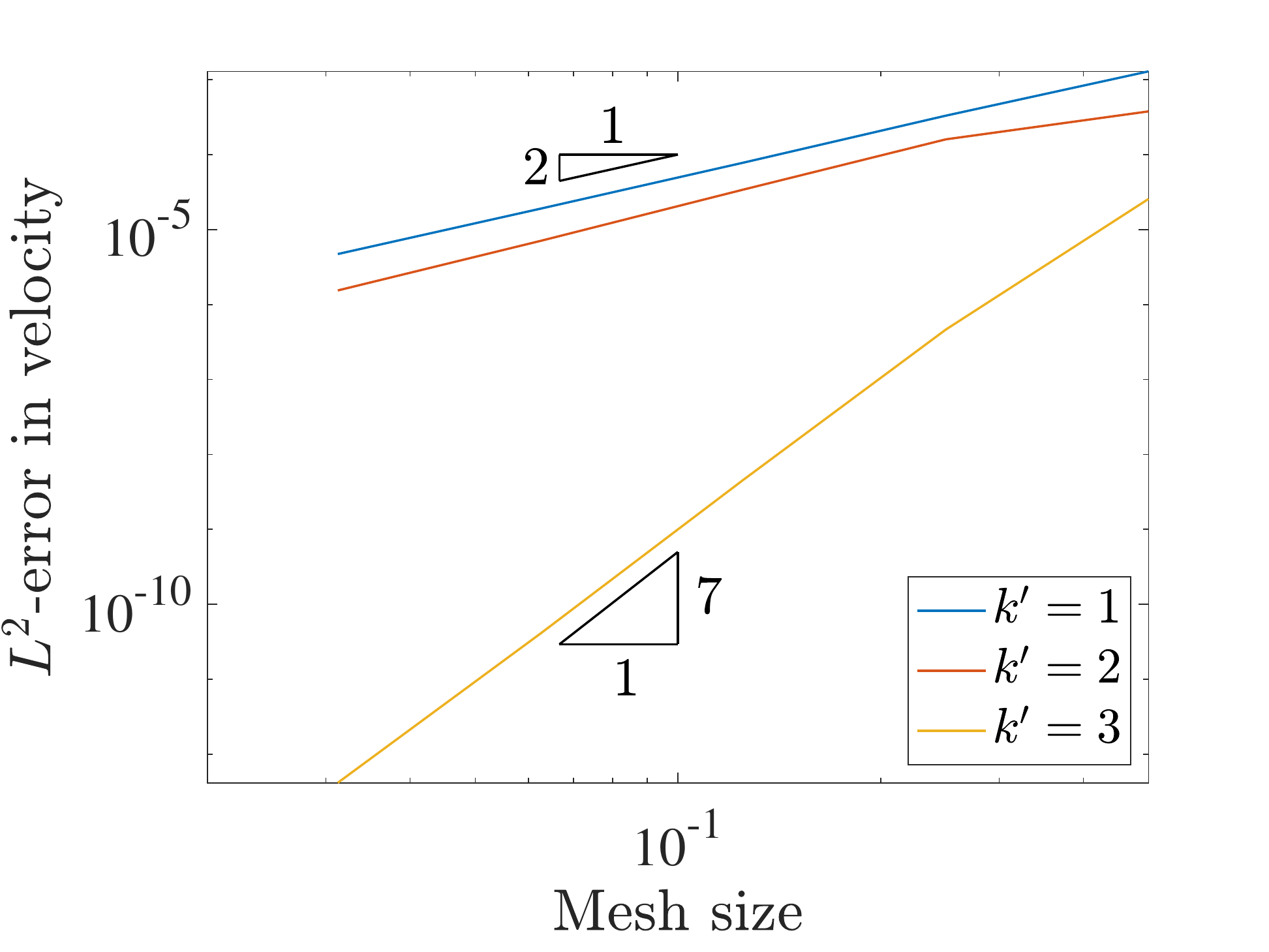}}\hfill
\subfloat[Velocity $H^1$ error]{\label{sfig:3D_conv_b}\includegraphics[width=.5\textwidth]{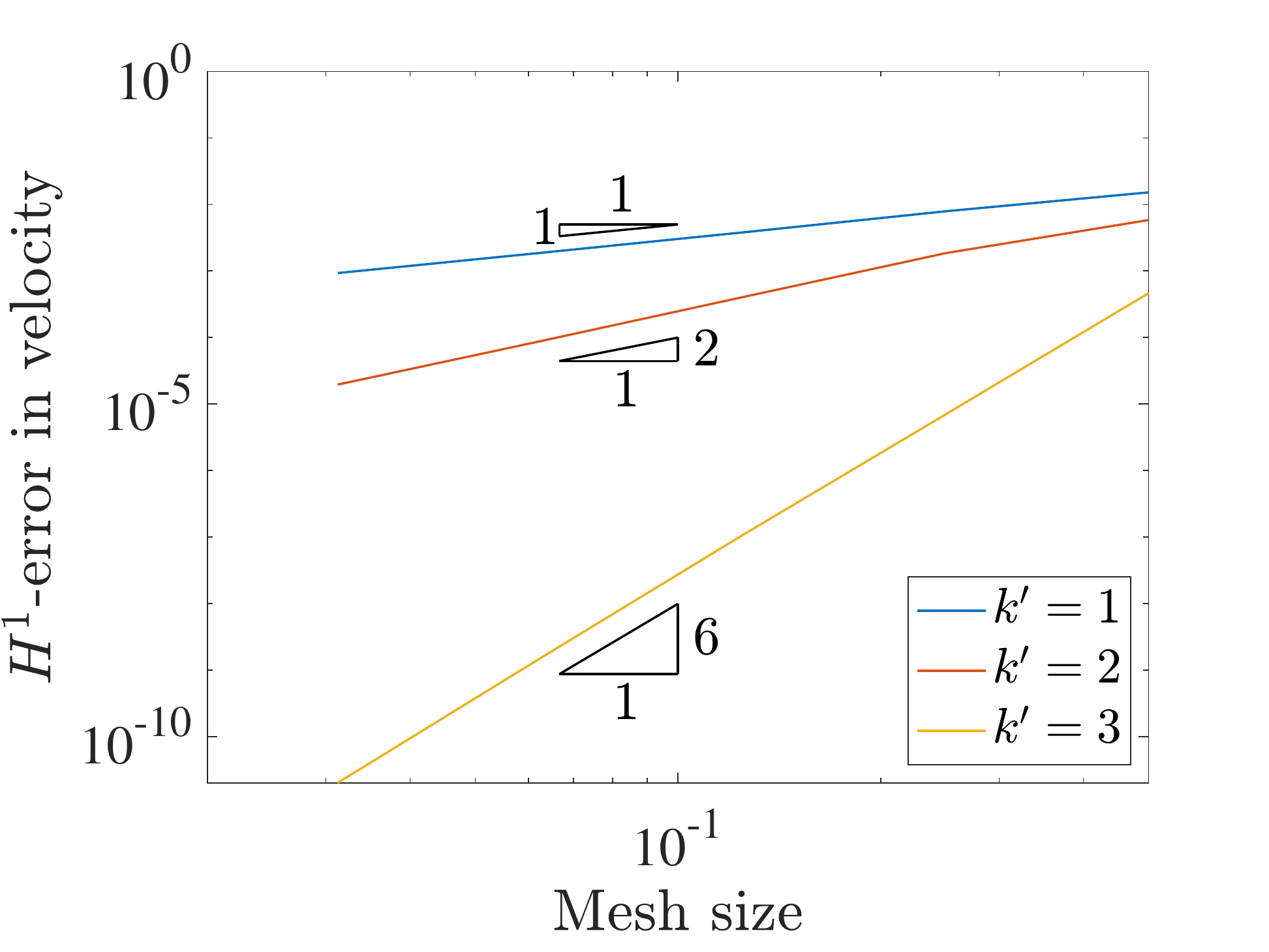}} \\
\subfloat[Pressure $L^2$ error]{\label{sfig:3D_conv_c}\includegraphics[width=.5\textwidth]{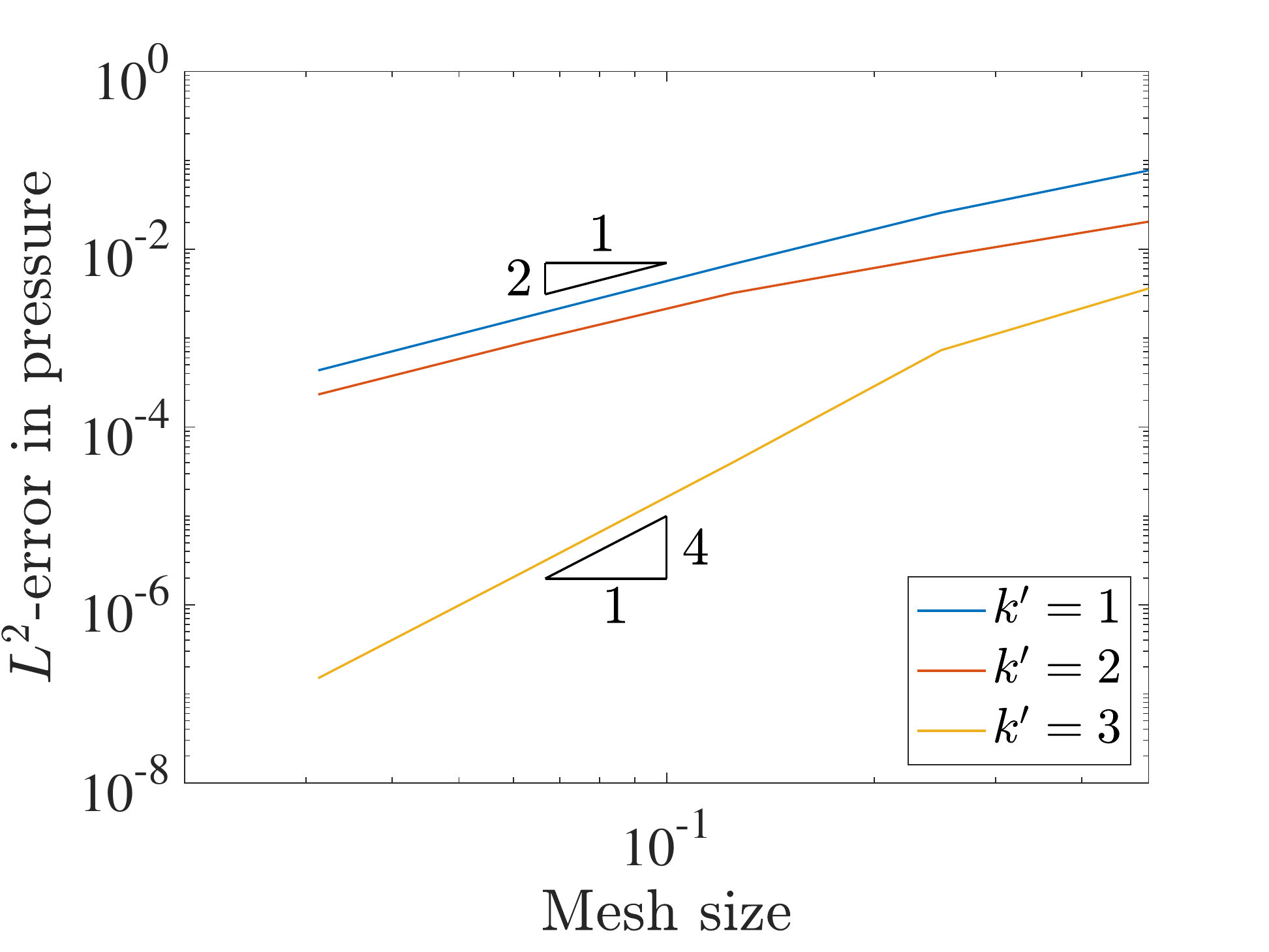}}\hfill
\subfloat[Pressure $H^1$ error]{\label{sfig:3D_conv_d}\includegraphics[width=.5\textwidth]{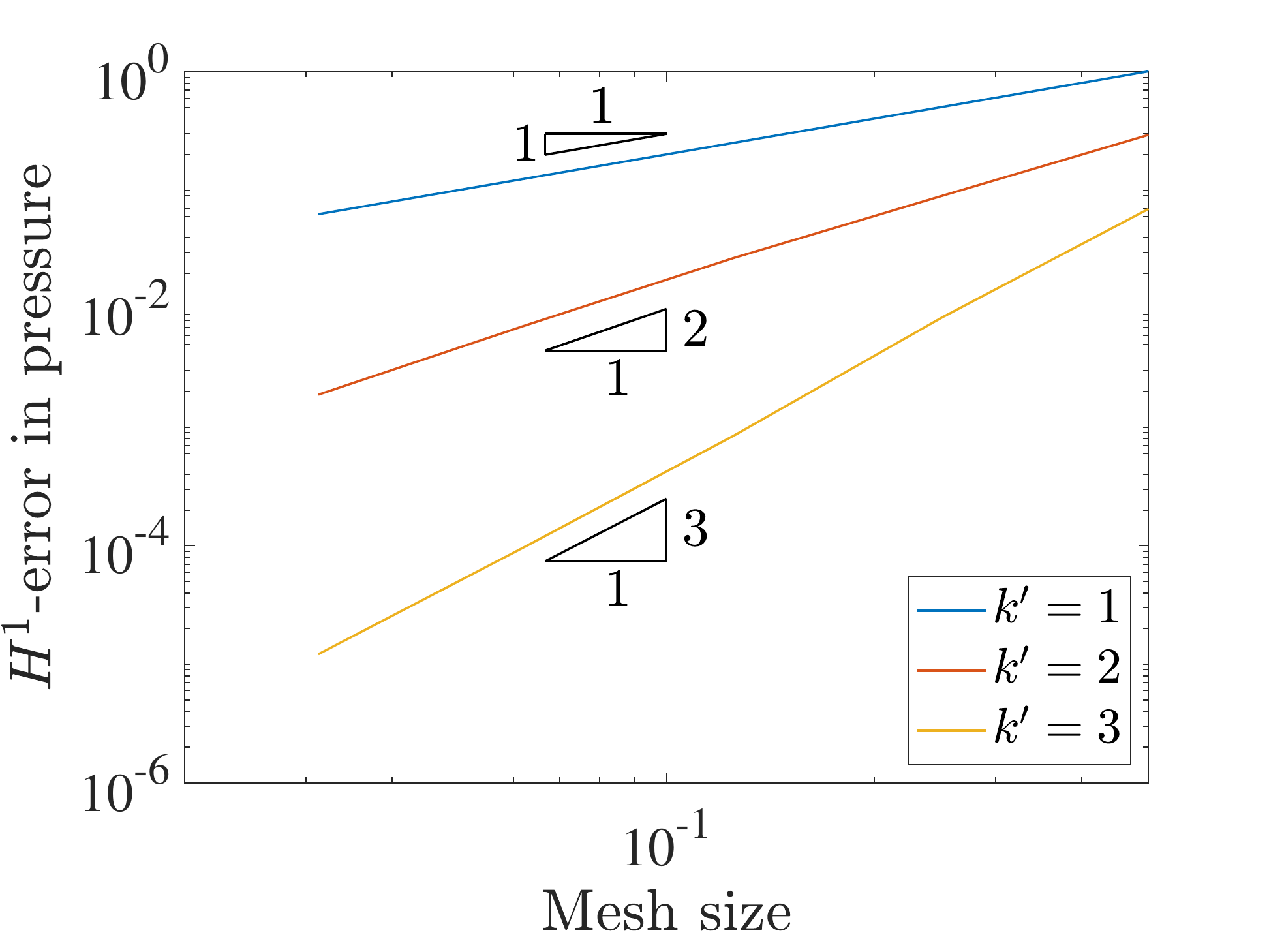}} \\
\subfloat[Vorticity $L^2$ error]{\label{sfig:3D_conv_e}\includegraphics[width=.5\textwidth]{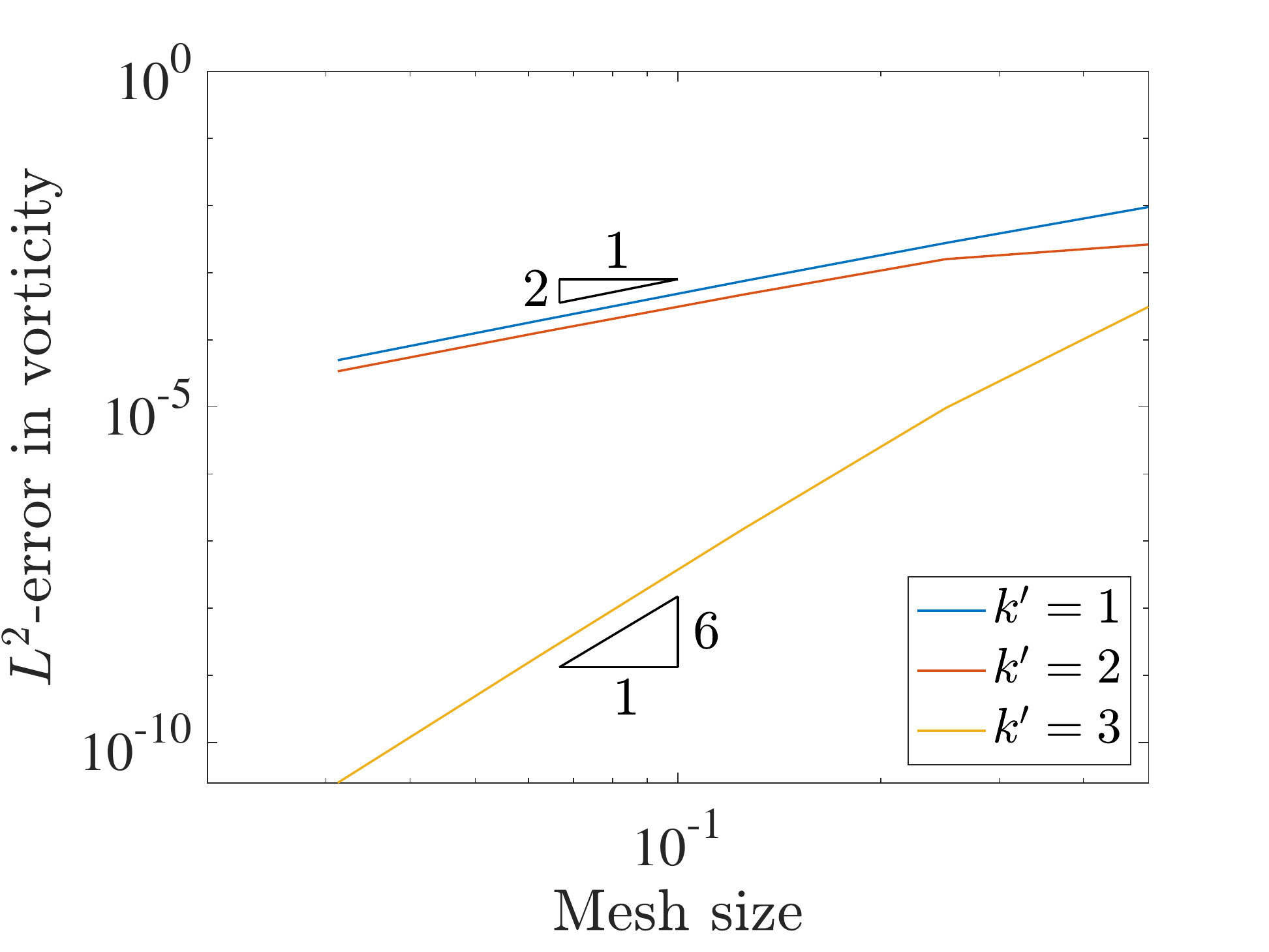}}\hfill
\subfloat[Vorticity $H^1$ error]{\label{sfig:3D_conv_f}\includegraphics[width=.5\textwidth]{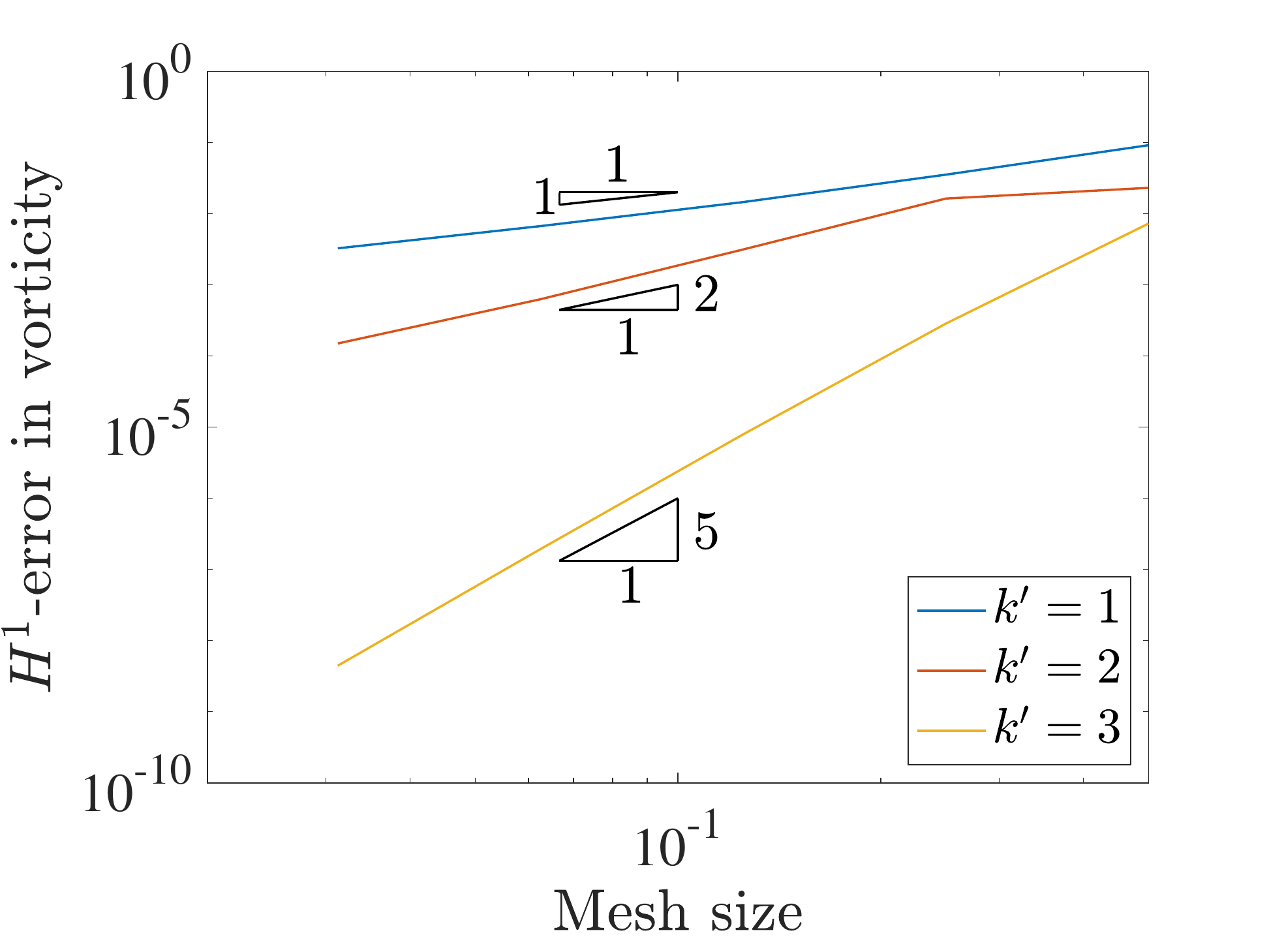}} \\
\caption{Errors in 3D manufactured vortex solution for vorticity-velocity-pressure formulation}
\label{fig:3D_conv}
\end{figure}

Similar to the 2D case, we set $Re = \frac{1}{\nu} = 1$ and measure the errors produced on a variety of grids in the $L^2$ norm and $H^1$ semi-norm. Figure \ref{fig:3D_conv_lap} shows the results for the velocity-pressure scheme while Figure \ref{fig:3D_conv} details the errors for the vorticity-velocity-pressure scheme.

We start by noting that when $k' < 3$ in both cases, everything behaves in the same manner as in the 2D setting. Once $k' \geq 3$ we start to see very fast convergence rates and some pre-asymptotic type behavior in the velocity errors produced by both schemes. This can be explained by talking a closer look at the exact velocity field for this problem. In fact, the exact velocity field is given by a quartic polynomial in each direction and this solution is actually contained within the discrete velocity approximation space for $k' \geq 3$. If we were using a pressure robust Galerkin method, the velocity error would be zero. Since the collocation scheme is not pressure robust, we obtain superconvergence rather than exactly zero error.

The pressure convergence results also show some interesting behavior. While the vorticity-velocity-pressure scheme seems to behave in the same manner as in 2D, the velocity-pressure scheme seems to be recovering the faster rates seen in the three field scheme. We believe that this is a consequence of the superconvergence of velocity.




\subsection{Three-Dimensional Lid-Driven Cavity}

The next numerical study that we perform is on the 3D lid-driven cavity flow. Consider again the cavity setup describing the 2D flow, but now extend the square cavity by unit length in the out-of-page direction, thus making it a cube. The point singularities of the 2D case now extend along the top edges of the cube and we expect to see more influence of 3D boundary effects \cite{ku1987_3Dcavity}. 

\begin{figure}
\centering
\subfloat[Velocity-pressure formulation]{\label{sfig:vvp_grid_mapped_a}\includegraphics[width=.5\textwidth]{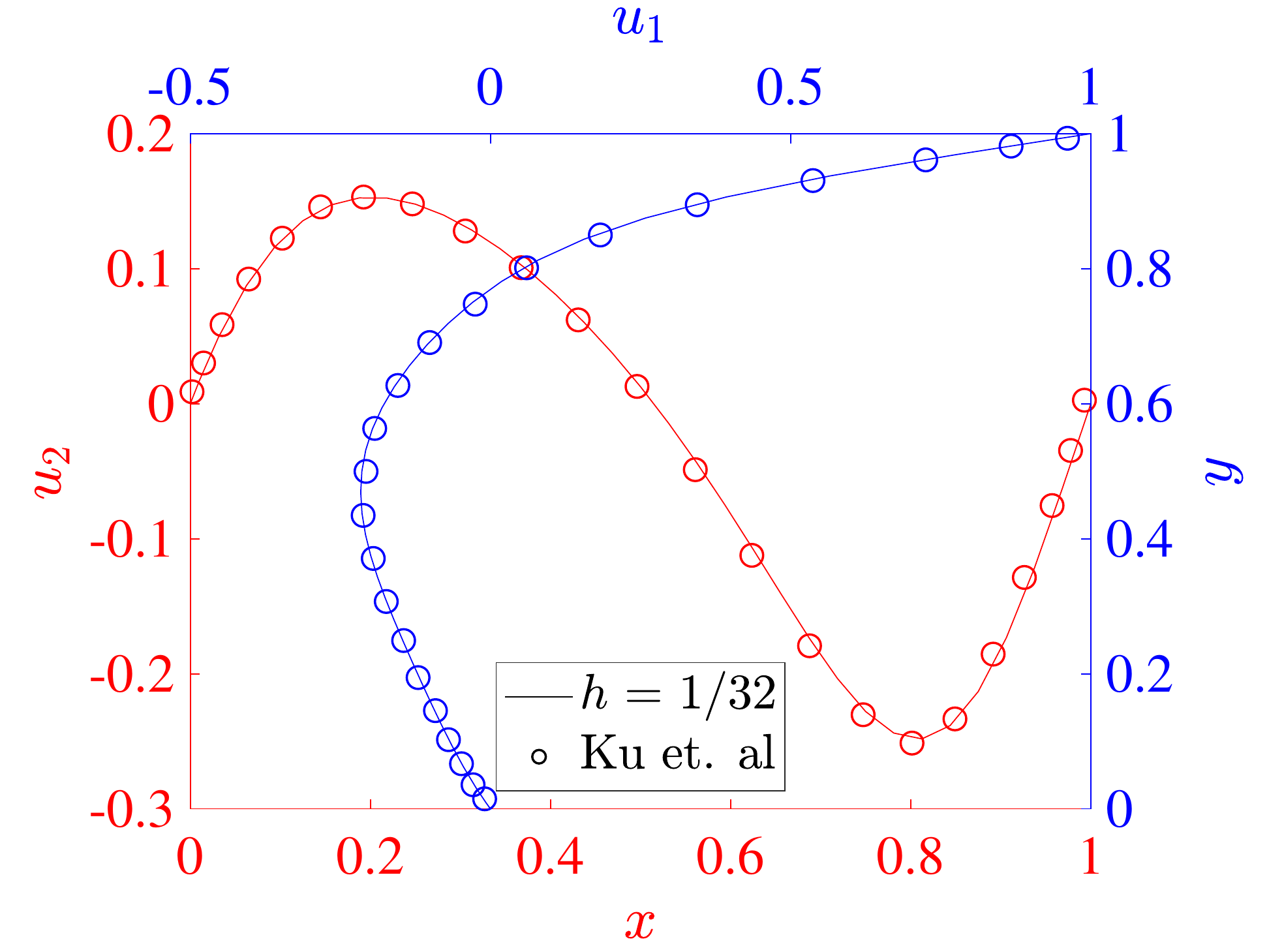}}\hfill
\subfloat[Vorticity-velocity-pressure formulation]{\label{sfig:vvp_grid_mapped_b}\includegraphics[width=.5\textwidth]{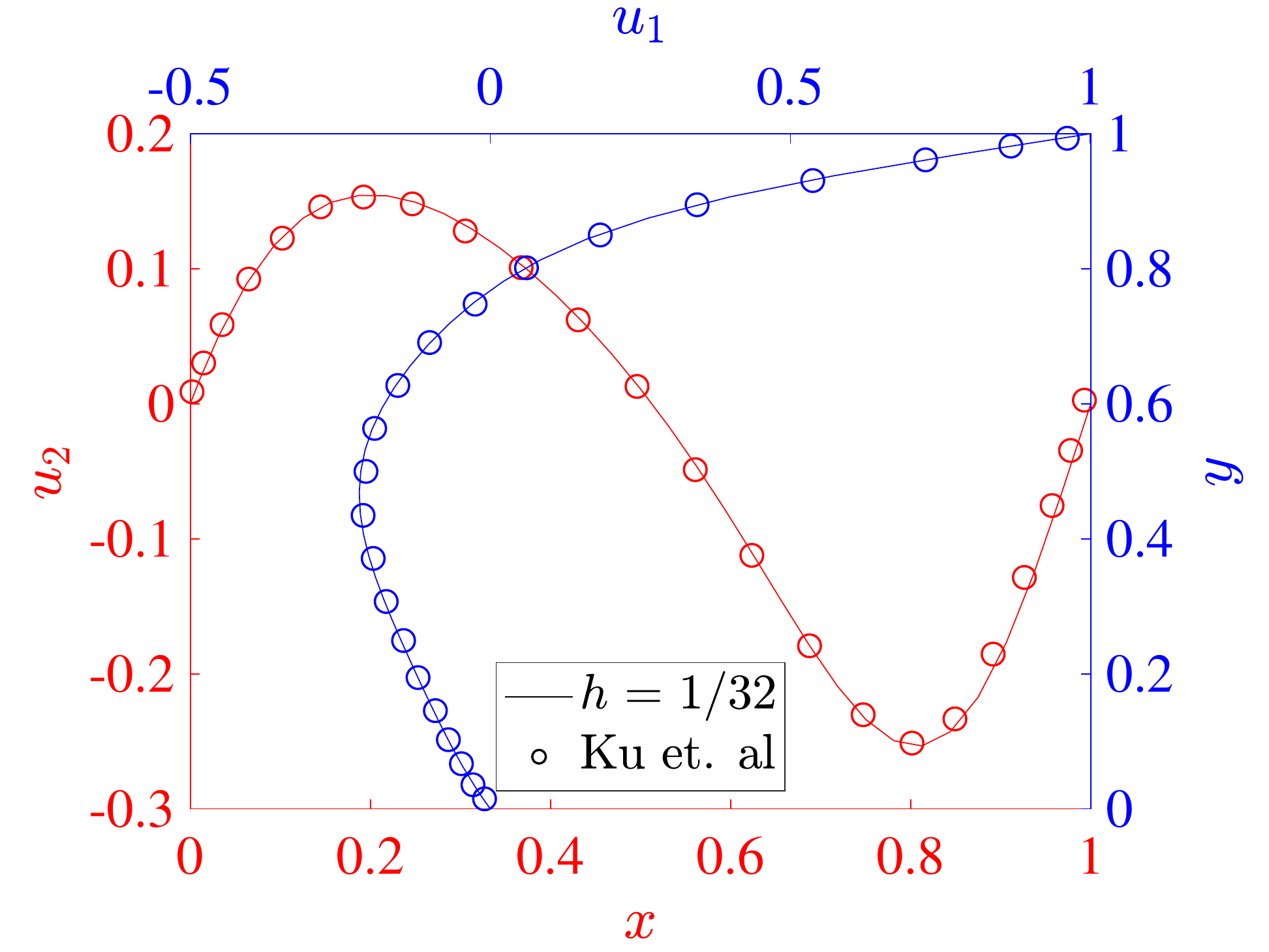}} \\
\caption{Centerline velocity profile for 3D lid-driven cavity using both formulations, $k' = 2$. Red curves and axes represent the vertical velocity along the horizontal centerline, while blue curves and axes represent the horizontal velocity along the vertical centerline.}
\label{fig:3d_cavity}
\end{figure}

In our tests we again set the wall speed $U = 1$, the side length $H = 1$, and consider $Re = \frac{UH}{1} = 100$. We use an unstretched mesh with 32 elements per side and $k' = 2$, and compare the $x$ velocity along the vertical centerline and the $y$ velocity along the horizontal centerline with the pseudospectral results from \cite{ku1987_3Dcavity}. Figure \ref{fig:3d_cavity} shows the results with each formulation. Once again the results match very well with the literature, and it seems as though the results from the three field formulation match with the reference results slightly better that the two field results. 



\section{Collocation Methods on Mapped Domains}

As the last main component of this paper we shift our focus to problems posed on more complicated domains. We will present some theory for both 2D and 3D problems, but for simplicity we will focus the development of numerical schemes for the 2D, linear Stokes equations. However, the results would generalize to the nonlinear, 3D setting as well. We will also focus on the rotational form of the equations, as the first order nature enables easier mappings between domains. 

The main idea of this section is the mapping back to a parametric reference domain, i.e. a square in 2D or a cube in 3D. The previous sections detail how to develop collocation schemes on these simple geometries, thus simply pulling the equations and unknowns back to the reference domain, collocating as before, and pushing the results forward to the physical domain gives our numerical solution. 

Let $\hat{\Omega}$ be the parametric domain (the unit square in 2D or the unit cube in 3D), and let $\Omega$ be the physical domain. We define the function $\mathbf{F}$ as mapping from $\hat{\Omega}$ to $\Omega$. Let $\mathbf{DF}$ be the Jacobian of the parametric mapping, and define

\begin{equation}
    J = \textup{Det}(\mathbf{DF}),
\end{equation}
\begin{equation}
    \mathbf{C} = (\mathbf{DF})^T(\mathbf{DF}).
\end{equation}

Next we can define the pull-back operators in 3D as 

\begin{align}
   \iota_{\Phi}(\phi) &= ( \phi \circ \mathbf{F}), \\
   \iota_{\boldsymbol{\omega}}(\boldsymbol{\psi}) &=(\mathbf{DF})^{T} ( \boldsymbol{\psi} \circ \mathbf{F}), \\
    \iota_{\mathbf{u}}(\mathbf{v}) &= J (\mathbf{DF})^{-1} ({\mathbf{v} \circ \mathbf{F}}), \\
    \iota_p(q) &= J (q \circ \mathbf{F}).
\end{align}

\noindent We define the pulled-back unknowns on the reference domain via the $\iota$ maps, specifically $\hat{\mathbf{u}} = \iota_{\mathbf{u}}(\mathbf{u}) $, $\hat{p} = \iota_p(p) $, and $\hat{\boldsymbol{\omega}} = \iota_{\boldsymbol{\omega}}(\boldsymbol{\omega})$. These are the unknowns for which we solve using collocation, and the physical domain solution is then obtained via the corresponding push-forward. Importantly, the push-forward of velocity as defined above maps divergences to divergences and preserves nullity of normal components. Similarly the push-forward of pressure preserves the nullity of the integral operator. These facts imply that the following commuting diagram exists:

\begin{equation}
    \begin{CD}
        \mathbb{R} @>>> \Phi @>\nabla>> \boldsymbol{\Psi} @>\nabla\times>> \boldsymbol{\mathcal{V}} @>\nabla \cdot>> \mathcal{Q} @>>> 0\\
        @. @VV\iota_{\Phi}V @VV\iota_{\boldsymbol{\omega}}V @VV\iota_{\mathbf{u}}V @VV\iota_{p}V @.\\
        \mathbb{R} @>>> \Phi @>\nabla>> \hat{\boldsymbol{\Psi}} @>\nabla\times>> \hat{\boldsymbol{\mathcal{V}}} @>\nabla \cdot>> \hat{\mathcal{Q}} @>>> 0,
    \end{CD}
    \label{eq:deRham3D_mapped}
\end{equation}

\noindent where now the hat spaces correspond to the ones defined over the parametric domain, and are identical to the ones used in the previous sections of this paper. Moreover, by composing the $\iota$ maps with the projectors from the de Rham complex in the square domain setting, we arrive at a new commuting diagram between the physical domain continuous spaces and the discrete spaces in the physical domain defined by the push-forward of the discrete spaces chosen for the unit square.  

For completeness we also define the 2D pull-back operators

\begin{align}
   \iota_{\omega}(\psi) &= \psi \circ \mathbf{F}, \\
    \iota_{\mathbf{u}}(\mathbf{v}) &= J (\mathbf{DF})^{-1} ({\mathbf{v} \circ \mathbf{F}}), \\
    \iota_p(q) &= J (q \circ \mathbf{F}).
\end{align}

\noindent In 2D a commuting diagram also exists:

\begin{equation}
    \begin{CD}
        \mathbb{R} @>>> \Psi @>\nabla^\perp>> \boldsymbol{\mathcal{V}} @>\nabla \cdot>> \mathcal{Q} @>>> 0\\
        @. @VV\iota_{\omega}V @VV\iota_{\mathbf{u}}V @VV\iota_{p}V @.\\
        \mathbb{R} @>>> \hat{\Psi} @>\nabla^\perp>> \hat{\boldsymbol{\mathcal{V}}} @>\nabla \cdot>> \hat{\mathcal{Q}} @>>> 0.
    \end{CD}
    \label{eq:deRham2D_mapped}
\end{equation}

Next we begin the process of mapping the governing equations back to the reference domain. We start with Equations \eqref{eq:vvp_mom_3D} - \eqref{eq:vvp_const_3D} for the rotational form of the 3D Navier-Stokes equations. The Stokes equations are recovered by simply removing the nonlinear term in the momentum equation, Equation \eqref{eq:vvp_mom_3D}, and noting now that the pressure becomes the standard kinematic pressure $p$. In the momentum equation, the viscous term is mapped back to the reference domain via

\begin{equation}
    (\nabla \times \boldsymbol{\omega}) \circ \mathbf{F} = J^{-1} (\mathbf{DF}) (\hat{\nabla} \times \iota_{\boldsymbol{\omega}} (\boldsymbol{\omega})) = J^{-1} (\mathbf{DF}) (\hat{\nabla} \times \hat{\boldsymbol{\omega}}),
\end{equation}

\noindent and the pressure term is mapped to 

\begin{equation}
    (\nabla p) \circ \mathbf{F} = (\mathbf{DF})^{-T} \hat{\nabla} (\iota_{\Phi}(p)) = (\mathbf{DF})^{-T} \hat{\nabla} (J^{-1}\iota_{p}(p)) = (\mathbf{DF})^{-T} \hat{\nabla} (J^{-1}\hat{p}).
\end{equation}

Within the continuity equation, Equation \eqref{eq:vvp_cont_3D}, the divergence is mapped via 

\begin{equation}
   ( \nabla \cdot \mathbf{u} ) \circ \mathbf{F} = J^{-1} \hat{\nabla} \cdot (\iota_{\mathbf{u}}(\mathbf{u})) = J^{-1} \hat{\nabla} \cdot \hat{\mathbf{u}}.
\end{equation}

Finally, in the constitutive law, Equation \eqref{eq:vvp_const_3D}, the curl term is mapped similarly to the viscous momentum term

\begin{equation}
\begin{split}
    (\nabla \times \mathbf{u} ) \circ \mathbf{F} =  J^{-1} (\mathbf{DF}) (\hat{\nabla} \times \iota_{\boldsymbol{\omega}} (\mathbf{u})) = J^{-1} (\mathbf{DF}) (\hat{\nabla} \times ((\mathbf{DF})^{T} \mathbf{u})) \\ = J^{-1} (\mathbf{DF}) (\hat{\nabla} \times ((\mathbf{DF})^{T} (J^{-1}(\mathbf{DF}) \hat{\mathbf{u}}))) \\ = J^{-1} (\mathbf{DF}) (\hat{\nabla} \times (J^{-1} \mathbf{C} \hat{\mathbf{u}})).
    \end{split}
\end{equation}

Now we pull each equation back to the reference domain via the corresponding $\iota$ map, so the momentum equations are pulled back via $\iota_{\mathbf{u}}$, the continuity equation is pulled back with $\iota_p$ and the constitutive law is pulled back with $\iota_{\boldsymbol{\omega}}$. For brevity, we will not state the full form of the mapped equations in 3D, but instead state just the 2D form. This arises in a similar way as the 2D rotational form of the Navier-Stokes equations was generated from the 3D equations. In particular we can simply write the equations out component-wise and note that z velocities as well as derivatives in the z direction are zero. This yields:

\bigskip

$$
\left\{ \hspace{5pt}
\parbox{6in}{
\noindent Given $\nu \in \mathbb{R}^+$, $\hat{\textbf{f}} : \hat{\Omega} \rightarrow \mathbb{R}^2$, and $\hat{\textbf{g}} : \partial \hat{\Omega} \rightarrow \mathbb{R}^2$, find $\hat{\textbf{u}} : \hat{\Omega} \rightarrow \mathbb{R}^2$, $\hat{p} : \hat{\Omega} \rightarrow \mathbb{R}$, and  $\hat{\omega} : \hat{\Omega} \rightarrow \mathbb{R}$ such that:
\begin{equation}
    \nu \frac{\partial \hat{\omega}}{\partial \hat{y}} + J C^{-1}_{11}\frac{\partial (J^{-1}\hat{p})}{\partial \hat{x}} + J C^{-1}_{12}\frac{\partial (J^{-1}\hat{p})}{\partial \hat{y}} = \hat{f}_1 \quad \textup{in} \quad \hat{\Omega} \label{eq:vvp_momx_2D_mapped}
\end{equation}
\begin{equation}
    - \nu \frac{\partial \hat{\omega}}{\partial \hat{x}} + J C^{-1}_{21}\frac{\partial (J^{-1}\hat{p})}{\partial \hat{x}} + J C^{-1}_{22}\frac{\partial (J^{-1}\hat{p})}{\partial \hat{y}} = \hat{f}_2 \quad \textup{in} \quad \hat{\Omega} \label{eq:vvp_momy_2D_mapped}
\end{equation}
\begin{equation}
    \hat{\nabla} \cdot \hat{\textbf{u}} = 0 \quad \textup{in} \quad \hat{\Omega} \label{eq:vvp_cont_2D_mapped}
\end{equation}
\begin{equation}
  \hat{\omega} - J^{-1}(\frac{\partial}{\partial \hat{x}}(J^{-1} (C_{21} \hat{u}_x + C_{22} \hat{u}_y)) - \frac{\partial}{\partial \hat{y}}(J^{-1} (C_{11} \hat{u}_x + C_{12} \hat{u}_y))) = 0 \quad \textup{in} \quad \hat{\Omega} \label{eq:vvp_const_2D_mapped}
\end{equation}
\begin{equation}
     \hat{\textbf{u}} =  \hat{\textbf{g}} \quad \textup{on} \quad \partial \hat{\Omega}, 
\end{equation}

\noindent  where $\hat{\mathbf{f}} = \iota_{\mathbf{u}}(\mathbf{f})$ and $\hat{\mathbf{g}} = \iota_{\mathbf{u}}(\mathbf{g})$.
}
\right.
$$

\bigskip

\begin{figure}
\centering
\subfloat[Before strong enforcement of no penetration conditions]{\label{sfig:vvp_grid_mapped_nobc}\includegraphics[width=.4\textwidth]{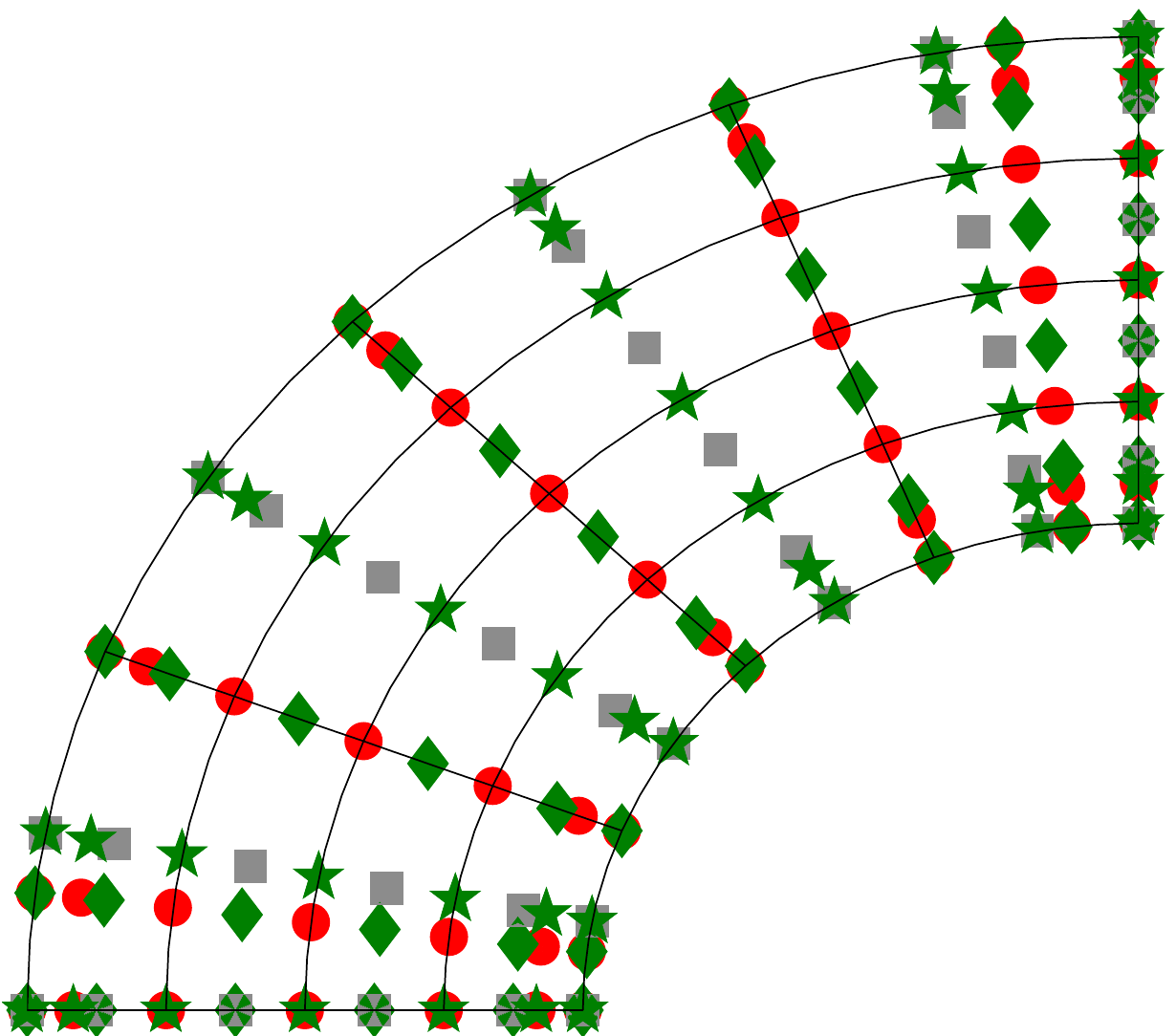}}\hfill
\subfloat[After strong enforcement of no penetration conditions]{\label{sfig:vvp_grid_mapped_bc}\includegraphics[width=.4\textwidth]{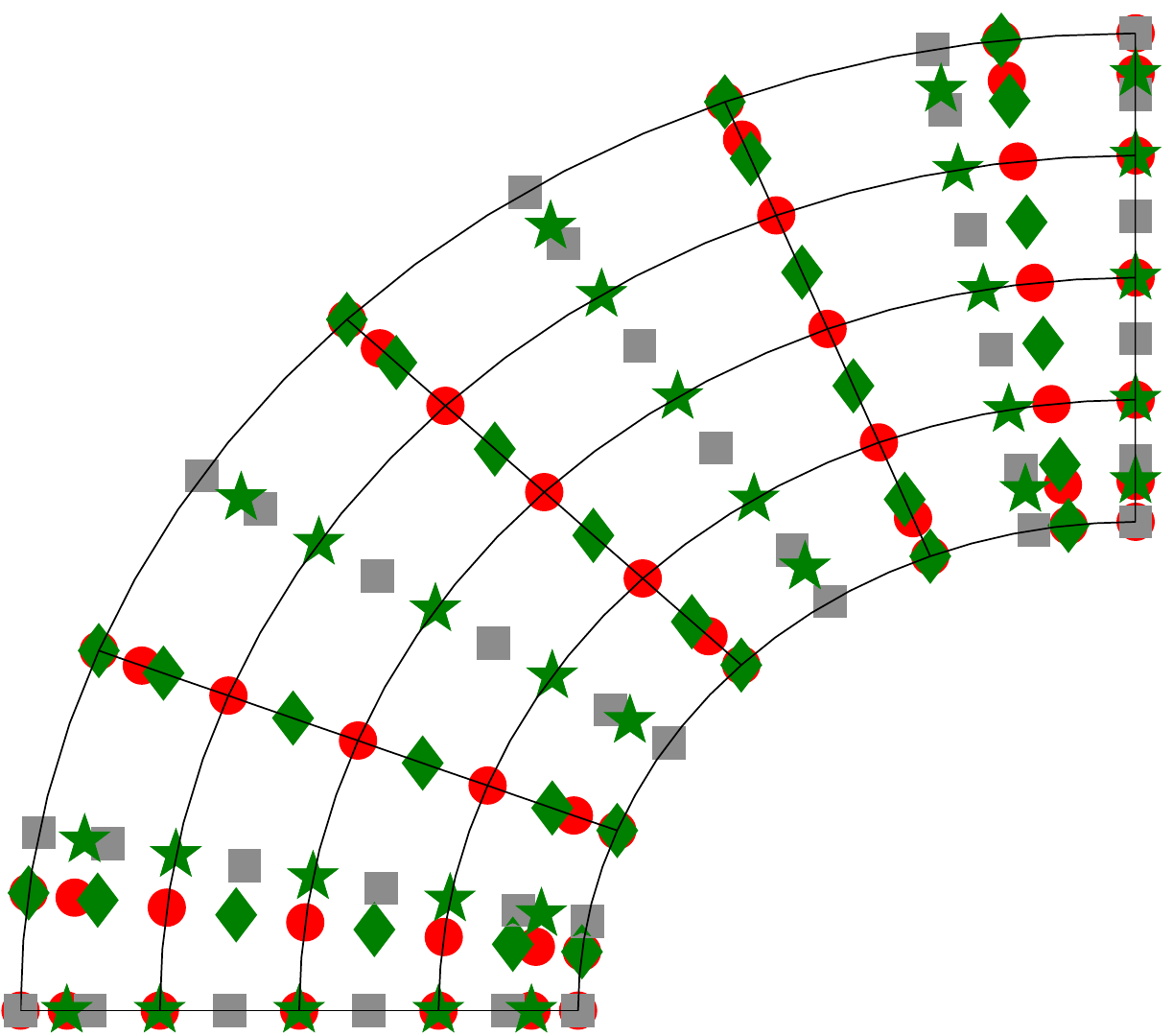}} \\
\caption{Example of collocation grid on a mapped domain for vorticity-velocity-pressure scheme}
\label{fig:vvp_coll_mapped}
\end{figure}

We collocate these equations in the same manner as in the previous sections to solve for the parametric domain variables $\hat{\mathbf{u}}$, $\hat{p}$, and $\hat{{\omega}}$. The collocation points are chosen as the Greville abscissae in the parametric domain, and an example of the resulting points pushed forward into the physical domain is shown in Figure \ref{fig:vvp_coll_mapped}. No penetration boundary conditions are enforced strongly and no slip boundary conditions are enforced weakly with a suitable penalty term. For brevity we omit the full statement of the discrete problem and simply note that it leads to a linear system of equations (as we are focused in this section on Stokes flow).

\section{Numerical Results on Mapped Domains}

In this penultimate section we verify the performance of the vorticity-velocity-pressure collocation scheme on non-square domains. We first consider linear Couette flow to confirm that the expected convergence rates are maintained and then move on to modified lid-driven cavity flows in non-square setups.  

\subsection{Cylindrical Couette Flow}

The first problem posed on a mapped domain that we consider is Couette flow. This models the behavior of a fluid between 2 concentric cylinders, with the outer fixed and the inner rotating at a constant rate. We solve the problem over a quarter circle domain as shown in Figure \ref{fig:vvp_coll_mapped}, enforcing homogeneous Dirichlet boundary conditions on the outer cylindrical wall, zero normal and unit tangential velocity on the inner cylindrical wall, and zero pressure gradient on the horizontal and vertical boundaries. The last Neumann boundary condition is enforced using the Enhanced Collocation approach \cite{deLorenzis_Neumann_contact}. 

The exact velocity field is given in polar coordinates as:

\begin{equation}
    \bar{\textbf{u}} = \left[
\begin{array}{c}
(Ar + B/r)\sin{\theta} \\
(Ar + B/r)\sin{\theta} 
\end{array}
\right],
\end{equation}

\noindent with $A = -\Omega_{in}\frac{\delta^2}{1-\delta^2}$, $B = \Omega_{in}\frac{r_{in}^2}{1-\delta^2}$, $\Omega_{in} = \frac{U}{r_{in}}$, $\delta = \frac{r_{in}}{r_{out}}$, $r_{in} = 1$ is the radius of the inner cylinder, $r_{out} = 2$ is the radius of the outer cylinder, and the velocity of the inner cylinder has magnitude $U = 1$. The exact pressure field is zero everywhere, and the exact vorticity is a constant equal to $2A$. We use a polar mapping to map between the parametric and physical domains:

\begin{equation}
    \textbf{F}(\xi_1, \xi_2) = \left[
\begin{array}{c}
((r_{out} - r_{in}) \xi_2 + r_{in})\sin({2\pi \xi_1}) \\
((r_{out} - r_{in}) \xi_2 + r_{in})\cos({2\pi \xi_1}) 
\end{array}
\right].
\end{equation}

In solving this problem with this mapping one can show analytically that the collocation approximation to the exact solution $\hat{\mathbf{u}}$ is a function of $\hat{y}$ only, the collocation approximation to $\hat{\mathbf{v}}$ is zero, the collocation approximation to $\hat{p}$ is zero, and the collocation approximation to $\hat{\omega}$ is a constant. However, we assemble and solve the full linear system without utilizing this structure.

\begin{figure}
\centering
\subfloat[Velocity $L^2$ error]{\label{sfig:couette_conv_a}\includegraphics[width=.5\textwidth]{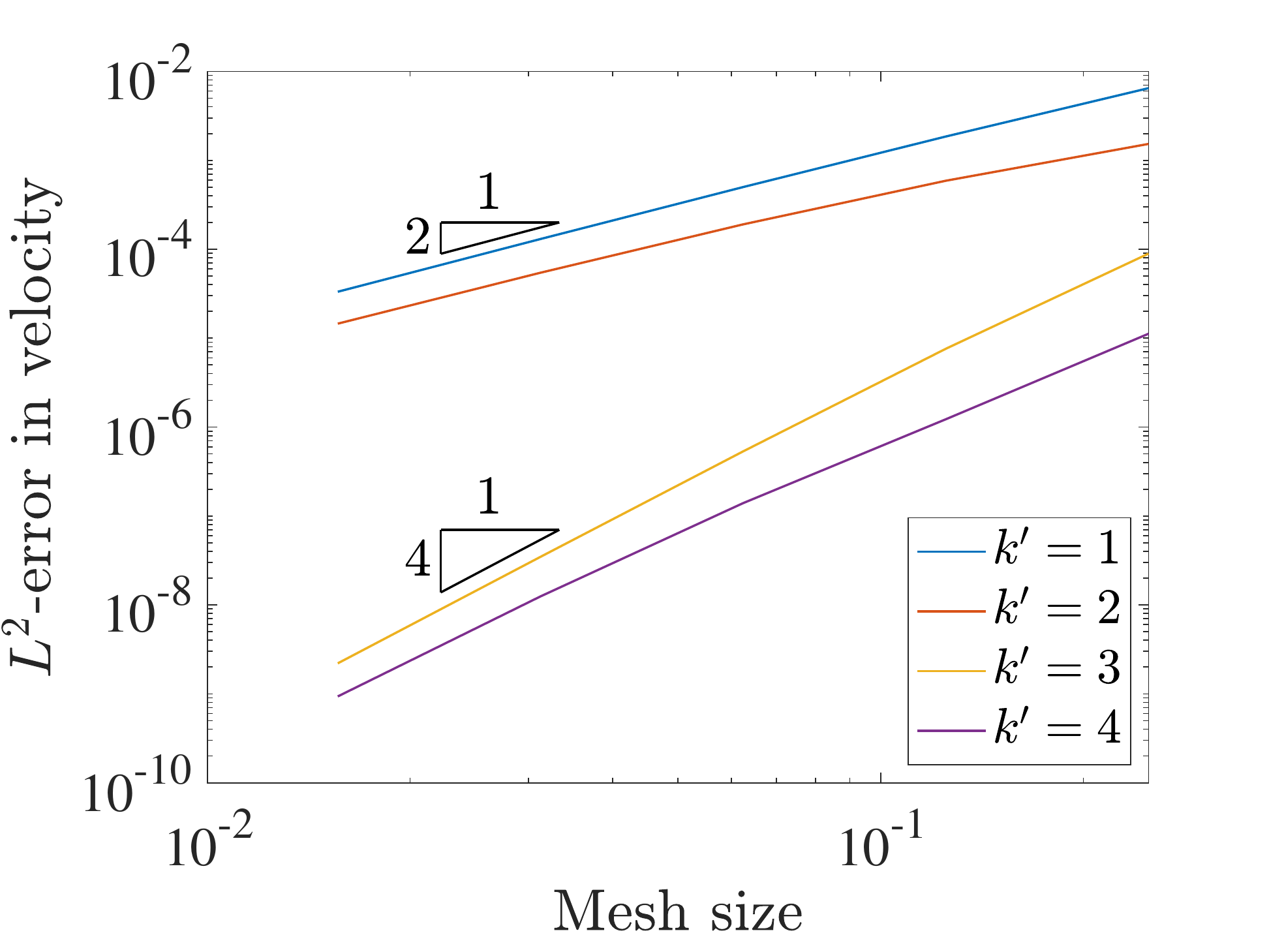}}\hfill
\subfloat[Velocity $H^1$ error]{\label{sfig:couette_conv_b}\includegraphics[width=.5\textwidth]{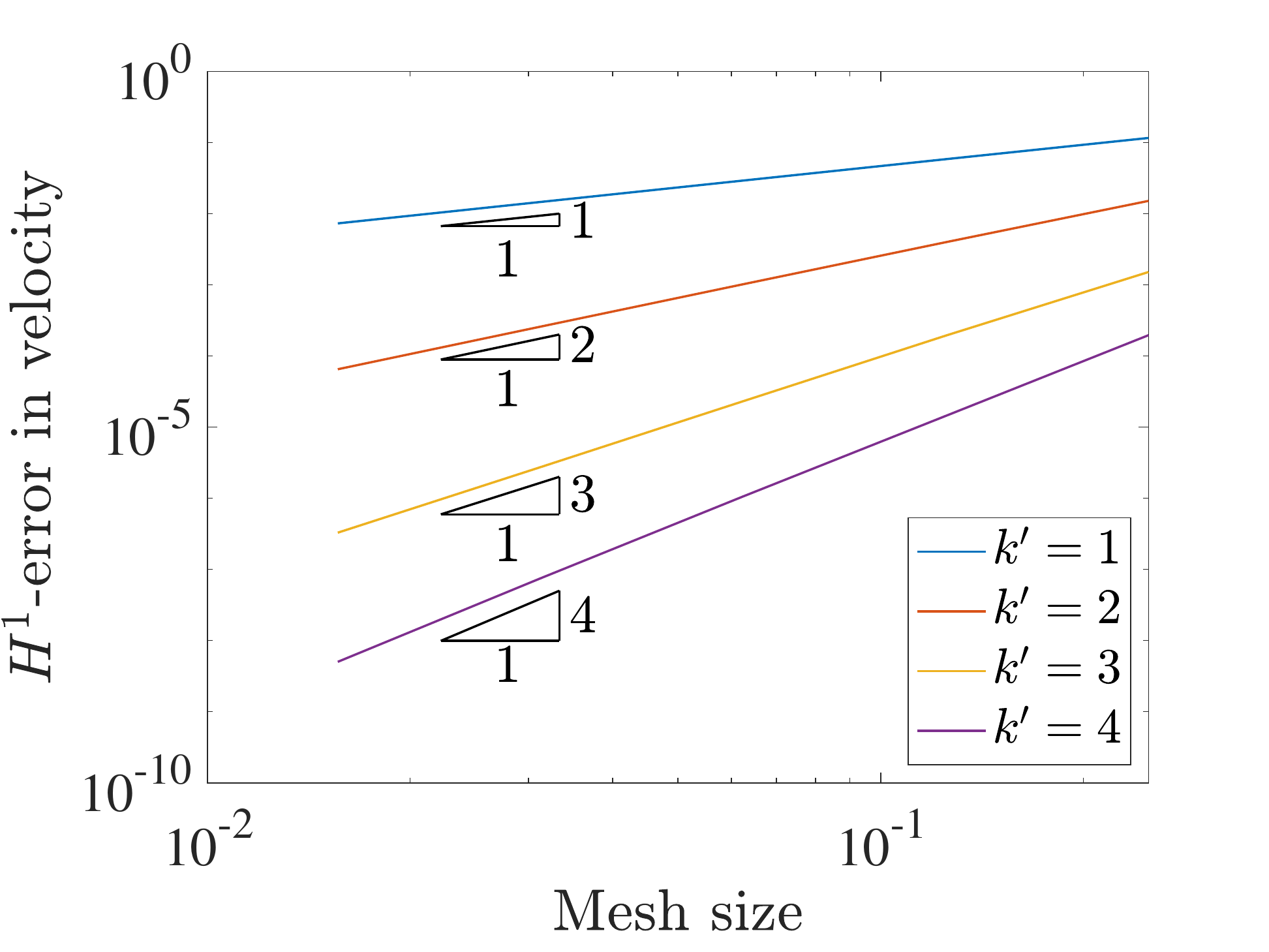}} \\
\subfloat[Pressure $L^2$ error]{\label{sfig:couette_conv_c}\includegraphics[width=.5\textwidth]{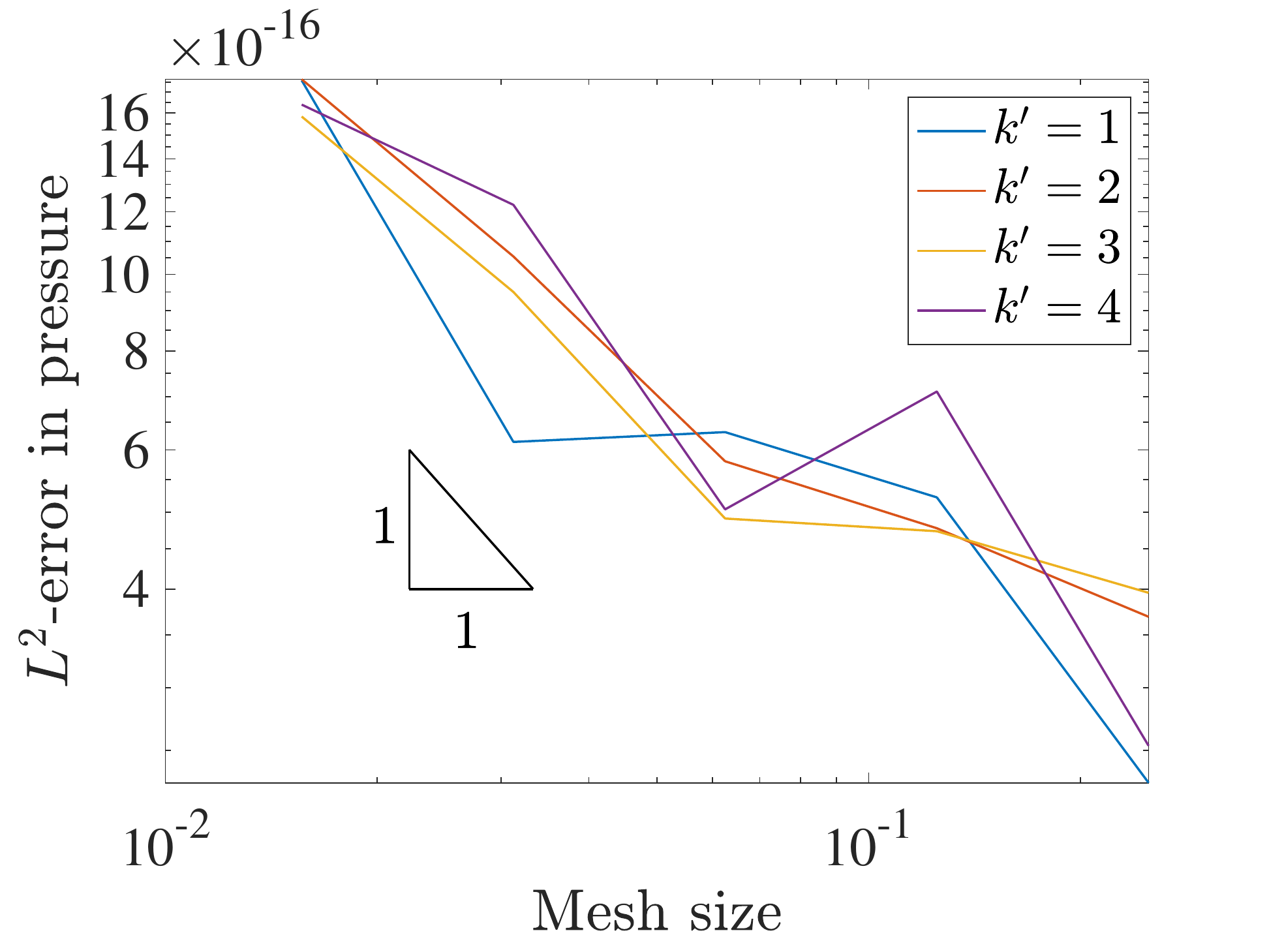}}\hfill
\subfloat[Pressure $H^1$ error]{\label{sfig:couette_conv_d}\includegraphics[width=.5\textwidth]{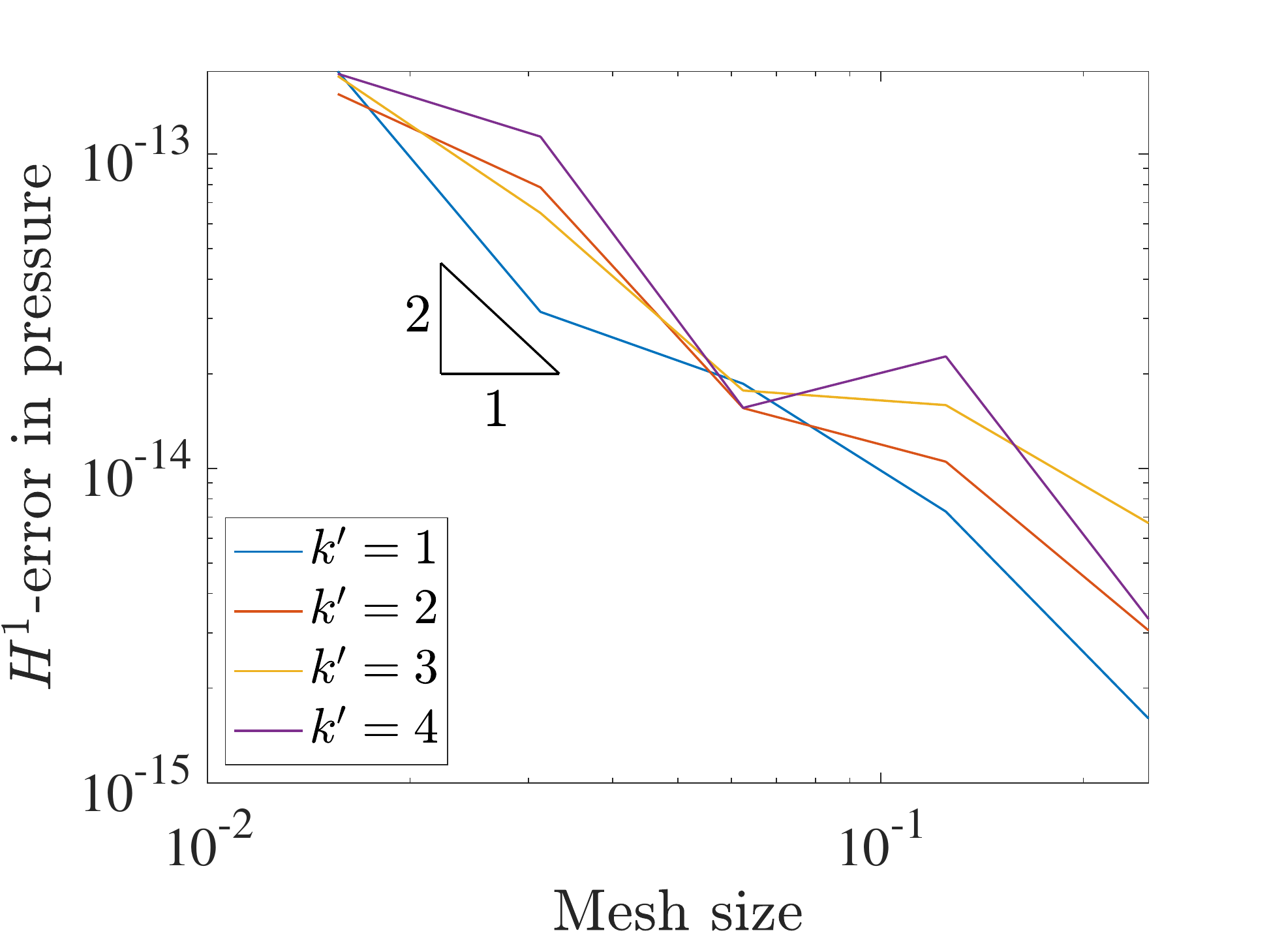}} \\
\subfloat[Vorticity $L^2$ error]{\label{sfig:couette_conv_e}\includegraphics[width=.5\textwidth]{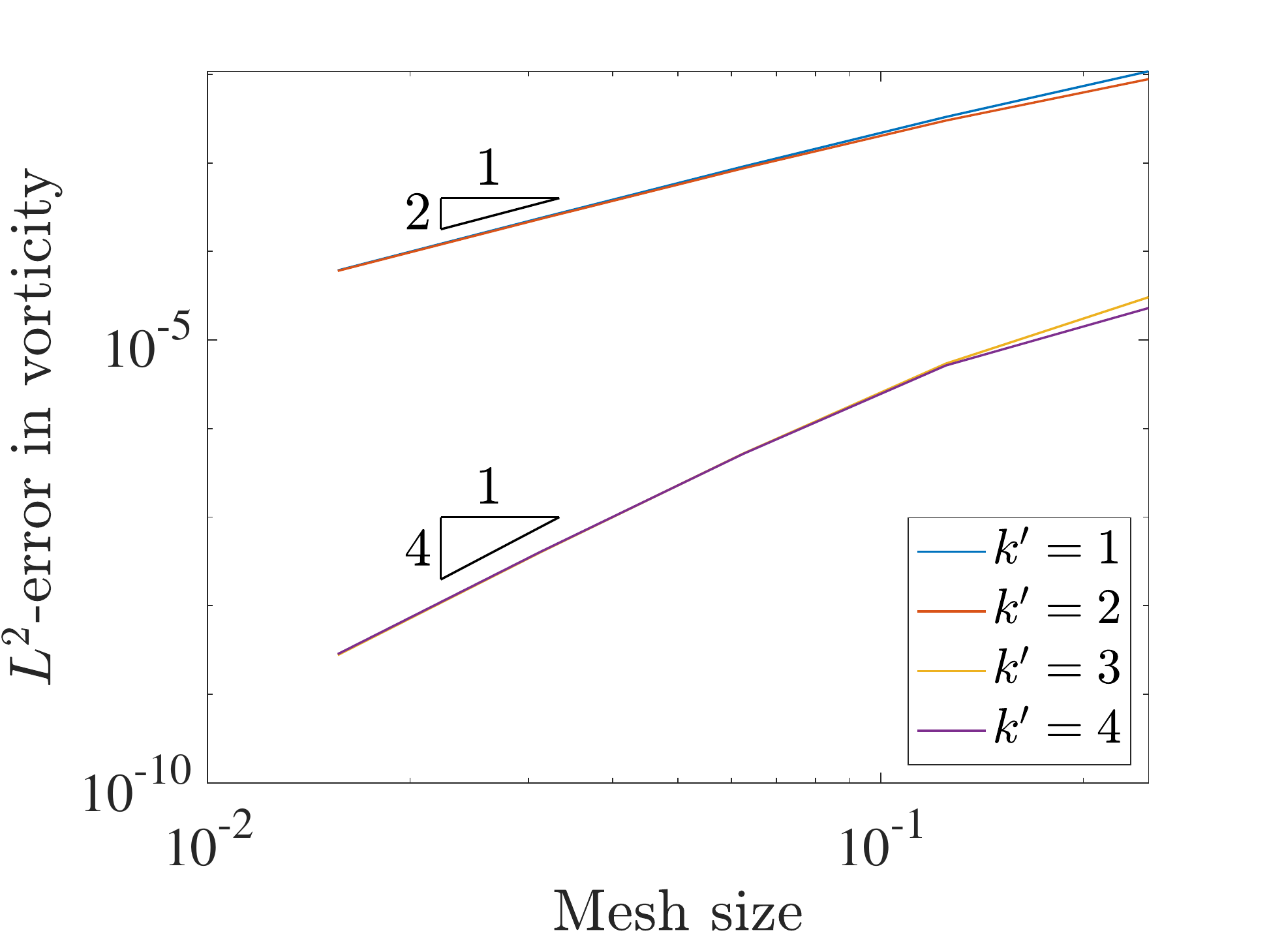}}\hfill
\subfloat[Vorticity $H^1$ error]{\label{sfig:couette_conv_f}\includegraphics[width=.5\textwidth]{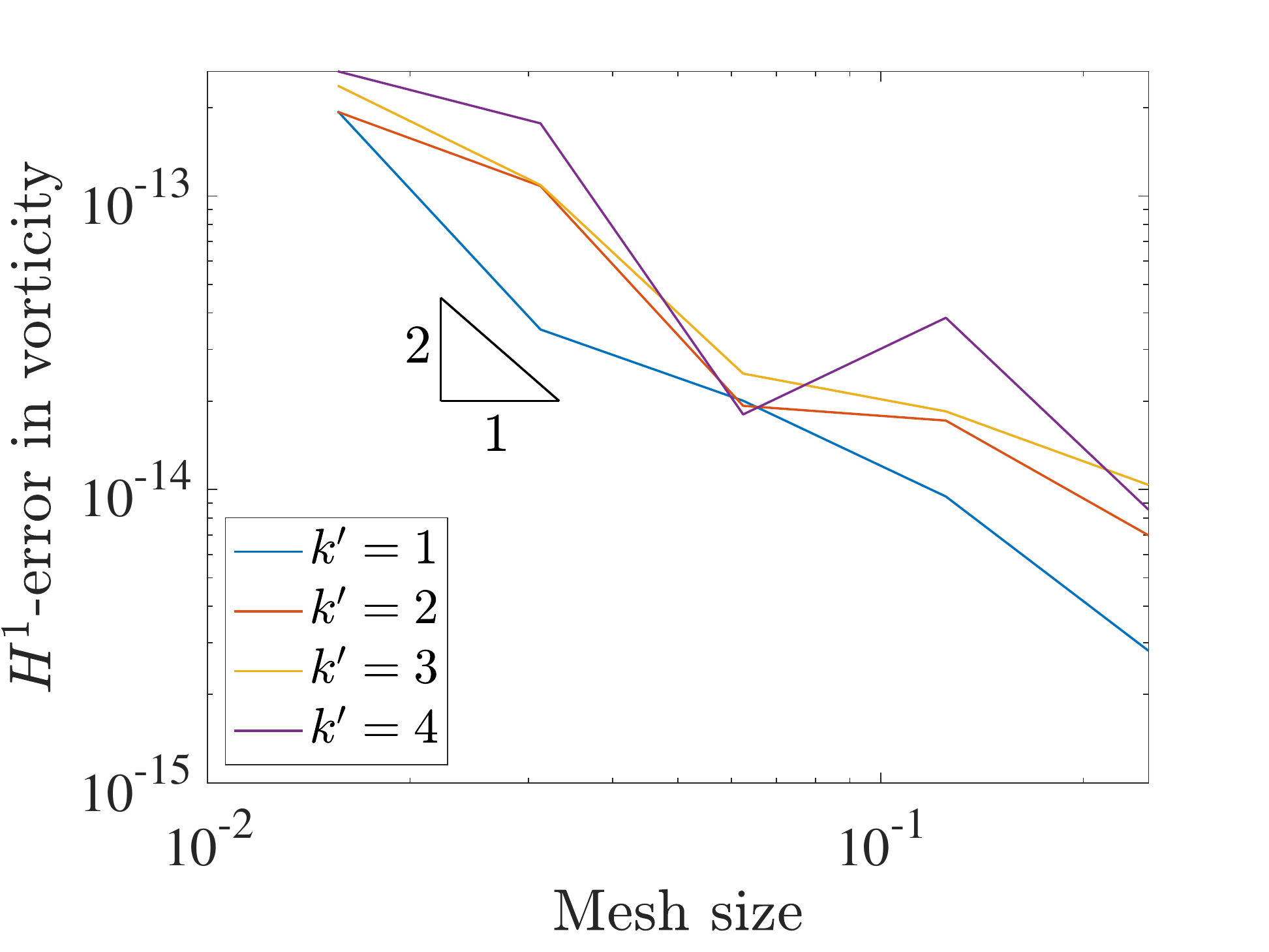}} \\
\caption{Errors in Couette flow solution for vorticity-velocity-pressure formulation}
\label{fig:couette_conv}
\end{figure}

Figure \ref{fig:couette_conv} shows the errors in the solution as a function of resolution. For the $L^2$ norm and $H^1$ semi-norm errors of velocity we recover the same rates are in the square domain setting. The collocation scheme also captures the zero pressure up to finite precision on the coarsest mesh as both the $L^2$ and $H^1$ errors are essentially zero. As the mesh is refined we see this error increase, which we attribute to worsening matrix conditioning and roundoff error effects. We also see the same rates as in square domains for the $L^2$ convergence of vorticity. Note that a constant vorticity is also recovered even on the coarsest mesh, as evidenced by the numerically zero $H^1$ semi-norm error. Like the pressure errors the $H^1$ error grows with mesh refinement, and we believe the explanation is the same. 

\subsection{Lid-Driven Cavity Over Wavy Wall}

Our final numerical test case concerns the Stokes flow in a 2D lid-driven cavity, similar to the square domain examples, but now with a non-flat bottom surface of the cavity. In particular, the mapping from parametric to physical domain is given by

\begin{equation}
    \textbf{F}(\xi_1, \xi_2) = \left[
\begin{array}{c}
 \xi_1\\
 A(B(1-\xi_2)\sin(C \pi \xi_1) + \xi_2)
\end{array}
\right],
\end{equation}

\noindent where $A$, $B$, and $C$ are constants which control the shape of the domain. We use three combinations in this paper, in particular $A = 1$, $B = 0.75$, and $C = 1$ gives a domain with one bump, $A = 0.25$, $B = 0.3$, and $C = 3$ gives a domain with two bumps, and $A = 0.25$, $B = 0.3$, and $C = 5$ gives a domain with three bumps.

\begin{figure}[t]
\centering
\includegraphics[width=\linewidth]{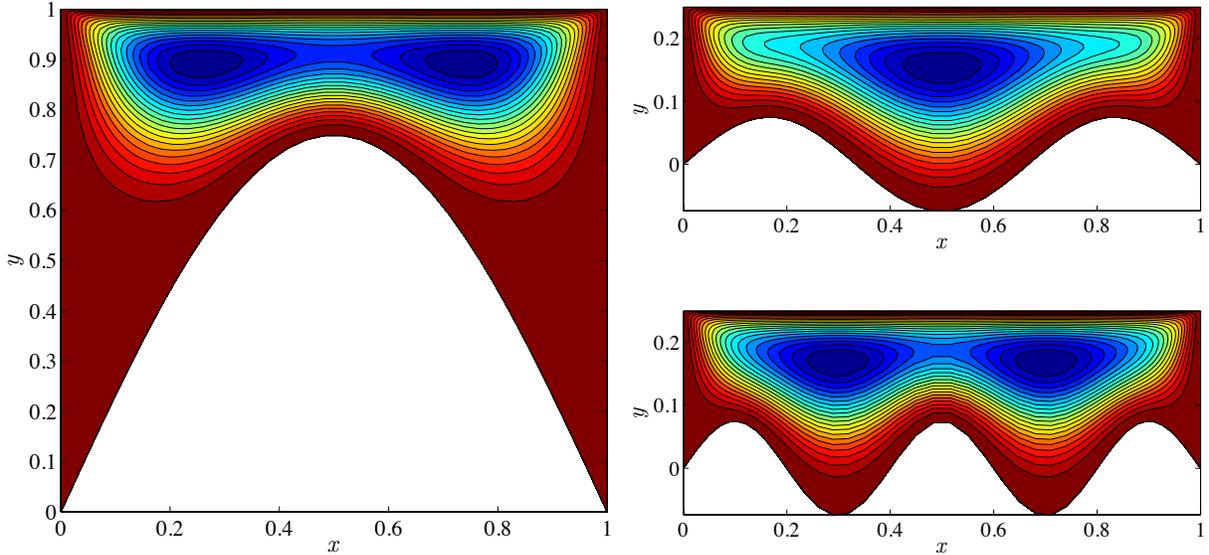}
\caption{Mapped Stokes results for lid-driven cavity with varying numbers of bumps}
\label{fig:bump_results}
\end{figure}

Figure \ref{fig:bump_results} shows the streamfunctions obtained with 64 elements and $k' = 2$. Clearly we are able to recover symmetric fields in all cases which are appropriate for Stokes flow. 

\FloatBarrier

\section{Conclusions}

In this paper, two divergence-conforming collocation methodologies have been presented for solution of the steady, incompressible Navier-Stokes equations using a velocity-pressure formulation and a vorticity-velocity-pressure formulation. By employing B-spline spaces that conform to the de Rham complex, these methods produce velocity fields which are exactly pointwise divergence free. Moreover, by the nature of collocation methods, these methods are much less computationally expensive than traditional Galerkin finite element formulations as no costly numerical integrations are required. By applying the discretizations to benchmark problems in two and three dimensions we have shown that the methods retain a high order of accuracy. Moreover, we have seen that by re-writing the equations in the vorticity-velocity-pressure form many convergence rates are improved compared to those obtained with a velocity-pressure scheme. However, useful properties of the corresponding divergence-conforming B-spline Galerkin method, such as pressure and Reynolds robustness, are not maintained in these collocation schemes. Finally, methods for problems posed in more complicated domains were created by mapping unknowns and equations between the physical and reference domains using structure-preserving transformations.

There are many interesting directions for future work. Collocation schemes that do retain pressure and Reynolds robustness properties would be useful, as would developing a strategy for stabilization of these types of collocation schemes in advection-dominated flow regimes. The schemes proposed in this paper could also be extended to the multi-patch setting to allow for simulations posed on even more complicated domains. The use of locally adaptive splines would also aid in maximizing the ratio of accuracy to cost in which collocation already excels. Finally, while collocation improves upon the cost of numerical integration, unsteady, incompressible Navier-Stokes solution strategies will still likely involve the solution of linear systems during each time step, and thus reducing cost of linear system solution is also very important.

\section*{Acknowledgements}

This material is based upon work supported by the National Science Foundation Graduate Research Fellowship Program under Grant No. DGE-1656518. Any opinions, findings, and conclusions or recommendations expressed in this material are those of the authors and do not necessarily reflect the views of the National Science Foundation.

\bibliographystyle{elsarticle-num}
\bibliography{ref}

\end{document}